\documentclass[nonblindrev]{informs3} 

\DoubleSpacedXI 

\usepackage{hyperref}
\usepackage{endnotes}
\let\footnote=\endnote
\let\enotesize=\normalsize
\def\notesname{Endnotes}%
\def\makeenmark{$^{\theenmark}$}
\def\enoteformat{\rightskip0pt\leftskip0pt\parindent=1.75em
  \leavevmode\llap{\theenmark.\enskip}}

\usepackage[caption=false]{subfig}
\usepackage{natbib}
 \bibpunct[, ]{(}{)}{,}{a}{}{,}%
 \def\bibfont{\small}%
\TheoremsNumberedThrough     
\ECRepeatTheorems

\EquationsNumberedThrough    

\begin{document}


\RUNAUTHOR{Esteban-P\'erez and Morales}

\RUNTITLE{Distributionally robust conditional stochastic programs based on trimmings}

\TITLE{Distributionally robust  stochastic programs with side information based on trimmings -- Extended version}

\ARTICLEAUTHORS{%
\AUTHOR{Adri\'an Esteban-P\'erez }
\AFF{Department of Applied Mathematics,  University of M\'alaga, M\'alaga, 29071, Spain \EMAIL{adrianesteban@uma.es}} 
\AUTHOR{Juan M. Morales}
\AFF{Department of Applied Mathematics,  University of M\'alaga, M\'alaga, 29071, Spain, \EMAIL{juan.morales@uma.es}}
} 

\ABSTRACT{%
We consider stochastic programs conditional on some covariate information, where the only knowledge of the possible relationship between the uncertain parameters and the covariates is reduced to a finite data sample of their joint distribution. By exploiting the close link between the notion of \emph{trimmings} of a probability measure and the \emph{partial} mass transportation problem, we construct a data-driven Distributionally Robust Optimization (DRO) framework to hedge the decision against the intrinsic error in the process of inferring conditional information from limited joint  data. We show that our approach is computationally as tractable as the standard (without side information) Wasserstein-metric-based DRO and enjoys performance guarantees. Furthermore, our DRO framework can be conveniently used to address data-driven decision-making problems under contaminated samples  and naturally produces distributionally robust versions of some local nonparametric predictive methods, such as  Nadaraya-Watson kernel regression and $K$-nearest neighbors, which are often used in the context of conditional stochastic optimization. Finally, the theoretical results are illustrated using a  single-item newsvendor problem and a portfolio allocation problem with side information.
}%



\maketitle

%


\section{Introduction}

Today's decision makers not only collect observations of the uncertainties directly affecting their decision-making processes, but also gather data about measurable exogenous variables that may have some predictive power on those uncertainties \citep{Ban2019}. In Statistics, Operations Research and Machine Learning, these variables are often referred to as \emph{covariates}, \emph{explanatory variables}, \emph{side information} or \emph{features} \citep{PangHo2019}.

In the framework of Optimization Under Uncertainty, the side information acts by changing the probability measure of the uncertainties. In fact, if the joint distribution of the features and the uncertainties were known, this measure change would correspond to conditioning that distribution on the side information given. Unfortunately, in practice, the decision maker only has an incomplete picture of such a joint distribution in the form of a finite data sample. The development of optimization methods capable of exploiting the side information to make improved decisions, in a context of limited knowledge of its explanatory power on the uncertainties, defines the ultimate purpose of the so-called \emph{Prescriptive Stochastic Programming} or \emph{Conditional Stochastic Optimization} paradigm. This paradigm has recently become very popular in the technical literature, see, for instance, \cite{Ban2019,Bertsimas2019b,PangHo2019} and references therein. More specifically, a data-driven approach to address the newsvendor problem, whereby the decision is explicitly modeled as a parametric function of the features, is proposed in \cite{Ban2019}. This approach thus seeks to optimize said function.
In contrast,
\cite{Bertsimas2019b} formulate and formalize the problem of minimizing the conditional expectation cost given the side information, and develop various schemes based on  machine learning methods (typically used for regression and prediction) to get data-driven solutions. Their approach is  \emph{non-parametric} in the sense that the optimal decision is not constrained to be a member of a certain family of the features' functions. The inspiring work of \cite{Bertsimas2019b} has been subject to further study and improvement in two principal directions, namely, the design of efficient algorithms to trim down the computational burden of the optimization \citep{Diao2020} and the development of strategies to reduce the variance and bias of the decision obtained and its associated cost (the pairing of both interpreted as a statistical estimator).
In the latter case, we can cite the work of
\cite{bertsimas2017bootstrap}, where they leverage ideas from bootstrapping and machine learning to confer robustness on the decision and acquire   asymptotic  performance guarantees.
Similarly, \cite{BertsimasMCCORD2018} and \cite{PangHo2019} propose regularization procedures to reduce the variance of the data-driven solution to the conditional expectation cost minimization problem, which is formalized and studied in \cite{Bertsimas2019b}. A scheme to robustify the data-driven methods introduced in this work is also proposed in \cite{Bertsimas2019} for dynamic decision-making.

 A different, but related thrust of research focuses on developing methods to construct predictions specifically tailored to the optimization problem that is to be solved and where those predictions are then used as input information. Essentially, the predictions are intended to yield decisions with a low disappointment or regret. This framework is known in the literature as (smart) \emph{Predict-then-Optimize}, see, e.g., \cite{Balghiti2019,Donti2017,Elmachtoub2021,munoz2020bilevel}, and references therein.

Our research, in contrast, builds upon Distributionally Robust Optimization (DRO), which is a powerful modeling paradigm to protect the task of decision-making against the ambiguity of the underlying probability distribution of the uncertainty \citep{Rahimian2019}.
 Nevertheless, the technical literature on the use of DRO to address Prescriptive or Conditional Stochastic Programming problems
 is still relatively scarce. We highlight
\cite{Bertsimas2019,Chensim2020,hanasusanto2013robust,kannan2021heteroscedasticity,kannan2020droresiduals,Nguyen2020,Nguyen2021}\footnote{The preprints \cite{kannan2020droresiduals,kannan2021heteroscedasticity,Nguyen2020,Nguyen2021} became  available online while this paper was under review in this journal.}, with \cite{Nguyen2021} being a generalization of \cite{Nguyen2020}. In  \cite{Chensim2020}, they resort to a scenario-dependent ambiguity set to exploit feature information in a DRO framework. However, their objective is to minimize a joint expectation and consequently, their approach cannot directly handle the Conditional Stochastic Optimization setting we consider here.   \cite{hanasusanto2013robust} deal with a stochastic control problem with time-dependent data. They extend the idea of \cite{hannah2010nonparametric} to a fully dynamic setting and robustify the control policy against the worst-case weight vector that is within a certain $\chi^2$-distance from the one originally given by the Nadaraya-Watson estimator. In the case of \cite{Bertsimas2019}, the authors propose using the conditional empirical distribution given by a local predictive method as the center of the Wasserstein ball that characterizes the DRO approach in \cite{MohajerinEsfahani2018}. This proposal, nonetheless, fails to explicitly account for the inference error associated with the local estimation. In \cite{kannan2021heteroscedasticity,kannan2020droresiduals}, the authors develop a two-step procedure whereby a regression model between the uncertainty and the features is first estimated and then a distributionally robust decision-making problem is formulated, considering a Wasserstein ball around the empirical distribution of the residuals. Finally, the authors in \cite{Nguyen2021} also consider a Wasserstein-ball ambiguity set as in \cite{Bertsimas2019,kannan2021heteroscedasticity,kannan2020droresiduals}, but centered at the empirical distribution of the joint data sample of the uncertainty and the features. In addition, they further constrain the ambiguity set by imposing that the worst-case distribution assigns some probability mass to the support of the uncertainty conditional on the values taken on by the features.

Against this background, our main contributions  are:
\begin{enumerate}
    \item \emph{Modeling power:} We develop a general framework
    to handle
   prescriptive stochastic programs within the DRO paradigm. Our DRO framework is based on  a new class of ambiguity sets that exploit the close and convenient connection between trimmings and the partial mass problem to immunize the decision against the error incurred in the process of inferring \emph{conditional} information from \emph{joint} (limited) data. 
  %
    %
    We also show that our approach serves as a natural framework for the application of DRO in data-driven decision-making under contaminated samples
   and naturally produces distributionally robust versions of some local nonparametric predictive methods such as Nadaraya-Watson kernel regression and $K$-nearest neighbors, which are used in the context of conditional stochastic optimization  \citep{bertsimas2017bootstrap,BertsimasMCCORD2018,Bertsimas2019,PangHo2019}.

    \item \emph{Computational tractability:} Our framework is as complex as the  Wasserstein-metric-based DRO approach proposed in \cite{MohajerinEsfahani2018} without side information. Therefore, we extend the mass-transportation approach to the realm of Conditional Stochastic Optimization while preserving its appealing tractability properties.


    \item \emph{Theoretical results and performance guarantees:}
  Leveraging theory from  probability trimmings and optimal transport,  we show that our DRO model enjoys a finite sample guarantee and is asymptotically consistent.


    \item \emph{Numerical results:} We evaluate our DRO approach on
    %
   the single-item newsvendor problem and the portfolio allocation problem, 
    and compare it with the KNN method described in \cite{Bertsimas2019b}, the robustified KNN proposed in \cite{Bertsimas2019}, and
    a KNN followed by the standard Wassertein-distance-based DRO model introduced in~\cite{MohajerinEsfahani2018}, as suggested in~\cite{Bertsimas2019} too. Unlike all these approaches, ours explicitly accounts for the cost impact of the potential error made when inferring conditional information from a joint sample of the uncertainty and the covariates. To this end, we minimize the worst-case cost over a Wasserstein ball of probability measures with \emph{an ambiguous center}.  

\end{enumerate}

The rest of the paper is organized as follows.   In Section \ref{ddro_section}, we  formulate our DRO framework to address decision-making problems under uncertainty in the presence of side information and  show that it is as tractable as the standard Wasserstein-metric-based DRO approach developed in \cite{MohajerinEsfahani2018}.
 In  Section \ref{case_alpha_pos_uni}, we
    deal with the case in which the side information corresponds to an event of known and positive probability and discuss its application to data-driven decision-making under contaminated samples.
    The situation in which the probability of such an event is positive, but unknown, is treated in
    Section~\ref{unknown-alpha}. Section \ref{case_alpha_0}
    elaborates on the case in which the side information reduces to a specific realization of the feature vector, more precisely, the instance where the side information represents an event of zero probability. Section~\ref{numerics} provides results from numerical experiments and, finally, Section \ref{conclusions} concludes the paper.

\textbf{Notation}.
We  use $\overline{\mathbb{R}}$ to represent
the extended real line, and adopt the
conventions of its associated arithmetic. Moreover,  $\mathbb{R}_+$ stands for the set of non-negative real numbers.
We employ lower-case bold face
letters to represent vectors. The inner
product of two vectors $ \mathbf{u},
\mathbf{v} $  is denoted as $\langle
\mathbf{u}, \mathbf{v}\rangle = \mathbf{u}^T
\mathbf{v}$ and by $\left\| \mathbf{u}\right\|$ we denote the norm of the vector $\mathbf{u}$.
 For a set $A$, the indicator function $\mathbb{I}_A(\mathbf{a})$ is defined through $\mathbb{I}_A(\mathbf{a})=1$
if $\mathbf{a}\in A$; $=0$ otherwise.
 The Lebesgue measure in $\mathbb{R}^d$ is denoted as $\lambda^d$.
 We use the symbol $\delta_{\boldsymbol{\xi}}$ to represent the Dirac distribution supported on $\boldsymbol{\xi}$. Additionally, we reserve the symbol ``$\;\widehat{}\;$'' for objects which are
  dependent  on the sample data. The $K$-fold product of a distribution $\mathbb{Q}$ will be denoted as $\mathbb{Q}^K$.  Finally, the symbols $\mathbb{E}$ and $\mathbb{P}$ denote, respectively, ``expectation'' and ``probability'' (the context will give us the measure under which that expectation or probability is taken).

\section{Data-driven distributionally robust
optimization  with side information}\label{ddro_section}

In this paper, we propose a general framework for data-driven distributionally robust optimization with side information that relies on two related tools, namely, the \emph{optimal mass transport theory} and the concept of  \emph{trimming of a probability measure}. Next, we introduce some preliminaries that  help  motivate our proposal. All the proofs that are missing in the main text are compiled in the Appendix.

\subsection{Preliminaries and motivation}\label{preliminaries}
Let $\mathbf{x}\in X \subseteq \mathbb{R}^{d_{\mathbf{x}}}$ be the decision variable vector and $\mathbf{y}$, with support set $\Xi_{\mathbf{y}} \subseteq \mathbb{R}^{d_{\mathbf{y}}}$, the random vector that models the uncertainty affecting the value of the decision. Let $\mathbf{z}$, with support set $\Xi_{\mathbf{z}} \subseteq \mathbb{R}^{d_{\mathbf{z}}}$, be the (random) feature vector and denote the objective function to be minimized as $f(\mathbf{x},\boldsymbol{\xi})$, where $\boldsymbol{\xi}:=(\mathbf{z},\mathbf{y})$.



 Given a new piece of information in the form of the event $\boldsymbol{\xi} \in \widetilde{\Xi}$, the decision maker seeks to compute the optimal decision that minimizes the (true) conditional
expected cost:
\begin{equation}\label{conditional_expectation_problem}
    J^*:=\inf_{\mathbf{x}\in X} \mathbb{E}_{\mathbb{Q}} \left[ f(\mathbf{x},\boldsymbol{\xi}) \; \mid \; \boldsymbol{\xi} \in \widetilde{\Xi} \right] = \inf_{\mathbf{x}\in X} \mathbb{E}_{\mathbb{Q}_{\widetilde{\Xi}}} \left[ f(\mathbf{x},\boldsymbol{\xi})\right]
\end{equation}
where $\mathbb{Q}$ is the true joint distribution of $\boldsymbol{\xi}:=(\mathbf{z},\mathbf{y})$  with support set $\Xi \subseteq \mathbb{R}^{d_{\mathbf{z}}+d_{\mathbf{y}}}$ and $\mathbb{Q}_{\widetilde{\Xi}}$ is the associated true distribution of $\boldsymbol{\xi}$ \emph{conditional on} $\boldsymbol{\xi} \in \widetilde{\Xi}$. Hence, we implicitly assume that $\mathbb{Q}_{\widetilde{\Xi}}$ is a regular conditional distribution and that the conditional expectation~\eqref{conditional_expectation_problem} is well defined.

An example of $\widetilde{\Xi}$ would be $\widetilde{\Xi}:=\{ \boldsymbol{\xi}=(\mathbf{z},\mathbf{y}) \in \Xi \; : \; \mathbf{z}\in \mathcal{Z}\}$, with $\mathcal{Z}\subseteq \Xi_{\mathbf{z}} $ being an uncertainty set built from the information on the features. We note that this definition includes the case in which $\mathcal{Z}$ reduces to a singleton $\mathbf{z}^*$ representing a particular realization of the features.

Unfortunately, when it comes to solving problem~\eqref{conditional_expectation_problem}, neither the true distribution $\mathbb{Q}$ nor ---even less so--- the conditional one $\mathbb{Q}_{\widetilde{\Xi}}$ are generally known to the decision maker. Actually, the decision maker typically counts only on a data sample consisting of $N$  observations $\widehat{\boldsymbol{\xi}}_i:=(\widehat{\mathbf{z}}_{i},\widehat{\mathbf{y}}_{i})$ for
$i=1,\ldots, N$, which we assume are i.i.d. Therefore, the solution to problem~\eqref{conditional_expectation_problem} \emph{per se} is, in practice, out of reach and the best the decision maker can do is to approximate the solution to~\eqref{conditional_expectation_problem} with some (probabilistic) performance guarantees. Within this context, \emph{Distributionally Robust Optimization} (DRO) emerges as a powerful modeling framework to achieve that goal. In brief, the DRO approach aims to find a decision $\mathbf{x}\in X$ that is \emph{robust} against all \emph{conditional} probability distributions that are somehow \emph{plausible} given the information at the decision maker's disposal. This is mathematically stated as follows:
\begin{align}\label{conditional_model1}
\inf_{\mathbf{x}\in X} \sup_{ Q_{\widetilde{\Xi}} \in \widehat{\mathcal{U}}_N}& \;  \mathbb{E}_{Q_{\widetilde{\Xi}}} \left[f(\mathbf{x},\boldsymbol{\xi}) \right]
\end{align}
where $\widehat{\mathcal{U}}_N$ is a so-called \emph{ambiguity set} that contains all those plausible conditional distributions. This ambiguity set must be built from the available information on $\boldsymbol{\xi}$, which, in our case, comprises the $N$ observations $\{\widehat{\boldsymbol{\xi}}_i\}_{i=1}^{N}$. The subscript $N$ in $\widehat{\mathcal{U}}_N$ is intended to underline this issue. Furthermore, the condition $Q_{\widetilde{\Xi}}(\widetilde{\Xi}) = 1$ for all $Q_{\widetilde{\Xi}} \in \widehat{\mathcal{U}}_N$ is implicit in the construction of that set.
In our setup, however, problem~\eqref{conditional_model1} poses a major challenge, which has to do with the fact that the observations $\{\widehat{\boldsymbol{\xi}}_i\}_{i=1}^{N}$
pertain to the true \emph{joint} distribution $\mathbb{Q}$, and \emph{not} to the conditional one $\mathbb{Q}_{\widetilde{\Xi}}$. 
Consequently,
we need to build an ambiguity set $\widehat{\mathcal{U}}_N$ for the plausible \emph{conditional} distributions from the limited joint information on $\mathbb{Q}$ provided
by the data $\{\widehat{\boldsymbol{\xi}}_i\}_{i=1}^{N}$.

At this point, we should note that there are several approaches in the technical literature to handle
the conditional stochastic optimization problem   \eqref{conditional_expectation_problem} for the particular case in which $\widetilde{\Xi}$ is defined as  $\widetilde{\Xi}:=\{ \boldsymbol{\xi}=(\mathbf{z},\mathbf{y}) \in \Xi \; : \; \mathbf{z}= \mathbf{z^*}\}$. For example,  \cite{Bertsimas2019b} approximate \eqref{conditional_expectation_problem} by the following conditional estimate
  \begin{equation}\label{bertsimaskallus_problem}
   \inf_{\mathbf{x}\in X}
   \sum_{i=1}^N w_N^i(\mathbf{z}^*)
   f(\mathbf{x},(\mathbf{z}^*,\widehat{\mathbf{y}}_{i}) )
\end{equation}
where $w_N^i(\mathbf{z}^*)$ is a weight function that can be given by various non-parametric machine learning methods such as $K$-nearest neighbors, kernel regression, CART, and random forests.
Formulation~\eqref{bertsimaskallus_problem} can be naturally interpreted as a (conditional) Sample-Average-Approximation (SAA) of problem \eqref{conditional_expectation_problem}.

\cite{BertsimasMCCORD2018} extend the work by~\cite{Bertsimas2019b} to accommodate the setting in which the outcome of the uncertainty $\mathbf{y}$ may be contingent on the taken decision $\mathbf{x}$. For this purpose, they work with an enriched data set comprising observations of the uncertainty $\mathbf{y}$, the decision $\mathbf{x}$ and the covariates $\mathbf{z}$, and allow the weights in~\eqref{bertsimaskallus_problem} to depend on $\mathbf{x}$ too. Besides, they add terms to the objective function of~\eqref{bertsimaskallus_problem} to penalize estimates of its variance and bias.
The case in which the weight function~\eqref{bertsimaskallus_problem} is given by the Nadaraya-Watson (NW) kernel regression estimator is considered in  \cite{hannah2010nonparametric} and \cite{PangHo2019}. In \cite{PangHo2019}, in addition, they leverage techniques from moderate deviations theory to design a regularization scheme that reduces the optimistic bias of the NW approximation and to provide insight into its out-of-sample performance. The work in \cite{bertsimas2017bootstrap}   focuses on conditional estimators~\eqref{bertsimaskallus_problem} where the weights are provided by the NW or KNN method. They use DRO, based on the relative entropy distance for discrete  distributions to get decisions from~\eqref{bertsimaskallus_problem} that perform well on a large portion of resamples \emph{bootstraped} from the empirical distribution of the available data set.

Finally, \cite{Bertsimas2019} provide a robustified version of the conditional estimator~\eqref{bertsimaskallus_problem}, which takes the following form
  \begin{equation}\label{bertsimas_sturt}
   \inf_{\mathbf{x}\in X}
   \sum_{i=1}^N w_N^i(\mathbf{z}^*) \sup_{\mathbf{y} \in \mathcal{U}^i_N }
      \left[f(\mathbf{x},(\mathbf{z}^*,\mathbf{y}) ) \right]
\end{equation}
    where $\mathcal{U}^i_N:= \{ \mathbf{y} \in \Xi_{\mathbf{y}}  : \left\| \mathbf{y}-\widehat{\mathbf{y}}_i \right\|_p \leqslant \varepsilon_N \}$.
This  problem can be seen as a robust SAA method capable of exploiting side information  and has also been used in~\cite{bertsimas2018data,BertsimasShtern2021}.

In our case, however, we  follow a different path to address the conditional stochastic optimization problem~\eqref{conditional_expectation_problem} by way of~\eqref{conditional_model1}. More precisely, we  leverage the notion of \emph{trimmings of a distribution} and the related theory of \emph{partial mass transportation}.

\subsection{The Partial Mass Transportation
Problem and Trimmings}\label{otp_section_trimmings}
This section introduces some concepts about trimmings and the partial mass transportation problem that  help us construct the ambiguity set $\widehat{\mathcal{U}}_N$ in \eqref{conditional_model1} from the sample data $\{\widehat{\boldsymbol{\xi}}_i\}_{i=1}^{N}$. For simplicity, we restrict ourselves to probability measures defined in $\mathbb{R}^d$.

If $\mathbb{Q}(\widetilde{\Xi}) = \alpha > 0$  (our analysis, though, will also cover the case $\alpha = 0$ later in Section~\ref{case_alpha_0})), problem~\eqref{conditional_expectation_problem} can be recast as
\begin{equation}\label{conditional_expectation_problem_transf}
    J^*:=\inf_{\mathbf{x}\in X} \frac{1}{\alpha}\mathbb{E}_{\mathbb{Q}} \left[f(\mathbf{x},\boldsymbol{\xi})\mathbb{I}_{\widetilde{\Xi}}(\boldsymbol{\xi})\right]
\end{equation}
which only requires that
$\mathbb{E}_{\mathbb{Q}}\left[ |f(\mathbf{x},\boldsymbol{\xi})\mathbb{I}_{\widetilde{\Xi}}(\boldsymbol{\xi})| \right]<\infty$ for all $\mathbf{x}\in X$, see   \cite[Eq. 6.2]{Gray2009}.

Now we introduce the notion of a \emph{trimming} of a distribution, which is at the core of our proposed DRO framework.

\begin{definition}[$(1-\alpha)$-trimmings, Definition 1.1 from \cite{DelBarrio2013a}] \label{def_trimming}
Given $0 \leqslant \alpha \leqslant 1$ and  probability measures $P,Q \in \mathbb{R}^d$, we say that $Q$ is an \emph{$(1-\alpha)$-trimming} of $P$ if $Q$ is absolutely continuous with
respect to $P$, and the Radon-Nikodym derivative satisfies $\frac{dQ}{dP} \leqslant \frac{1}{\alpha}$. The \emph{set of
all $(1-\alpha)$-trimmings} (or trimming set of level $1-\alpha$) of $P$  will be denoted by $\mathcal{R}_{1-\alpha}(P)$.
\end{definition}

As extreme cases, we have that for $\alpha = 1,\; \mathcal{R}_0(P)$ is just $P$, while, for $\alpha = 0$, $\mathcal{R}_1(P)$ is the set of
all probability measures absolutely continuous with respect to $P$.
Given a probability $P$ on $\mathbb{R}^d$, if $\alpha_1 \leqslant \alpha_2$, then  $\mathcal{R}_{1-\alpha_2}(P) \subset \mathcal{R}_{1-\alpha_1}(P) $. Especially useful is the fact that a trimming set is a convex set, which is, besides,   compact under the topology of weak convergence.
We refer the reader to \cite[Proposition 2.7]{Alvarez-Esteban2011}  for other interesting
properties about the set $\mathcal{R}_{1-\alpha}(P)$.




%

Consider now the following minimization problem:
\begin{equation}\label{partial_mass_problem}
    \inf_{Q \in \mathcal{R}_{1-\alpha}(P)} D(Q,R)
\end{equation}
where $D$ is a probability metric.

Problem~\eqref{partial_mass_problem} is known as the $(D,1-\alpha)-$\emph{partial (or incomplete) mass problem} \citep{DelBarrio2013a}. While there is a variety of probability metrics we could choose from to play the role of $D$ in~\eqref{partial_mass_problem}, here we  work with the space $\mathcal{P}_p(\mathbb{R}^d)$ of probability distributions supported on $\mathbb{R}^d$ with finite $p$-th moment and restrict ourselves to the $p-$Wasserstein metric, $\mathcal{W}_p$, for its tractability and theoretical advantages. In such a case (i.e., when $D=\mathcal{W}_p$), problem~\eqref{partial_mass_problem} is referred to as a partial mass \emph{transportation} problem and interpolates between the classical optimal mass transportation problem (when
$\alpha=1$) and the random quantization problem (when $\alpha=0$).

Intuitively, the partial optimal transport problem goes as follows. We have an excess of offer of a certain quantity of mass at origin (supply) and
a mass that needs to be satisfied at destination (demand), so that it is not necessary to serve all the mass (demand$= \alpha\times$supply). In other words,
some $(1-\alpha)$-fraction of the mass at origin can be left non-served.
The goal is to perform this task at the cheapest transportation cost. If we represent the demand at destination by a target probability distribution $R$, we can model the supply at origin as $\frac{P}{\alpha}$, where $P$ is another probability distribution and the mass required at destination is $\alpha$ times the mass at origin. This way, a \emph{partial optimal transportation plan} is a probability measure $\Pi$
on $\mathbb{R}^d \times \mathbb{R}^d$ with first marginal in $\mathcal{R}_{1-\alpha}(P)$ and with second marginal equal to $R$, which solves the following cost minimization problem:
$$\mathcal{W}_p(\mathcal{R}_{1-\alpha}(P),R):=\min_{Q \in
\mathcal{R}_{1-\alpha}(P)}\mathcal{W}_p(Q,R) $$

 The following lemma allows us to characterize the connection between the joint distribution $\mathbb{Q}$ and the conditional distribution $\mathbb{Q}_{\widetilde{\Xi}}$ in problem~\eqref{conditional_expectation_problem} above in terms of the partial mass  problem.

\begin{lemma}\label{Lema_PMT}
Let $Q$ be a probability on $\mathbb{R}^d$  such  that $Q(\widetilde{\Xi}) = \alpha >0$ and let $Q_{\widetilde{\Xi}}$ be the $Q$-conditional probability distribution given the event $\boldsymbol{\xi} \in \widetilde{\Xi}$. Also, assume that for a given probability metric $D$, $\mathcal{R}_{1-\alpha}(Q)$ is closed for $D$ over an appropiate set of probability distributions. Then,  $Q_{\widetilde{\Xi}}$ is the unique distribution that satisfies $Q_{\widetilde{\Xi}}(\widetilde{\Xi})=1$ and $D\left(\mathcal{R}_{1-\alpha}(Q), Q_{\widetilde{\Xi}}\right) = 0$.
\end{lemma}

By way of Lemma~\eqref{Lema_PMT}, we can reformulate Problem~\eqref{conditional_expectation_problem} as follows:
\begin{subequations}\label{DRO_reform}
\begin{align}
 \inf_{\mathbf{x}\in X}  \sup_{Q_{\widetilde{\Xi}}} &\;\;
     \mathbb{E}_{Q_{\widetilde{\Xi}}} \left[f(\mathbf{x},\boldsymbol{\xi}) \right] \\
 & \text{s.t.} \;   \mathcal{W}_p^p(\mathcal{R}_{1-\alpha}( \mathbb{Q}),Q_{\widetilde{\Xi}}) = 0 \label{DRO_reform_wass}\\
 &\phantom{\text{s.t.}}\; Q_{\widetilde{\Xi}}(\widetilde{\Xi})=1
\end{align}
\end{subequations}
which now presents a form which is much more suited to our purpose, that is, to get to the DRO-type of problem~\eqref{conditional_model1} we propose. The change, nonetheless, has been essentially cosmetic, because  problem~\eqref{DRO_reform} still relies on the true \emph{joint} distribution $\mathbb{Q}$ and therefore, is of no use in practice as it stands right now. To make it practical, we need to rewrite it not in terms of the unknown $\mathbb{Q}$, but in terms of the information available to the decision maker, i.e., the sample data $\{\widehat{\boldsymbol{\xi}}_i\}_{i=1}^{N}$. For that purpose, it seems sensible and natural to replace $\mathbb{Q}$ in~\eqref{DRO_reform_wass} with its best approximation taken directly from the data, namely, the empirical measure of the sample, $\widehat{\mathbb{Q}}_N$. Logically, to accommodate the approximation, we will need to introduce a \emph{budget} $\widetilde{\rho}$ in equation~\eqref{DRO_reform_wass}, that is,
\begin{subequations}
\begin{align}\label{DRO_datadriven}
{\rm (P)} \enskip \inf_{\mathbf{x}\in X}  \sup_{Q_{\widetilde{\Xi}}} &\;\;
     \mathbb{E}_{Q_{\widetilde{\Xi}}} \left[f(\mathbf{x},\boldsymbol{\xi}) \right]  \\
 & \text{s.t.} \;   \mathcal{W}_p^p(\mathcal{R}_{1-\alpha}( \widehat{\mathbb{Q}}_N),Q_{\widetilde{\Xi}}) \leq \widetilde{\rho} \label{DRO_reform_W}\\
 &\phantom{\text{s.t.}}\; Q_{\widetilde{\Xi}}(\widetilde{\Xi})=1\label{constr_p_2}
\end{align}
\end{subequations}
Hereinafter we will use  $\widehat{\mathcal{U}}_{{N}}(\alpha,\widetilde{\rho})$ to denote the ambiguity set defined by constraints \eqref{DRO_reform_W}--\eqref{constr_p_2}. Under certain conditions, this uncertainty set enjoys nice topological properties,  as we state in Proposition \ref{th_compactness} in Appendix~\ref{proof_th_compactness}. 

Now we define what we call the \emph{minimum transportation budget}, which plays an important role in the selection of budget $\widetilde{\rho}$ in problem (P).
\begin{definition}[Minimum transportation budget]\label{MTB}
Given $\alpha >0$ in problem $\left({\rm P}\right)$, the \emph{minimum transportation budget}, which we denote as $\underline{\epsilon}_{N\alpha}$, is the $p$-Wasserstein distance  between the set $\mathcal{P}_p(\widetilde{\Xi})$ and the $(1-\alpha)$-trimming of the empirical distribution $\widehat{\mathbb{Q}}_N$ that is the \emph{closest} to that set, i.e., $ \inf \{ \mathcal{W}_p(P,Q) \;:\; P \in \mathcal{R}_{1-\alpha}( \widehat{\mathbb{Q}}_N), \; Q \in \mathcal{P}_p(\widetilde{\Xi})\}$, which is given by
\begin{equation}\label{def_MTB}
    \underline{\epsilon}_{N\alpha} =\left( \frac{1}{N\alpha}\sum_{k=1}^{\lfloor N\alpha \rfloor}{\textrm{dist}(\boldsymbol{\xi}_{k:N}, \widetilde{\Xi})^p} + \left(1-\frac{\lfloor N\alpha \rfloor}{N\alpha}\right) \textrm{dist}(\boldsymbol{\xi}_{\lceil N\alpha\rceil:N}, \widetilde{\Xi})^p\right)^{\frac{1}{p}}
\end{equation}
where $\boldsymbol{\xi}_{k:N}$ is the $k$-\emph{th} nearest data point from the sample to set $\widetilde{\Xi}$ and
$
\textrm{dist}
    (\boldsymbol{\xi}_j,
    \widetilde{\Xi}):=
    \inf_{\boldsymbol{\xi}\in
    \widetilde{\Xi}} \textrm{dist}
    (\boldsymbol{\xi}_j,\boldsymbol{\xi}) = \inf_{\boldsymbol{\xi}\in
    \widetilde{\Xi}}||\boldsymbol{\xi}_j-\boldsymbol{\xi}||.
$
If $\alpha = 0$, then $ \underline{\epsilon}_{N0} = \textrm{dist}(\boldsymbol{\xi}_{1:N}, \widetilde{\Xi})$.
\end{definition}
Importantly, the minimum transportation budget to the power of $p$, i.e., $\underline{\epsilon}_{N\alpha}^{p}$, is the minimum value of $\widetilde{\rho}$ in $\left({\rm P}\right)$ for this problem to be feasible. Furthermore, $\underline{\epsilon}_{N\alpha}$ is random, because it depends on the available
data sample, but realizes before the decision $\mathbf{x}$ is to be made. It constitutes, therefore, input data to problem  $\left({\rm P}\right)$.

We note that, if the random vector $\mathbf{y}$ takes values in a set that is independent of the feature vector $\mathbf{z}$, i.e., for all $\mathbf{z}^* \in \Xi_{\mathbf{z}} $, $\{\mathbf{y} \in \Xi_{\mathbf{y}} : \boldsymbol{\xi} = (\mathbf{z}^*, \mathbf{y}) \in \Xi\} = \Xi_{\mathbf{y}} $,
then $\textrm{dist}(\boldsymbol{\xi}_j,\widetilde{\Xi}) = \inf_{\boldsymbol{\xi}\in \widetilde{\Xi}}||\boldsymbol{\xi}_j-\boldsymbol{\xi}|| = \inf_{\boldsymbol{\xi} = (\mathbf{z},\mathbf{y}) \in \widetilde{\Xi}}||\mathbf{z}_j-\mathbf{z}||$.

Furthermore, in what follows, we assume that $\textrm{dist}(\boldsymbol{\xi}_j,\widetilde{\Xi})$ (interpreted as a random variable) conditional on $\xi_j \notin \widetilde{\Xi}$ has a continuous distribution function. This ensures that, in the case $\mathbb{Q}(\widetilde{\Xi}) = 0$, which we study in Section~\ref{case_alpha_0}, there will be exactly $K$ nearest data points to $\widetilde{\Xi}$ with probability one.

Next we present an interesting result, which deals with the inner supremum of problem (P) and adds more meaning to this problem by linking it to an alternative formulation more in the style of the Wasserstein data-driven DRO approach proposed in \cite{MohajerinEsfahani2018}, where, however, no side information is taken into account. {In fact, the distributionally robust approach to conditional stochastic optimization that is proposed in \cite{Nguyen2021} is based on this alternative formulation (see Proposition A.4 in that work)}\footnote{  Proposition~\ref{reform_potp} in this paper predates the publication of preprint~\cite{Nguyen2021}.}.  A proof of the following result can be found in Appendix~\ref{reform_potp_proof}.

\begin{proposition}\label{reform_potp}
Given $N\geqslant 1$, $\mathbb{Q}(\widetilde{\Xi}) = \alpha >0$, and any positive value of $\widetilde{\rho}$, problem \emph{(SP2)} is a relaxation of \emph{(SP1)}, where \emph{(SP1)} and \emph{(SP2)} are given by



\begin{align*}
 {\rm (SP1)}\;\; \left\{ \begin{array}{cl}  \sup_{Q} &{} \mathbb{E}_Q \left[f(\mathbf{x},\boldsymbol{\xi}) \;{\mid}\; \boldsymbol{\xi} \in \widetilde{\Xi} \right] \\ \text {s.t.}&{} \mathcal{W}_p^p(Q, \widehat{\mathbb{Q}}_N)\leqslant \widetilde{\rho}\cdot  \alpha \\ &{} Q(\widetilde{\Xi})=\alpha \end{array} \right.,  \;\; {\rm (SP2)} \left\{ \begin{array}{cl} \sup_{Q_{\widetilde{\Xi}}} &{} \mathbb{E}_{Q_{\widetilde{\Xi}}} \left[f(\mathbf{x},\boldsymbol{\xi}) \right] \\ \text {s.t.}&{}  \mathcal{W}_p^p(\mathcal{R}_{1-\alpha}( \widehat{\mathbb{Q}}_N),Q_{\widetilde{\Xi}})\leqslant \widetilde{\rho} \\ &{} Q_{\widetilde{\Xi}}(\widetilde{\Xi})=1 \end{array} \right. \end{align*}
 and where by ``relaxation'' it is meant that any solution $Q$ feasible in \emph{(SP1)} can be mapped into a solution $Q_{\widetilde{\Xi}}$ feasible in \emph{(SP2)} with the same objective function value.

Moreover, if  $\widehat{\mathbb{Q}}_{N}(\widetilde{\Xi}) = 0$ or $\alpha = 1$, then \emph{(SP1)} and \emph{(SP2)} are equivalent.

\end{proposition}

Among other things, Proposition~\ref{reform_potp} reveals that parameter $\widetilde{\rho}$ in problem (SP2), and hence in problem (P), can be understood as a cost budget \emph{per unit of transported mass}. Likewise, parameter $\alpha$ can be interpreted as the minimum amount of mass (in per unit) of the empirical distribution $\widehat{\mathbb{Q}}_N$ that must be transported to the support $\widetilde{\Xi}$. This interpretation of parameters $\widetilde{\rho}$ and $\alpha$ will be useful to follow the rationale behind the DRO solution approaches that we develop later on.

On the other hand, despite the connection between problems (SP1) and (SP2) that Proposition~\ref{reform_potp} unveils, the latter is qualitatively more amenable to further generalization and analysis. Examples of this are given by the relevant cases $\alpha = 0$, for which problem (SP1) is \emph{ill-posed}, while problem (SP2) is not, and $\alpha$ unknown, for which the use of trimming sets in (SP2) allows for a more straightforward treatment. We will deal with both cases in  Sections~\ref{case_alpha_0} and~\ref{unknown-alpha}, respectively. Before that, we provide an implementable reformulation of the proposed DRO problem (P).

\subsection{Towards a tractable reformulation of the partial mass transportation problem}
In this section, we put the proposed DRO problem (P) in a form more suited to tackle its computational implementation and solution. For this purpose, we first need to introduce a technical result whereby we characterize the trimming sets of an empirical probability measure.
\begin{lemma}\label{lemma_trimming}
Consider the sample data $\{\widehat{\boldsymbol{\xi}}_i\}_{i=1}^{N}$ and their associated empirical measure $\widehat{\mathbb{Q}}_{N} = \frac{1}{N}\sum_{i=1}^{N}{\delta_{\widehat{\boldsymbol{\xi}}_i}}$.
 If $\alpha > 0$, the set of all $(1-\alpha)$-trimmings of $\widehat{\mathbb{Q}}_{N}$ is given by all probability distributions in the form $\sum_{i=1}^{N}{b_i\delta_{\widehat{\boldsymbol{\xi}}_i}}$ such that $0\leq b_i\leq \frac{1}{N\alpha}$, $\forall i =1, \ldots, N$, and $\sum_{i=1}^{N}{b_i}=1$.  Furthermore, if $\alpha = 0$, the set $\mathcal{R}_{1-\alpha}(\widehat{\mathbb{Q}}_{N})$ of $(1-\alpha)$-trimmings of $\widehat{\mathbb{Q}}_{N}$ becomes $\mathcal{R}_{1}(\widehat{\mathbb{Q}}_{N}) = \{\sum_{i=1}^{N}{b_i\delta_{\widehat{\boldsymbol{\xi}}_i}}$ such that $b_i \geqslant 0$, $\forall i =1, \ldots, N$, and $\sum_{i=1}^{N}{b_i}=1\}$.
\end{lemma}
\proof{Proof.} If $\alpha > 0$,  the form of any $(1-\alpha)$-trimming of $\widehat{\mathbb{Q}}_{N}$ as $\sum_{i=1}^{N}{b_i\delta_{\widehat{\boldsymbol{\xi}}_i}}$, along with the condition $b_i\leq \frac{1}{N\alpha}$, follows directly from Definition~\ref{def_trimming} of a $(1-\alpha)$-trimming. Naturally, $b_{i} \geq 0$ and $\sum_{i=1}^{N}{b_i}=1$ are then required because any $(1-\alpha)$-trimming is a probability distribution.

On the other hand,  if $\alpha= 0$, the resulting trimming set $\mathcal{R}_{1}(\widehat{\mathbb{Q}}_{N})$ is simply the family of all probability distributions supported on the data points $\{\widehat{\boldsymbol{\xi}}_i\}_{i=1}^{N}$. \qed

In short, Lemma~\ref{lemma_trimming} tells us that trimming a data sample of size $N$ with level $1-\alpha$ involves reweighting the empirical distribution of such data by giving a new weight less than or equal to $\frac{1}{N\alpha}$ to each data point. Therefore, we can recast constraint $\mathcal{W}_p^p(\mathcal{R}_{1-\alpha}( \widehat{\mathbb{Q}}_N),Q_{\widetilde{\Xi}})\leqslant \widetilde{\rho}$ in problem (P) as %
 \begin{align*}
&\min_{b_i, \forall i \leqslant N} \mathcal{W}_p\left(\sum_{i=1}^N b_i \delta_{\widehat{\boldsymbol{\xi}}_i},Q_{\widetilde{\Xi}}\right)\leqslant \widetilde{\rho}^{1/p}\\
&\hspace{1cm} \text{s.t.} \ 0  \leqslant b_i \leqslant \frac{1}{N \alpha}, \enskip \forall i\leqslant N\\
&\hspace{1cm} \phantom{s.t.} \ \sum_{i=1}^N b_i=1
\end{align*}
We are now ready to introduce the main result of this section.
\begin{theorem}[Reformulation based on strong duality]\label{Theorem_partial_mass_reformulation_duality}
For  $\alpha > 0$ and any value of $\widetilde{\rho}\geqslant\underline{\epsilon}^p_{N\alpha}$, subproblem \emph{(SP2)} is equivalent to the following one:

      \begin{align*}
 {\rm (SP2')}\;\;  &  \inf_{\lambda \geqslant 0; \overline{\mu}_i, \forall i\leqslant N; \theta \in \mathbb{R}} \quad   \lambda \widetilde{\rho}+\theta +\dfrac{1}{N\alpha}\sum_{i=1}^N \overline{\mu}_i   \\
& \hspace{2cm}  \text{s.t.} \  \overline{\mu}_i+\theta\geqslant\!\!\!\! \sup_{ (\mathbf{z},\mathbf{y}) \in  \widetilde{\Xi}}\!\!\left(
f(\mathbf{x},(\mathbf{z},\mathbf{y}))-\lambda  \left\|(\mathbf{z},\mathbf{y})-(\widehat{\mathbf{z}}_{i},\widehat{\mathbf{y}}_{i})\right\|^p  \right)\!\!,  \forall i \leqslant N\\
& \hspace{2cm}  \phantom{s.t.} \ \overline{\mu}_i \geqslant 0,\enskip \forall i \leqslant N
\end{align*}

\end{theorem}

Surely the most important takeaway message of Theorem~\ref{Theorem_partial_mass_reformulation_duality} is that problem (P) is \emph{as tractable as} the standard Wasserstein-metric-based DRO formulation proposed in \cite{MohajerinEsfahani2018} and \cite{Kuhn2019}. In these two approaches, conditions under which the inner supremum in  ${\rm (SP2')}$ can be recast in a more tractable form are provided. As an example,  in Theorem~\ref{Theorem_partial_mass_reformulation} in Appendix~\ref{sup_reformulation}, we provide a more refined reformulation of  ${\rm (SP2')}$, whereby the problems we solve in Section~\ref{numerics} can be directly handled. 

In the following section, we show that problem (P) works, under certain conditions, as a statistically meaningful surrogate decision-making model for the target conditional stochastic program~\eqref{conditional_expectation_problem}.


 \section{Finite sample guarantee and asymptotic consistency}

Next we argue that the worst-case optimal expected cost provided by problem (P) for a fixed sample size $N$ and a suitable choice of parameters $(\alpha, \widetilde{\rho})$ (dependent on $N$) leads to an upper confidence bound on the out-of-sample performance  attained by the optimizers of (P) (\emph{finite sample guarantee}) and that those optimizers almost surely converge  to an optimizer of the true optimal expected cost
as $N$ grows to infinity (\emph{asymptotic consistency}).

To be more precise, the \emph{out-of-sample performance} of a given data-driven candidate solution  $\widehat{\mathbf{x}}_N$ to problem~\eqref{conditional_expectation_problem} is defined as $\mathbb{E}_{\mathbb{Q}} [f(\widehat{\mathbf{x}}_N ,\boldsymbol{\xi})\mid \boldsymbol{\xi} \in \widetilde{\Xi}] = \mathbb{E}_{\mathbb{Q}_{\widetilde{\Xi}}} [f(\widehat{\mathbf{x}}_N ,\boldsymbol{\xi})]$. We say that a data-driven method built to address problem~\eqref{conditional_expectation_problem} enjoys a \emph{finite sample guarantee}, if it produces pairs $(\widehat{\mathbf{x}}_N, \widehat{J}_{N})$ satisfying a relation in the form
   \begin{equation}\label{bounds_finite_sample_guarantee_1}
      \mathbb{Q}^N\Big[ \mathbb{E}_{\mathbb{Q}} [f(\widehat{\mathbf{x}}_N ,\boldsymbol{\xi})\mid \boldsymbol{\xi} \in \widetilde{\Xi}]\leqslant \widehat{J}_N \Big]
      \geqslant 1-\beta
  \end{equation}
   and $\widehat{J}_N$ is a \emph{certificate} for the out-of-sample
performance of $\widehat{\mathbf{x}}_N $ (i.e., an upper bound that is generally contingent on the data sample). The probability on the right-hand side
of \eqref{bounds_finite_sample_guarantee_1}, i.e.,  $1-\beta$, is known as  the \emph{reliability} of $(\widehat{\mathbf{x}}_N, \widehat{J}_N )$ and can be understood as a confidence level.

Our analysis relies on the lemma below, which immediately follows from setting $P_1:=\widehat{\mathbb{Q}}_N,Q:=\mathbb{Q}_{\widetilde{\Xi}}, P_2:=\mathbb{Q}$ in Lemma 3.13 on probability trimmings in~\cite{AgulloAntolin2018}.



\begin{lemma}\label{lemma_unifying}
Assume that $\mathbb{Q}_{\widetilde{\Xi}}, \mathbb{Q} \in \mathcal{P}_p(\mathbb{R}^d)$, and take $p \geqslant 1$, then
\begin{equation}\label{ineq_unify}
  \mathcal{W}_p(\mathcal{R}_{1-\alpha}(\widehat{\mathbb{Q}}_N),\mathbb{Q}_{\widetilde{\Xi}})\leqslant     \mathcal{W}_p(\mathcal{R}_{1-\alpha}(\mathbb{Q}),\mathbb{Q}_{\widetilde{\Xi}})+\frac{1}{\alpha^{1/p}}\mathcal{W}_p(\widehat{\mathbb{Q}}_N,\mathbb{Q})
\end{equation}
\end{lemma}
%

We notice that the term $\mathcal{W}_p(\mathcal{R}_{1-\alpha}(\mathbb{Q}),\mathbb{Q}_{\widetilde{\Xi}})$ in \eqref{ineq_unify} is not random and depends exclusively on the true distributions $\mathbb{Q}_{\widetilde{\Xi}}$, $\mathbb{Q}$, and the trimming level $\alpha$. It is, therefore, independent of the data sample (unlike the other two terms involved).

Inequality~\eqref{ineq_unify} reveals an interesting trade-off. On the one hand, the distance $\mathcal{W}_p(\mathcal{R}_{1-\alpha}(\mathbb{Q}),\mathbb{Q}_{\widetilde{\Xi}})$ diminishes as $\alpha$ decreases to zero, because the trimming set $\mathcal{R}_{1-\alpha}(\mathbb{Q})$ grows in size. On the other, the term  $\frac{1}{\alpha^{1/p}}\mathcal{W}_p(\widehat{\mathbb{Q}}_N,\mathbb{Q})$ becomes larger as $\alpha$ approaches zero. As we will see later on, controlling this trade-off is key to endowing problem (P) with performance guarantees. To this end, we will make use of Proposition~\ref{prop_measure_conc_uni} below.
\begin{assumption}\label{light_tailed_assumption_joint}
     Suppose that the true joint probability distribution $\mathbb{Q}$ is light-tailed, i.e., there exists a constant $a > p \geqslant 1$ such that $\mathbb{E}_{\mathbb{Q}}\left[ \exp(\|   \boldsymbol{\xi}\|^{a})\right]< \infty$.
\end{assumption}
\begin{proposition}[Concentration tail inequality]\label{prop_measure_conc_uni}
    Suppose that Assumption~\ref{light_tailed_assumption_joint} holds. Then, there are constants
     $c,C>0$ such that, for all $\epsilon > 0, \alpha >0$, and $N \geqslant 1$, it holds
     \begin{equation}\label{ineq_measure_conc_uni}
     \mathbb{Q}^{N} \left[ \mathcal{W}_p\left(\mathcal{R}_{1-\alpha}(\widehat{\mathbb{Q}}_N),\mathbb{Q}_{\widetilde{\Xi}} \right)  \geqslant \mathcal{W}_p(\mathcal{R}_{1-\alpha}(\mathbb{Q}),\mathbb{Q}_{\widetilde{\Xi}}) + \epsilon\right] \leqslant \beta_{p,\epsilon,\alpha}(N)
      \end{equation}
      where
  \begin{align}\label{beta_fournier_uni}
\beta_{p,\epsilon,\alpha}(N)\!=\!
  &\,\mathbb{I}_{\left\{\epsilon \leqslant 1/\alpha^{1/p}\right\}}C\! \left\{ \begin{array}{ll} \exp (-c N \, \alpha^{2}\, \epsilon^{2p}) &{}\ \hbox {if }\,p>d/2,\\ \exp (-cN(\alpha \epsilon^p/\log (2+1/\alpha\epsilon^p))^2) &{}\ \hbox {if }\,p=d/2, \\ \exp (-cN\, \alpha^{d/p}\,\epsilon^{d}) &{}\  \hbox {if }\,p\in [1,d/2), \, d > 2 \end{array}\right.\\
  &+C\exp(-cN\,\alpha^{a/p}\,\epsilon^{a})
  \mathbb{I}_{\{\epsilon > 1/\alpha^{1/p}\}}\notag
\end{align}
with $d = d_{\mathbf{z}} + d_{\mathbf{y}}$.

\end{proposition}
\proof{Proof.}
Because of Lemma~\ref{lemma_unifying} we have
$$
 \mathbb{Q}^{N}\left( \mathcal{W}_p(\mathcal{R}_{1-\alpha}(\widehat{\mathbb{Q}}_N),\mathbb{Q}_{\widetilde{\Xi}}) - \mathcal{W}_p(\mathcal{R}_{1-\alpha}(\mathbb{Q}),\mathbb{Q}_{\widetilde{\Xi}}) \geqslant \epsilon\right) \leqslant \mathbb{Q}^{N}\left(\mathcal{W}_{p}^p\left(\widehat{\mathbb{Q}}_N, \mathbb{Q}\right) \geqslant \alpha \epsilon^p\right)$$
where the right-hand side of this inequality is upper bounded by~\eqref{beta_fournier_uni} according to  \cite[Theorem 2]{Fournier2015}. \qed

Assuming $p\neq d/2$ , if we equate $\beta$ to $\beta_{p,\epsilon,\alpha}(N)$  and solving for $\epsilon$ we get:
\begin{align}\label{radius}
\epsilon_{N,p,\alpha}(\beta) {:=}\left\{
\begin{array}{ll} \Big ({\log (C\beta ^{-1}) \over cN}
\Big )^{1/2p}\frac{1}{\alpha^{1/p}} &{} \quad \text {if } N \geq {\log (C\beta ^{-1}) \over c},\quad  p>d/2, \\ \Big ({\log (C\beta ^{-1}) \over cN} \Big )^{1/d}\frac{1}{\alpha^{1/p}} &{} \quad \text {if } N \geq {\log (C\beta ^{-1}) \over c},\quad p\in [1,d/2), \, d > 2\\
\Big ({\log (C\beta ^{-1}) \over cN} \Big )^{1/a}\frac{1}{\alpha^{1/p}} &{} \quad \text {if } N < {\log (C\beta ^{-1}) \over c} \end{array}
\right.
\end{align}

In what follows, we distinguish three general setups that may appear in the real-life use of Conditional Stochastic Optimization, namely, the case $\mathbb{Q}(\widetilde{\Xi}) = \alpha > 0$ with $\alpha$ known, the case $\mathbb{Q}(\widetilde{\Xi}) = \alpha > 0$ with $\alpha$ unknown, and  the case $\mathbb{Q} \ll \lambda^{d}$ with $\mathbb{Q}(\widetilde{\Xi}) = \alpha = 0$.

\subsection{Case $\mathbb{Q}(\widetilde{\Xi}) = \alpha > 0$. Applications in data-driven decision making under contaminated samples}\label{case_alpha_pos_uni}
%

When $\mathbb{Q}(\widetilde{\Xi}) = \alpha > 0$ and \emph{known}, we can solve the following DRO problem:
\begin{subequations}
\begin{align}
{\rm (P_{(\alpha, \widetilde{\rho}_N)})} \enskip \inf_{\mathbf{x}\in X}  \sup_{Q_{\widetilde{\Xi}}} &\;\;
     \mathbb{E}_{Q_{\widetilde{\Xi}}} \left[f(\mathbf{x},\boldsymbol{\xi}) \right]\\
 & \text{s.t.} \;   \mathcal{W}_p^p(\mathcal{R}_{1-\alpha}( \widehat{\mathbb{Q}}_N),Q_{\widetilde{\Xi}}) \leq \widetilde{\rho}_N \\
 &Q_{\widetilde{\Xi}}(\widetilde{\Xi})=1
\end{align}
\end{subequations}
As we show below, problem $P_{(\alpha, \widetilde{\rho}_N)}$ enjoys a finite sample guarantee and produces solutions that are asymptotically consistent, i.e., that converge to the true solution (under complete information) given by problem~\eqref{conditional_expectation_problem}. This is somewhat hinted at by the connection between problems (SP1) and (SP2) highlighted in Proposition~\ref{reform_potp}.

\begin{theorem}[Case $\alpha > 0$: Finite sample guarantee]\label{finite_sample_alpha_pos}
  Suppose that the assumptions of Proposition~\ref{prop_measure_conc_uni} hold  and take $p \neq d/2$. Given $N \geqslant 1$ and $\alpha > 0$, choose $\beta \in \left(0,1\right)$, and determine $\epsilon_{N,p,\alpha}(\beta)$ through~\eqref{radius}. Then, for all $\widetilde{\rho}_N \geqslant \max(\epsilon^p_{N,p,\alpha}(\beta),\underline{\epsilon}_{N\alpha}^p)$, where $\underline{\epsilon}_{N\alpha}^p$ is the minimum transportation budget as in Definition~\ref{MTB},  the pair $(\widehat{\mathbf{x}}_N$, $\widehat{J}_N)$ that is solution to problem $\left({\rm P}_{(\alpha,\widetilde{\rho}_N)}\right)$ enjoys the finite sample guarantee~\eqref{bounds_finite_sample_guarantee_1}.

    \end{theorem}
\proof{Proof.} For problem $\left({\rm P}_{(\alpha,\widetilde{\rho}_N)}\right)$ to be feasible, we must have $\widetilde{\rho}_N \geqslant \underline{\epsilon}_{N\alpha}^p$.
Furthermore,  $\mathcal{W}_p(\mathcal{R}_{1-\alpha}(\mathbb{Q}),\mathbb{Q}_{\widetilde{\Xi}}) = 0$ in \eqref{ineq_measure_conc_uni} because of Lemma~\ref{Lema_PMT}. Hence, Proposition~\ref{prop_measure_conc_uni} ensures that
$\mathbb{Q}^{N}\left(\mathbb{Q}_{\widetilde{\Xi}} \in \widehat{\mathcal{U}}_N(\alpha,\widetilde{\rho}_N)\right)\geqslant 1-\beta$ for any $\widetilde{\rho}_N \geqslant \epsilon^p_{N,p,\alpha}(\beta)$. It follows then
\begin{align*}
 \mathbb{E}_{\mathbb{Q}} [f(\widehat{\mathbf{x}}_N ,\boldsymbol{\xi})\;\mid\; \boldsymbol{\xi} \in \widetilde{\Xi}]&=\mathbb{E}_{\mathbb{Q}_{\widetilde{\Xi}}}[f(\widehat{\mathbf{x}}_N
      ,\boldsymbol{\xi})]\\
      &\leqslant
      \widehat{J}_N:= \sup_{Q_{\widetilde{\Xi}}}\left\{\mathbb{E}_{Q_{\widetilde{\Xi}}}[f(\widehat{\mathbf{x}}_N
      ,\boldsymbol{\xi})]\; : \; Q_{\widetilde{\Xi}} \in \widehat{\mathcal{U}}_N(\alpha,\widetilde{\rho}_N)\right\}
\end{align*}
with probability at least $1-\beta$. \qed

We point out that, in the case $\alpha >0$, data points may fall into the set $\widetilde{\Xi}$. Logically, the contribution of these points to the minimum transportation budget $\underline{\epsilon}_{N\alpha}^p$ is null and their order (the way their tie is broken) is irrelevant to our purpose.

Now we show that the solutions of the distributionally robust optimization problem $\left({\rm P}_{(\alpha,\widetilde{\rho}_N)}\right)$ converge to the solution of the target conditional stochastic program~\eqref{conditional_expectation_problem} as $N$ increases, for a careful choice of the budget $\widetilde{\rho}_N$. This result is underpinned by the fact that, under that selection of  $\widetilde{\rho}_N$, any distribution in
$\widehat{\mathcal{U}}_N(\alpha,\widetilde{\rho}_N)$
converges to the true conditional distribution $\mathbb{Q}_{\widetilde{\Xi}}$. This is formally stated in the following lemma.

\begin{lemma}[Case $\alpha >0$: Convergence of conditional distributions]\label{lemma_convergence_distributions_alpha_pos}
 Suppose that the assumptions of Proposition~\ref{prop_measure_conc_uni} hold. Choose a sequence $\beta_N \in \left(0,1\right)$, $N \in \mathbb{N}$, such that $\sum_{N=1}^{\infty}{\beta_N} < \infty$ and $\lim_{N\rightarrow \infty}{\epsilon_{N,p,\alpha}(\beta_N)}\rightarrow 0$. Then,
 $$\mathcal{W}_p(Q^N_{\widetilde{\Xi}}, \mathbb{Q}_{\widetilde{\Xi}}) \rightarrow 0 \enskip a.s.$$
 for any sequence $Q^N_{\widetilde{\Xi}}$,  $N \in \mathbb{N}$, such that $Q^N_{\widetilde{\Xi}}\in \widehat{\mathcal{U}}_N(\alpha,\widetilde{\rho}_N)$  with $\widetilde{\rho}_N = \max(\epsilon^p_{N,p,\alpha}(\beta_N),\underline{\epsilon}_{N\alpha}^p)$.
\end{lemma}
\proof{Proof.} Take $N$ large enough and let $\widehat{Q}_{N/\widetilde{\Xi}}$ be the conditional probability distribution of $\widehat{\mathbb{Q}}_N$ given $\boldsymbol{\xi} \in \Xi$. We have
$$\mathcal{W}_p(Q^N_{\widetilde{\Xi}}, \mathbb{Q}_{\widetilde{\Xi}}) \leqslant \mathcal{W}_p(Q^N_{\widetilde{\Xi}},\widehat{Q}_{N/\widetilde{\Xi}})+ \mathcal{W}_p(\widehat{Q}_{N/\widetilde{\Xi}},\mathbb{Q}_{\widetilde{\Xi}})$$
We show that the two terms on the right-hand side of the above inequality vanish with probability one as $N$ grows to infinity. We start with $\mathcal{W}_p(\widehat{Q}_{N/\widetilde{\Xi}},\mathbb{Q}_{\widetilde{\Xi}})$.

Let $I$ denote the subset of observations  $\widehat{\boldsymbol{\xi}}_i:=(\widehat{\mathbf{z}}_{i},\widehat{\mathbf{y}}_{i})$ for
$i=1,\ldots, N$, such that $\widehat{\boldsymbol{\xi}}_i \in \widetilde{\Xi}$. It follows from the Strong Law of Large Numbers that
$\widehat{\mathbb{Q}}_N(\widetilde{\Xi}) = \frac{|I|}{N} = \alpha_N \rightarrow \alpha$ almost surely.
Besides, since the sequence $\beta_N, N \in \mathbb{N}$ is summable and $\lim_{N\rightarrow \infty}{\epsilon_N(\beta_N)}\rightarrow 0$, the Borel-Cantelli Lemma and Proposition~\ref{prop_measure_conc_uni} implies
$$\mathcal{W}_{p}\left(\mathcal{R}_{1-\alpha}(\widehat{\mathbb{Q}}_N), \mathbb{Q}_{\widetilde{\Xi}}\right) \rightarrow 0 \ a.s.$$

Then, from Lemma~\ref{Lema_PMT}, we deduce that $\mathcal{W}_p(\widehat{Q}_{N/\widetilde{\Xi}},\mathbb{Q}_{\widetilde{\Xi}}) \rightarrow 0$ with probability one.

We can deal with the term $\mathcal{W}_p(Q^N_{\widetilde{\Xi}},\widehat{Q}_{N/\widetilde{\Xi}})$ in a similar fashion, except for the subtle difference that, in this case, we require $\widetilde{\rho}_N = \max(\epsilon^p_{N,p,\alpha}(\beta_N),\underline{\epsilon}_{N\alpha}^p)$, so that, for \emph{all} $N \in \mathbb{N}$, problem $P_{(\alpha, \widetilde{\rho}_N)}$ delivers a feasible $Q^N_{\widetilde{\Xi}}$ in the sequence. Hence, in order to prove that $\mathcal{W}_p(Q^N_{\widetilde{\Xi}},\widehat{Q}_{N/\widetilde{\Xi}}) \rightarrow 0$ almost surely, we need  to show that $\lim_{N\rightarrow \infty}\underline{\epsilon}_{N\alpha} = 0$ with probability one. This is something that can be directly deduced from the definition of $\underline{\epsilon}_{N\alpha}$, namely,
 \begin{align}
\underline{\epsilon}_{N\alpha}^p&:= \mathcal{W}_p^p(\mathcal{R}_{1-\alpha}(\widehat{\mathbb{Q}}_N), \mathcal{P}_p(\widetilde{\Xi}))=\min_{Q' \in \mathcal{P}_p(\widetilde{\Xi})}\mathcal{W}_p^p(\mathcal{R}_{1-\alpha}(\widehat{\mathbb{Q}}_N), Q')\\
&\leqslant \mathcal{W}^p_{p}\left(\mathcal{R}_{1-\alpha}(\widehat{\mathbb{Q}}_N), \mathbb{Q}_{\widetilde{\Xi}}\right) \rightarrow 0 \ a.s.
\end{align}
 \qed
%

Note that, by Equation~\eqref{def_MTB} in Definition~\ref{MTB}, we have that $\underline{\epsilon}_{N\alpha} > 0$ if and only if
$$\lceil N\alpha \rceil > |I| \Leftrightarrow \frac{\lceil N\alpha \rceil}{N} > \frac{|I|}{N} = \alpha_N = \widehat{\mathbb{Q}}_N(\widetilde{\Xi}) \Leftrightarrow \alpha > \alpha_N $$

Once the convergence of $Q^{N}_{\widetilde{\Xi}}$ to the true conditional distribution $\mathbb{Q}_{\widetilde{\Xi}}$ in the $p$-Wasserstein metric has been established by the previous lemma,  the following asymptotic consistency result, which is analogous to that of \cite[Theorem 3.6]{MohajerinEsfahani2018}, can also be derived.
\begin{theorem}[Asymptotic consistency]\label{theorem_consistency_trimmed_alphapos} Consider that the conditions of Theorem~\ref{finite_sample_alpha_pos} hold.
Take  a sequence $\widetilde{\rho}_N$ as in Lemma~\ref{lemma_convergence_distributions_alpha_pos}.   Then, we have
 \begin{enumerate}
 \item[(i)] If  for any fixed value $\mathbf{x}\in X$,  $f(\mathbf{x},\boldsymbol{\xi})$ is continuous in $\boldsymbol{\xi}$ and there is  $L\geqslant 0$ such that
     $|f(\mathbf{x},\boldsymbol{\xi})|\leqslant L(1+ \left\|\boldsymbol{\xi} \right\|^p) $  for all $\mathbf{x}\in X$ and $\boldsymbol{\xi} \in \widetilde{\Xi}$, then we have that $\widehat{J}_N \rightarrow J^*$ almost surely when $N$ grows to infinity.
     \item[(ii)]  If the assumptions in (i) are satisfied,  $f(\mathbf{x},\boldsymbol{\xi})$ is   lower semicontinuous on $X$ for any fixed $\boldsymbol{\xi}\in \widetilde{\Xi}$, and the feasible set $X$ is closed,  then we have that any accumulation point of the sequence $\{ \widehat{\mathbf{x}}_N \}_N$ is almost surely an optimal solution of problem~\eqref{conditional_expectation_problem}.

 \end{enumerate}
\end{theorem}

\proof{Proof.}
We omit the proof, because it is essentially the same as the one in \cite[Theorem 3.6]{MohajerinEsfahani2018}, except that, since we are working with $p \geqslant 1$, we additionally require that $f(\mathbf{x},\boldsymbol{\xi})$ be continuous in $\boldsymbol{\xi}$ so that we can make use of Theorem 7.12 from \cite{Villani2003}.

In the following remark, we show how problem $P_{(\alpha, \widetilde{\rho}_N)}$ can be used to make distributionally robust decisions in a context where the data available to the decision maker is contaminated.

\begin{remark}[Data-driven decision-making under contaminated samples]
Suppose that the dataset $\widehat{\boldsymbol{\xi}}_i:=(\widehat{\mathbf{z}}_{i},\widehat{\mathbf{y}}_{i})$ for
$i=1,\ldots, N$ is composed of \emph{correct} and \emph{contaminated} samples. The decision maker only knows that a sample is correct with probability $\alpha$ and contaminated with probability $1-\alpha$, but  does not know which type each sample belongs to. Thus, the data have been generated from a mixture distribution given by $P=\alpha Q^*+(1-\alpha)R$, where $Q^*$ is the correct distribution and $R$ a contamination.

In our context, this is equivalent to stating that
$Q^*\in \mathcal{R}_{1-\alpha}(P)$, which, in turn, can be formulated as
$\mathcal{W}_p(\mathcal{R}_{1-\alpha}(P), Q^*)=0$. Since we only have limited information on $P$ in the form of the empirical distribution $\widehat{P}_N$, we propose to solve  problem ${\rm P}_{(\alpha, \widetilde{\rho}_N)}$, that is,
\begin{subequations}
\begin{align}\label{DRO_contamination1}
\inf_{\mathbf{x}\in X}  \sup_{Q} &\;\;
     \mathbb{E}_{Q} \left[f(\mathbf{x},\boldsymbol{\xi}) \right]\\
 & \text{s.t.} \;   \mathcal{W}_p^p(\mathcal{R}_{1-\alpha}( \widehat{P}_N),Q) \leq \widetilde{\rho}_N \label{DRO_contamination2}
\end{align}
\end{subequations}
where we have assumed that the correct distribution $Q^*$, the contamination $R$ and the data-generating distribution $P$ are all supported on $\Xi$.

The decision maker can profit from the finite sample guarantee that the solution to  problem~\eqref{DRO_contamination1}--\eqref{DRO_contamination2} satisfies as per Theorem~\ref{finite_sample_alpha_pos}, with $\widetilde{\rho}_N \geqslant \epsilon^p_{N,p,\alpha}(\beta)$, $\beta \in (0,1)$, since $\underline{\epsilon}^p_{N\alpha} = 0$ in this case. Furthermore, if we choose a summable sequence of $\beta_N \in (0,1)$, $N \in \mathbb{N}$, such that $\lim_{N\rightarrow \infty}\epsilon_{N}(\beta_N) = 0$, then we have that
\begin{equation}\label{almost_sure_alpha_pos}
 P^{\infty}\left(\lim_{N\rightarrow \infty}\mathcal{W}_p\left(\mathcal{R}_{1-\alpha}(\widehat{P}_N),Q^*\right)=0\right) = 1
\end{equation}

In plain words, for $N$ large enough, the decision vector $\mathbf{x}$ is being optimized by way of problem~\eqref{DRO_contamination1}--\eqref{DRO_contamination2} over the ``smallest'' ambiguity set that almost surely contains the correct distribution $Q^*$ of the data (in the absence of any other information on $Q^*$). In fact, this means our DRO approach  deals with contaminated samples in a way that is distinctly more convenient than that of
\cite{RuidiChenphdthesis} and \cite{Farokhi2021}.  Essentially, they suggest optimizing over a 1-Wasserstein ball centered at $\widehat{P}_N$ of radius $\widetilde{\rho}$, that is,
\begin{subequations}
\begin{align}\label{DRO_contamination_farokhi}
 \enskip \inf_{\mathbf{x}\in X}  \sup_{Q} &\;\;
     \mathbb{E}_{Q} \left[f(\mathbf{x},\boldsymbol{\xi}) \right]\\
 & \text{s.t.} \;   \mathcal{W}_1( \widehat{P}_N,Q) \leq \widetilde{\rho} \label{DRO_stand_budget}
\end{align}
\end{subequations}
under the argument that for $\rho$ sufficiently large, the Wasserstein ball contains the true distribution of the data $Q^*$ with a certain confidence level. For instance, the author of \cite{Farokhi2021} uses the triangle inequality and the convexity property of the Wasserstein distance to establish that $\mathcal{W}_1(\widehat{P}_{N}, Q^*) \leqslant \mathcal{W}_1(\widehat{P}_{N}, P) + (1-\alpha)\mathcal{W}_1(R, Q^*)$, so that the extra budget $(1-\alpha)\mathcal{W}_1(R, Q^*)$ would ensure that $Q^*$ is within the Wasserstein ball with a given confidence level (a similar argument is made in \cite{RuidiChenphdthesis}). In practice, though, this extra budget as such cannot be computed, because neither the correct distribution $Q^*$ nor the contamination $R$ are known to the decision maker. However, our approach naturally encodes it in the ambiguity set~\eqref{DRO_contamination2}. Indeed, for $N$ large enough, result~\eqref{almost_sure_alpha_pos} tells us that the correct distribution $Q^*$ belongs, almost surely, to the $(1-\alpha)$-trimming set of the empirical distribution $\widehat{P}_N$. It follows precisely from this and Proposition~\ref{Prop_trim_cont} in Appendix~\ref{appen:proofs} that $\mathcal{W}_{p}(\widehat{P}_N, Q^*) \rightarrow \mathcal{W}_{p}(\alpha Q^*+(1-\alpha)R, Q^*) \leqslant \alpha \mathcal{W}_{p}(Q^*, Q^*) + (1-\alpha)\mathcal{W}_p(R,Q^*)$, i.e., $\mathcal{W}_{p}(\widehat{P}_N, Q^*) \leqslant (1-\alpha)\mathcal{W}_p(R,Q^*)$.

In short, our approach offers probabilistic guarantees in the finite-sample regime and, in the asymptotic one, naturally exploits all the information we have on $Q^*$, namely, $Q^* \in \mathcal{R}_{1-\alpha}(P)$, to robustify the decision $\mathbf{x}$ under contamination.

\end{remark}

\subsection{The case of unknown $\mathbb{Q}(\widetilde{\Xi}) = \alpha > 0$.}\label{unknown-alpha}
In this section, we  discuss how we can use the proposed DRO approach to deal with the case in which $\mathbb{Q}(\widetilde{\Xi}) = \alpha > 0$ is unknown. For this purpose, we first introduce a proposition that will allows us to design a distributionally robust strategy to tackle problem~\eqref{conditional_expectation_problem} by means of problem (P).

\begin{proposition}\label{prop_equivalencemass}
Suppose that $\mathbb{Q}(\widetilde{\Xi}) = \alpha >0$. Take $0 < \alpha' < \alpha$ and any positive value of $\widetilde{\rho}$. Given $N \geqslant 1$, the following problem %
\begin{subequations}
  \begin{align*}
 \emph{(SP3)}\;\;     \sup_{Q_{\widetilde{\Xi}}} &\;\;
     \mathbb{E}_{Q_{\widetilde{\Xi}}} \left[f(\mathbf{x},\boldsymbol{\xi}) \right]\\
 & \text{s.t.} \;   \mathcal{W}_p^p(\mathcal{R}_{1-\alpha'}( \widehat{\mathbb{Q}}_N),Q_{\widetilde{\Xi}})\leqslant \widetilde{\rho}\\
 & \phantom{s.t.} \; Q_{\widetilde{\Xi}}(\widetilde{\Xi}) = 1
\end{align*}
\end{subequations}
is either fully equivalent to \emph(SP2), if $\frac{1}{N} \geqslant \alpha$ or a relaxation otherwise.
\end{proposition}
\proof{Proof.} The proof of the proposition is trivial and directly follows from the fact that $\mathcal{R}_{1-\alpha}(\widehat{\mathbb{Q}}_N) \subset \mathcal{R}_{1-\alpha'}(\widehat{\mathbb{Q}}_N)$, if $\alpha' \leqslant \alpha$, and that $\mathcal{R}_{1-\alpha}(\widehat{\mathbb{Q}}_N) = \mathcal{R}_{1-\alpha'}(\widehat{\mathbb{Q}}_N)$ if, besides, $\frac{1}{N\alpha} \geqslant 1$. \qed

Based on Proposition~\ref{prop_equivalencemass}, we could use the following two-step \emph{safe} strategy to handle the case of unknown $\mathbb{Q}(\widetilde{\Xi}) = \alpha > 0$:
\begin{enumerate}
    \item First, solve the following uncertainty quantification problem (see \cite{Gao2016,MohajerinEsfahani2018} for further details),
\begin{equation}\label{UQ}
\alpha_N:=\inf_{Q \in \mathbb{B}_{\epsilon_N}(\widehat{\mathbb{Q}}_N)} Q(\boldsymbol{\xi}\in \widetilde{\Xi})=1-\sup_{Q \in \mathbb{B}_{\epsilon_N}(\widehat{\mathbb{Q}}_N)} Q(\boldsymbol{\xi}\notin \widetilde{\Xi})
\end{equation}
where the radius $\epsilon_N$ of the Wasserstein ball has been chosen so that $\alpha_N$ represents the minimum probability that the joint true distribution $\mathbb{Q}$ of the data assigns to the event $\boldsymbol{\xi} \in \widetilde{\Xi}$ with confidence $1-\beta_N$, $\beta_N \in (0,1)$.

\item Next, solve problem $(P_{(\alpha_N, \widetilde{\rho}_N)})$, that is,
\begin{subequations}
\begin{align}
 \inf_{\mathbf{x}\in X}  \sup_{Q_{\widetilde{\Xi}}} &\;\;
     \mathbb{E}_{Q_{\widetilde{\Xi}}} \left[f(\mathbf{x},\boldsymbol{\xi}) \right]\\
 & \text{s.t.} \   \mathcal{W}_p^p(\mathcal{R}_{1-\alpha_N}( \widehat{\mathbb{Q}}_N),Q_{\widetilde{\Xi}}) \leq \widetilde{\rho}_N \label{pW_const}\\
 &\phantom{\text{s.t.}}\ Q_{\widetilde{\Xi}}(\widetilde{\Xi})=1
\end{align}
\end{subequations}
with $\widetilde{\rho}_N \geqslant \epsilon^p_{N}(\beta_N)/\alpha_N$.
%
%
\end{enumerate}

Now suppose that $\mathbb{Q} \in \mathbb{B}_{\epsilon_N(\beta_N)}(\widehat{\mathbb{Q}}_N)$ and therefore, $\alpha_N \leqslant \alpha$ (this is a random event that occurs with probability at least $1-\beta_N$). According to Lemma~\ref{lemma_unifying}, we have
\begin{align*}
 \alpha^{1/p} \mathcal{W}_{p}\left(\mathcal{R}_{1-\alpha}(\widehat{\mathbb{Q}}_N), \mathbb{Q}_{\widetilde{\Xi}}\right) &\leqslant \mathcal{W}_{p}\left(\widehat{\mathbb{Q}}_N, \mathbb{Q}\right) \leqslant \epsilon_N(\beta_N)\\
 \mathcal{W}^{p}_{p}\left(\mathcal{R}_{1-\alpha_N}(\widehat{\mathbb{Q}}_N), \mathbb{Q}_{\widetilde{\Xi}}\right) \leqslant \mathcal{W}^{p}_{p}\left(\mathcal{R}_{1-\alpha}(\widehat{\mathbb{Q}}_N), \mathbb{Q}_{\widetilde{\Xi}}\right) &\leqslant \frac{\epsilon_N^p(\beta_N)}{\alpha}\leqslant \frac{\epsilon_N^p(\beta_N)}{\alpha_N} = \widetilde{\rho}_N\\
\end{align*}
Hence, $\mathbb{Q}_{\widetilde{\Xi}} \in \widehat{\mathcal{U}}_N(\alpha_N, \widetilde{\rho}_N)$ with probability at least $1-\beta_N$. In other words, the two-step procedure here described does not degrade the reliability of the DRO solution. Furthermore, the minimum transportation budget $\underline{\epsilon}_{N\alpha_N}$ that makes problem $(P_{(\alpha_N, \widetilde{\rho}_N)})$ feasible is always zero here, if the event $\boldsymbol{\xi} \in \widetilde{\Xi}$ has been observed at least once. This is so because the uncertainty quantification problem of step 1 ensures that $\alpha_N$ is lower than or equal to the fraction of training data points falling in $\widetilde{\Xi}$. Moreover, when $N$ grows to infinity, this uncertainty quantification problem reduces to computing such a fraction of points, which, by the Strong Law of Large Numbers converges to the real $\alpha$, i.e., $\alpha_N \rightarrow \alpha$ with probability one. Therefore, in the asymptotic regime, this case resembles that of known $\alpha >0$.

 We notice, however, that, in practice, setting $\widetilde{\rho}_N \geqslant \epsilon^p_{N}(\beta_N)/\alpha_N$ may result in too large budgets  $\widetilde{\rho}_N$, and thus, in overly conservative solutions, because, as $\epsilon_{N}$ is increased,
  $\alpha_N$ decreases to zero. 
  For this reason, in the numerical experiments of the electronic companion,
 we describe an alternative data-driven procedure to address the case $\alpha >0$, in which we simply set $\alpha_N = \widehat{\mathbb{Q}}_N(\widetilde{\Xi})$ in problem $(P_{(\alpha_N, \widetilde{\rho}_N)})$ and use the data to tune parameter $\widetilde{\rho}_N$.

\subsection{The case $\mathbb{Q} \ll \lambda^{d}$ and $\mathbb{Q}(\widetilde{\Xi}) = \alpha = 0$.}\label{case_alpha_0}
Suppose that the true joint distribution $\mathbb{Q}$ governing the random vector $\boldsymbol{\xi}:=(\mathbf{z},\mathbf{y})$ admits a density function with respect to the Lebesgue measure $\lambda^d$, with $d = d_{\mathbf{z}}+d_{\mathbf{y}}$. Without loss of generality, consider the event $\boldsymbol{\xi} \in \widetilde{\Xi}$, where $\widetilde{\Xi}$ is defined as
$\widetilde{\Xi}=\{ \boldsymbol{\xi}=(\mathbf{z},\mathbf{y}) \in \Xi \; : \; \mathbf{z}= \mathbf{z^*}\}$. This means that $\mathbb{Q}(\widetilde{\Xi}) = \alpha = 0$.

Therefore, our focus in this case is on the particular variant of problem~\eqref{conditional_expectation_problem} given by
\begin{equation}\label{conditional_expectation_problem_Case1}
    J^*:=\inf_{\mathbf{x}\in X} \mathbb{E}_{\mathbb{Q}} \left[ f(\mathbf{x},(\mathbf{z},\mathbf{y})) \;{\mid}\; \mathbf{z} =  \mathbf{z^*} \right]
\end{equation}
Problem~\eqref{conditional_expectation_problem_Case1} has become a central object of study in what has recently come to be known as \emph{Prescriptive Stochastic Programming or Conditional Stochastic Optimization}, (see, e.g.,
\cite{Ban2019,BertsimasMCCORD2018,Bertsimas2019,Bertsimas2019b,bertsimas2017bootstrap,Diao2020,PangHo2019,sen2018learning}, and references therein).

Devising a  DRO approach to problem~\eqref{conditional_expectation_problem_Case1} using the standard Wasserstein ball $\mathcal{W}_p(\widehat{\mathbb{Q}}_N, Q)\leqslant \varepsilon$ is of no use here, because any point from the support of $\widehat{\mathbb{Q}}_N$ with an \emph{arbitrarily small} mass can be transported to the set $\widetilde{\Xi}$ at an arbitrarily small cost in terms of $\mathcal{W}_p(\widehat{\mathbb{Q}}_N, Q)$. This way, one could always place this arbitrarily small particle at a point $(\mathbf{z^*}, \mathbf{y'}) \in \underset{(\mathbf{z},\mathbf{y}) \in \widetilde{\Xi}}{\argmax}\enskip f(\mathbf{x}, (\mathbf{z},\mathbf{y}))$.
In contrast, problem (P), which is based on \emph{partial} mass transportation, offers a richer framework to  seek for a distributional robust solution to~\eqref{conditional_expectation_problem_Case1}. To see this, consider again the inequality~\eqref{ineq_unify}. If we could set $\alpha = 0$, the term $\mathcal{W}_p(\mathcal{R}_{1-\alpha}(\mathbb{Q}),\mathbb{Q}_{\widetilde{\Xi}})$ would vanish, because we could take random variables $\boldsymbol{\xi} \sim \mathbb{Q}_{\widetilde{\Xi}}$, $\boldsymbol{\xi}_m \sim \mathbb{Q}_m \in \mathcal{R}_{1}(\mathbb{Q}), m \in \mathbb{N}$, such that $\mathcal{W}_p(\mathbb{Q}_m, \mathbb{Q}_{\widetilde{\Xi}}) \rightarrow 0$. Unfortunately, fixing $\alpha$ to zero is not a real option due to the term $\frac{1}{\alpha^{1/p}}\mathcal{W}_p(\widehat{\mathbb{Q}}_N,\mathbb{Q})$ in the inequality. Therefore, what we propose instead is to solve a sequence of optimization problems in the form
\begin{subequations}
 \begin{align}
\left({\text P}_{(\alpha_N,\widetilde{\rho}_N)}\right) \enskip \inf_{\mathbf{x}\in X}  \sup_{Q_{\widetilde{\Xi}}} &\;\;
     \mathbb{E}_{Q_{\widetilde{\Xi}}} \left[f(\mathbf{x},\boldsymbol{\xi}) \right]\label{DRO_datadriven_alpha_N_OF}\\
 & \text{s.t.} \;   \mathcal{W}_p^p(\mathcal{R}_{1-\alpha_N}( \widehat{\mathbb{Q}}_N),Q_{\widetilde{\Xi}}) \leq \widetilde{\rho}_N\\
 &Q_{\widetilde{\Xi}}(\widetilde{\Xi})=1
\end{align}
\end{subequations}
with both $\alpha_N$ and $\widetilde{\rho}_N$ tending to zero appropriately as $N$ increases. Next we show that, under certain conditions, problem $\left({\text P}_{(\alpha_N,\widetilde{\rho}_N)}\right)$ enjoys a finite sample guarantee and is asymptotically consistent.

\begin{assumption}[Condition (3.6)  from \cite{falk}]\label{assump_hellinger}
Let $B(\mathbf{z}^*, r):= \{\mathbf{z} \in \Xi_{\mathbf{z}} : ||\mathbf{z}-\mathbf{z}^*|| \leqslant r\}$ denote the closed ball in $\mathbb{R}^{d_{\mathbf{z}}}$ with center $\mathbf{z}^*$ and radius $r$.
The random vector $\boldsymbol{\xi}:= (\mathbf{z},\mathbf{y})$ has a joint density $\phi$
that verifies the following for some $r_0 > 0$.
\begin{enumerate}
\item It admits uniformly for $r \in [0, r_0]$ and
$\mathbf{y}\in \mathbb{R}^{d_{\mathbf{y}}}$
the following expansion:
\begin{equation}
    \phi(\mathbf{z}^*+r \,\mathbf{u}, \mathbf{y})=
     \phi(\mathbf{z}^*, \mathbf{y})\left[1+r \langle  \mathbf{u}, \ell_1(\mathbf{y})\rangle+ O(r^2\ell_2(\mathbf{y}))\right]
\end{equation}
where $\mathbf{u} \in \mathbb{R}^{d_{\mathbf{z}}}$ with $||\mathbf{u}||= 1$, and where $\ell_1:\mathbb{R}^{d_{\mathbf{y}}}\rightarrow \mathbb{R}^{d_{\mathbf{z}}}$ and
$\ell_2:\mathbb{R}^{d_{\mathbf{y}}}\rightarrow \mathbb{R}$ satisfy
$\int(||\ell_1(\mathbf{y})||^2
+|\ell_2(\mathbf{y})|^2) \phi(\mathbf{z}^*,\mathbf{y}) d\mathbf{y}<\infty$.
\item The marginal density of $\mathbf{z}$ is bounded away from zero in $B(\mathbf{z}^*, r_0)$.

\end{enumerate}
\end{assumption}

\begin{assumption}[Regularity and boundedness]\label{assumption_knn}
 We assume that
\begin{enumerate}
\item There exists $\widetilde{C}>0$ and $r_0 >0$ such that $\mathbb{P}(\|\mathbf{z}^*-\mathbf{z} \| \leqslant r) \geqslant \widetilde{C} r^{d_{\mathbf{z}}}$, for all
$0 < r \leqslant r_0$.

\item The uncertainty $\mathbf{y}$ is bounded, that is, $\|\mathbf{y}\| \leqslant M$ a.s. for some constant $M >0$.
\end{enumerate}
\end{assumption}

We note that Assumption~\ref{assumption_knn}.1 is automatically implied by Assumption~\ref{assump_hellinger}, but we explicitly state it here for  ease of readability.  Furthermore, under the boundedness condition  established in Assumption~\ref{assumption_knn}.2, Assumption 2 is satisfied, for example,  by a twice differentiable joint density $\phi(\mathbf{z},\mathbf{y})$  with continuous and bounded partial derivatives in $B(\mathbf{z}^*,r) \times \Xi_{\mathbf{y}}$ and bounded away from zero in that set. These are standard regularity conditions in the technical literature on kernel density estimation and regression \citep{PangHo2019}.

\begin{theorem}[Case $\alpha = 0$: Finite sample guarantee]\label{finite_sample_theorem_conditional_trimmed_uni}
 Suppose that Assumptions~\ref{assump_hellinger},~\ref{assumption_knn} and those of Proposition~\ref{prop_measure_conc_uni} hold. Set $\alpha_0:= \widetilde{C}r_0^{d_{\mathbf{z}}}$.  Given $N \geqslant 1$, choose $\alpha_{N} \in (0,\alpha_0]$, $\beta \in \left(0,1\right)$, and determine $\epsilon_{N,p,\alpha_N}(\beta)$ through~\eqref{radius}.

 Then, for all
 \begin{equation}\label{rho_finitesampleguarantee}
   \widetilde{\rho}_N \geqslant \max\left[\left(\epsilon_{N,p,\alpha_N}(\beta)+ O\left(\alpha_N^{\min\{1,\ 2/p\}/d_{\mathbf{z}}}\right)\right)^p,\, \underline{\epsilon}^p_{N\alpha_N}\right]
 \end{equation}
 we have that  the pair $(\widehat{\mathbf{x}}_N$, $\widehat{J}_N)$ delivered by problem $\left({\rm P}_{(\alpha_N,\widetilde{\rho}_N)}\right)$ with parameters $\widetilde{\rho}_N $ and $\alpha_N$  enjoys the finite sample guarantee~\eqref{bounds_finite_sample_guarantee_1}.
\end{theorem}
\proof{Proof.}
For problem $\left({\text P}_{(\alpha_N,\widetilde{\rho}_N)}\right)$ to be feasible, we need $\widetilde{\rho}_N \geqslant \underline{\epsilon}^p_{N\alpha_N}$.

The proof essentially relies on upper bounding the term $\mathcal{W}_p(\mathcal{R}_{1-\alpha}(\mathbb{Q}),\mathbb{Q}_{\widetilde{\Xi}})$ that appears in Equation~\eqref{ineq_measure_conc_uni} of Proposition~\ref{prop_measure_conc_uni}. To that end,
define $\alpha(r)= \widetilde{C}r^{d_{\mathbf{z}}}$, for all
$0 < r \leqslant r_0$. Set $\alpha_0:= \alpha(r_0)$. Let  $\mathbb{Q}_{B(\mathbf{z}^*, r)\times \Xi_{\mathbf{y}}}$ be  the probability measure of $(\mathbf{z},\mathbf{y})$ conditional on  $(\mathbf{z},\mathbf{y}) \in B(\mathbf{z}^*, r)\times \Xi_{\mathbf{y}}$ and let $\mathbb{Q}_{B(\mathbf{z}^*, r)}$ be its $\mathbf{y}$-marginal. Note that, by Assumption~\ref{assumption_knn}.1, $\mathbb{Q}_{B(\mathbf{z}^*, r) \times \Xi_{\mathbf{y}}} \in \mathcal{R}_{1-\alpha(r)}(\mathbb{Q})$ provided that $0 < r \leqslant r_0$.

Furthermore, according to Theorem 3.5.2 in \cite{falk}, there exists a positive constant $A$ such that
$${\rm Hell}(\mathbb{Q}_{B(\mathbf{z}^*, r)}, \mathbb{Q}_{\widetilde{\Xi}}) \leqslant Ar^2$$
uniformly for $0 < r < r_0$, where Hell stands for \emph{Hellinger distance}.

From Equation (5.1) in \cite{Santambrogio2015} and Assumption~\ref{assumption_knn}.2 we know that
$$\mathcal{W}_p(\mathbb{Q}_{B(\mathbf{z}^*, r)}, \mathbb{Q}_{\widetilde{\Xi}})\leqslant M^{\frac{p-1}{p}}\mathcal{W}_1(\mathbb{Q}_{B(\mathbf{z}^*, r)}, \mathbb{Q}_{\widetilde{\Xi}})^{1/p}$$
In turn, from \cite{gibbs} we have that $\mathcal{W}_1(\mathbb{Q}_{B(\mathbf{z}^*, r)}, \mathbb{Q}_{\widetilde{\Xi}})\leqslant M\cdot {\rm Hell}(\mathbb{Q}_{B(\mathbf{z}^*, r)}, \mathbb{Q}_{\widetilde{\Xi}})$. Hence,
\begin{align*}
\mathcal{W}^p_p(\mathbb{Q}_{B(\mathbf{z}^*, r)}, \mathbb{Q}_{\widetilde{\Xi}})&\leqslant M^p {\rm Hell}(\mathbb{Q}_{B(\mathbf{z}^*, r)}, \mathbb{Q}_{\widetilde{\Xi}})\\
\mathcal{W}_p(\mathbb{Q}_{B(\mathbf{z}^*, r)},\mathbb{Q}_{\widetilde{\Xi}})& \leqslant M A^{1/p}r^{2/p}, \enskip 0 < r \leqslant r_0
\end{align*}
Thus,
$$\mathcal{W}_p(\mathbb{Q}_{B(\mathbf{z}^*, r) \times \Xi_{\mathbf{y}}},\mathbb{Q}_{\widetilde{\Xi}}) \leqslant r+ M A^{1/p}r^{2/p}, \enskip 0 < r \leqslant r_0$$

Since $\mathbb{Q}_{B(\mathbf{z}^*, r)\times \Xi_{\mathbf{y}}} \in \mathcal{R}_{1-\alpha(r)}(\mathbb{Q})$ for all $0<r\leqslant r_0$, it holds
$$\mathcal{W}_p(\mathcal{R}_{1-\alpha(r)}(\mathbb{Q}),\mathbb{Q}_{\widetilde{\Xi}}) \leqslant \mathcal{W}_p(\mathbb{Q}_{B(\mathbf{z}^*, r) \times \Xi_{\mathbf{y}}},\mathbb{Q}_{\widetilde{\Xi}}) \leqslant r+ M A^{1/p}r^{2/p}$$
which we can express in terms of $\alpha$ as
\begin{align*}
 \mathcal{W}_p(\mathcal{R}_{1-\alpha}(\mathbb{Q}),\mathbb{Q}_{\widetilde{\Xi}}) &\leqslant \frac{\alpha^{1/d_{\mathbf{z}}}}{\widetilde{C}^{1/d_{\mathbf{z}}}} + A^{1/p}M \frac{\alpha^{2/(pd_{\mathbf{z}})}}{\widetilde{C}^{2/(pd_{\mathbf{z}})}}  \\
 \mathcal{W}_p(\mathcal{R}_{1-\alpha}(\mathbb{Q}),\mathbb{Q}_{\widetilde{\Xi}})&= O\left(\alpha^{\min\{1,\ 2/p\}/d_{\mathbf{z}}}\right)
\end{align*}
provided that $0 < \alpha \leqslant \alpha_0$.
\qed

\begin{remark}
There are conditions on the smoothness of the true joint distribution $\mathbb{Q}$ around $\mathbf{z} =\mathbf{z}^*$, other than those stated in Assumptions~\ref{assump_hellinger} and~\ref{assumption_knn}, for which we can also upper bound the distance $\mathcal{W}_p(\mathcal{R}_{1-\alpha}(\mathbb{Q}),\mathbb{Q}_{\widetilde{\Xi}})$. We provide below two examples of these conditions, which have been invoked in  \cite{kannan2020saaresiduals,kannan2020droresiduals} and \cite{Bertsimas2019}, respectively, and neither of which requires the boundedness of the uncertainty $\mathbf{y}$.
\end{remark}
\begin{example}\label{ex_guzi}
Suppose that the true data-generating model is given by $\mathbf{y} = f^*(\mathbf{z}) +\mathbf{e}$, where $f^*(\mathbf{z}'):=\mathbb{E}[ \mathbf{y} \mid \mathbf{z}=\mathbf{z}']$ is the regression function and $\mathbf{e}$ is a zero-mean random error. Furthermore, suppose that Assumption~\ref{assumption_knn}.1 holds and there exists a positive constant $L$ such that $\|f^*(\mathbf{z}')-f^*(\mathbf{z})\| \leqslant L\|\mathbf{z}'-\mathbf{z}\|$, for all  $0 \leqslant \|\mathbf{z}'-\mathbf{z}\| \leqslant r_0$.

Take $\alpha(r)= \widetilde{C}r^{d_{\mathbf{z}}}$, for all
$0 < r \leqslant r_0$ and set $\alpha_0:= \alpha(r_0)$. With  abuse of notation, we can write for any event within $B(\mathbf{z}^*, r) \times \Xi_{\mathbf{y}}$
$$\mathbb{Q}_{B(\mathbf{z}^*, r) \times \Xi_{\mathbf{y}}}(d\mathbf{z}, d\mathbf{y}) = \frac{1}{\mathbb{P}(B(\mathbf{z}^*, r))} \mathbb{Q}(d\mathbf{z},d\mathbf{y}) = \frac{1}{\mathbb{Q}_{\mathbf{z}}(B(\mathbf{z}^*, r))} \mathbb{Q}_{\mathbf{z}=\mathbf{z}'}(d\mathbf{y}) \mathbb{Q}_{\mathbf{z}}(d\mathbf{z}')$$
where $\mathbb{Q}_{\mathbf{z}}$ is the probability law of the feature vector $\mathbf{z}$ and $\mathbb{Q}_{\mathbf{z}=\mathbf{z}'}$ is the conditional measure of $\mathbb{Q}$ given that $\mathbf{z}=\mathbf{z}'$.

Since $\mathbb{Q}_{B(\mathbf{z}^*, r)\times \Xi_{\mathbf{y}}} \in \mathcal{R}_{1-\alpha(r)}(\mathbb{Q})$ for all $0<r\leqslant r_0$, by the convexity of the Wasserstein distance, we have
\begin{align*}
\mathcal{W}_p(\mathcal{R}_{1-\alpha}&(\mathbb{Q}),\mathbb{Q}_{\widetilde{\Xi}}) \leqslant \mathcal{W}_p(\mathbb{Q}_{B(\mathbf{z}^*, r) \times \Xi_{\mathbf{y}}},\mathbb{Q}_{\widetilde{\Xi}}) \\
& \leqslant \int_{B(\mathbf{z}^*, r)}{\left[\|\mathbf{z}'-\mathbf{z}^*\|+\mathcal{W}_p(\mathbb{Q}_{\mathbf{z}=\mathbf{z}'},\mathbb{Q}_{\widetilde{\Xi}})\right] \frac{\mathbb{Q}_{\mathbf{z}}(d\mathbf{z}')}{\mathbb{Q}_{\mathbf{z}}(B(\mathbf{z}^*, r))}}\\
& =  \int_{B(\mathbf{z}^*, r)}{\left[\|\mathbf{z}'-\mathbf{z}^*\|+\mathcal{W}_p(f^*(\mathbf{z}') +\mathbf{e},f^*(\mathbf{z}^*) +\mathbf{e})\right] \frac{\mathbb{Q}_{\mathbf{z}}(d\mathbf{z}')}{\mathbb{Q}_{\mathbf{z}}(B(\mathbf{z}^*, r))}}\\
& \leqslant \int_{B(\mathbf{z}^*, r)}{\left[\|\mathbf{z}'-\mathbf{z}^*\|+\|f^*(\mathbf{z}')-f^*(\mathbf{z}^*)\|\right] \frac{\mathbb{Q}_{\mathbf{z}}(d\mathbf{z}')}{\mathbb{Q}_{\mathbf{z}}(B(\mathbf{z}^*, r))}}\\
& \leqslant (1+L) \int_{B(\mathbf{z}^*, r)}{\|\mathbf{z}'-\mathbf{z}^*\| \frac{\mathbb{Q}_{\mathbf{z}}(d\mathbf{z}')}{\mathbb{Q}_{\mathbf{z}}(B(\mathbf{z}^*, r))}} = (1+L) O(r) = O(\alpha^{1/d_{\mathbf{z}}})
\end{align*}
for all $0 < \alpha \leqslant \alpha_0$.
\end{example}

\begin{example}\label{ex_bert}
Take $p=1$. Suppose that there exists a positive constant $L$ such that $\mathcal{W}_{1}(\mathbb{Q}_{\mathbf{z}=\mathbf{z}'}, \mathbb{Q}_{\mathbf{z}=\mathbf{z}^*}) \leqslant L\|\mathbf{z}'-\mathbf{z}^*\|$, for all  $0 \leqslant \|\mathbf{z}'-\mathbf{z}\| \leqslant r_0$ and that Assumption~\ref{assumption_knn}.1 holds.

Following a line of reasoning that is parallel to that of the previous example, we also get
$$\mathcal{W}_1(\mathcal{R}_{1-\alpha}(\mathbb{Q}),\mathbb{Q}_{\widetilde{\Xi}}) = O(\alpha^{1/d_{\mathbf{z}}})$$ for all $0 < \alpha \leqslant \alpha_0$, with $\alpha_0:= \alpha(r_0)$.
\end{example}

Equation~\eqref{rho_finitesampleguarantee} and Examples~\ref{ex_guzi} and~\ref{ex_bert} reveal that our finite sample guarantee is affected by the \emph{curse of dimensionality}. Recently, powerful ideas to break this curse have been introduced in~\cite{Gao2020} under the standard Wasserstein-metric-based DRO scheme. In our setup, however, we also need distributional robustness against the (uncertain) error incurred when inferring conditional information from a sample of the true \emph{joint} distribution. This implies increasing the robustness budget in our approach by an amount linked to the term $\mathcal{W}_p(\mathcal{R}_{1-\alpha}(\mathbb{Q}),\mathbb{Q}_{\widetilde{\Xi}})$. Consequently, we might need stronger assumptions on the data-generating model to break the dependence of this term with the dimension of the feature vector and thus extend the ideas in~\cite{Gao2020} to the realm of conditional stochastic optimization.

Now we state the conditions under which the sequence of problems $\left({\rm P}_{(\alpha_N,\widetilde{\rho}_N)}\right)$, $N \rightarrow \infty$, is asymptotically consistent.

\begin{lemma}[Convergence of conditional distributions]\label{lemma_convergence_trimmed_distributions_uni}
Suppose that the support $\Xi$ of the true joint distribution $\mathbb{Q}$ is compact and that Assumptions \ref{assump_hellinger} and \ref{assumption_knn}.1 hold.
Take $(\alpha_N, \widetilde{\rho}_N)$ such that  $\alpha_N
\rightarrow 0$, $\frac{N\alpha_N^2}{\log(N)}\rightarrow \infty$,
and $\widetilde{\rho}_N \downarrow \underline{\epsilon}^p_{N\alpha_N}$, where $\underline{\epsilon}_{N\alpha_N}$ is the minimum transportation budget as in Definition~\ref{MTB}.
%
Then, we have that
$$\mathcal{W}_p( Q^N_{\widetilde{\Xi}}, \mathbb{Q}_{\widetilde{\Xi}}) \rightarrow 0 \enskip a.s.$$
where $Q^N_{\widetilde{\Xi}}$ is any distribution from the ambiguity set $ \widehat{\mathcal{U}}_N(\alpha_N, \widetilde{\rho}_N)$.
\end{lemma}

\proof{Proof.}
First, we need to provide conditions under which $\mathcal{W}_p\left(\mathcal{R}_{1-\alpha}(\widehat{\mathbb{Q}}_N),\mathbb{Q}_{\widetilde{\Xi}} \right) \rightarrow 0$ a.s.
Since $\Xi$ is compact and $\mathcal{W}_{p-1}\left(\mathcal{R}_{1-\alpha}(\widehat{\mathbb{Q}}_N),\mathbb{Q}_{\widetilde{\Xi}} \right) \leqslant \mathcal{W}_p\left(\mathcal{R}_{1-\alpha}(\widehat{\mathbb{Q}}_N),\mathbb{Q}_{\widetilde{\Xi}} \right)$, we can take $p > d/2$ and $\alpha_N$ such that $\frac{N\alpha_N^2}{\log(N)}\rightarrow \infty$, so that the probabilities \eqref{ineq_measure_conc_uni} becomes summable over $N$ for any arbitrarily small $\epsilon$. In this way, we can choose a sequence $\beta_N \in \left(0,1\right)$, $N \in \mathbb{N}$, such that $\sum_{N=1}^{\infty}{\beta_N} < \infty$ and $\lim_{N\rightarrow \infty}{\epsilon_{N,p,\alpha_N}(\beta_N)}\rightarrow 0$. With this choice, we have
     \begin{align*}
     \mathbb{Q}^{\infty} \bigg[\lim_{N\rightarrow \infty} \mathcal{W}_p\left(\mathcal{R}_{1-\alpha_N}(\widehat{\mathbb{Q}}_N),\mathbb{Q}_{\widetilde{\Xi}} \right)  &- \mathcal{W}_p\left(\mathcal{R}_{1-\alpha_N}(\mathbb{Q}),\mathbb{Q}_{\widetilde{\Xi}}\right) =0\bigg] \\
     &=\mathbb{Q}^{\infty} \left[\lim_{N\rightarrow \infty} \mathcal{W}_p\left(\mathcal{R}_{1-\alpha_N}(\widehat{\mathbb{Q}}_N),\mathbb{Q}_{\widetilde{\Xi}} \right)=0   \right] = 1
      \end{align*}
because $\mathcal{W}_p\left(\mathcal{R}_{1-\alpha_N}(\mathbb{Q}),\mathbb{Q}_{\widetilde{\Xi}}\right) = O\left(\alpha_N^{2/pd_{\mathbf{z}}}\right) \rightarrow 0$ for $\alpha_N \rightarrow 0$.

Since, $\mathbb{Q}_{\widetilde{\Xi}} \in \mathcal{R}_{1-\alpha_N}(\widehat{\mathbb{Q}}_N)$ a.s. in the limit and, by definition, $\mathbb{Q}_{\widetilde{\Xi}}(\widetilde{\Xi}) = 1$, we have that $\mathbb{Q}_{\widetilde{\Xi}} \in \widehat{\mathcal{U}}_N(\alpha_N, \widetilde{\rho}_N)$ for $N$ sufficiently large, with both $\alpha_N, \widetilde{\rho}_N \rightarrow 0$.

For its part, because $Q^N_{\widetilde{\Xi}} \in \widehat{\mathcal{U}}_N(\alpha_N, \widetilde{\rho}_N)$, this means that $\mathcal{W}_p\left(\mathcal{R}_{1-\alpha_N}(\widehat{\mathbb{Q}}_N),Q^N_{\widetilde{\Xi}}\right) \leqslant \widetilde{\rho}_N$. Take $N$ large enough, set $\widetilde{\rho}_N$ arbitrarily close to $\underline{\epsilon}^p_{N\alpha_N}$ and notice that $\widehat{\mathcal{U}}_N(\alpha_N, \underline{\epsilon}^p_{N\alpha_N})$ boils down to one single probability measure, the one made up of the $N\alpha_N$ data points of  $\widehat{\mathbb{Q}}_N$ that are the closest to $\widetilde{\Xi}$. In addition, we have $\underline{\epsilon}^p_{N\alpha_N} \rightarrow 0$ with probability one. To see this, take $K:= \lceil N\alpha_N\rceil$ and note that  $$\underline{\epsilon}^p_{N\alpha_N} \leqslant \textrm{dist}
    (\widehat{\boldsymbol{\xi}}_{K:N},
    \widetilde{\Xi}) \rightarrow \|\mathbf{\widehat{z}}_{K:N}-\mathbf{z}^*\| \rightarrow 0$$
    almost surely provided that $\alpha_N \rightarrow 0$ (see   \cite[Lemmas 2.2 and 2.3]{Biau2015}), where $\mathbf{\widehat{z}}_{K:N}$ is the $\mathbf{z}$-component of the $K$-\emph{th} nearest neighbor to $\mathbf{z}^*$ after reordering the data sample $\{\widehat{\boldsymbol{\xi}}_i:=(\mathbf{\widehat{z}}_i, \mathbf{\widehat{y}}_i)\}_{i=1}^{N}$ in terms of $\|\mathbf{\widehat{z}}_{i}-\mathbf{z}^*\|$ only.

    Therefore, it must hold that $\mathcal{W}_p( Q^N_{\widetilde{\Xi}}, \mathbb{Q}_{\widetilde{\Xi}}) \rightarrow 0$ a.s.
\qed

\begin{remark}
The compactness of the support set $\Xi$ is assumed here just to simplify the proof. In fact, in Appendix~\ref{case_alpha_0_supplementary}, we use results from nearest neighbors to show that the convergence of conditional distributions can be attained under the less restrictive condition $\frac{N\alpha_N}{\log(N)}\rightarrow \infty$ even in some cases for which the uncertainty $\mathbf{y}$ and the feature vector $\mathbf{z}$ are unbounded. In addition, we also make use of those results to demonstrate that distributionally robust versions of some local nonparametric predictive methods, such as Nadaraya-Watson kernel regression and $K$-nearest neighbors, naturally emerge from our approach.
\end{remark}

\begin{remark}
The convergence of conditional distributions allows us to establish an asymptotic consistency result analogous to that of Theorem~\ref{theorem_consistency_trimmed_alphapos}, by simply replacing ``Theorem~\ref{finite_sample_alpha_pos}'', ``$\widetilde{\rho}_N$'' and ``Lemma~\ref{lemma_convergence_distributions_alpha_pos}'' with ``Theorem~\ref{finite_sample_theorem_conditional_trimmed_uni}'', ``$(\alpha_N, \widetilde{\rho}_N)$'' and ``Lemma~\ref{lemma_convergence_trimmed_distributions_uni}'', respectively.
\end{remark}

\begin{remark}\label{remark:alpha0}
Suppose that the event $\widetilde{\Xi}$ on which we condition problem~\eqref{conditional_expectation_problem} is given by $\widetilde{\Xi}:=\{ \boldsymbol{\xi}=(\mathbf{z}_1,\mathbf{z}_2,\mathbf{y}) \in \Xi \; : \; \mathbf{z}_1 = \mathbf{z}^*_1,\ \mathbf{z}_2 \in \mathcal{Z}_2\}$, with $\mathbb{Q}(\widetilde{\Xi})=0$ and $\mathbb{P}(\mathbf{z}_2 \in \mathcal{Z}_2) > 0$. Let $\mathbb{Q}_{\mathcal{Z}_2}$ be the probability measure of $(\mathbf{z}_1, \mathbf{y})$ conditional on $\mathbf{z}_2 \in \mathcal{Z}_2$. If we have that there is $\widetilde{C}>0$ and $r_0 >0$ such that $\mathbb{P}(\|\mathbf{z}^*_1-\mathbf{z}_1 \| \leqslant r) \geqslant \widetilde{C} r^{d_{\mathbf{z}_1}}$, for all
$0 < r \leqslant r_0$, and that $\mathbb{Q}_{\mathcal{Z}_2}$ satisfies the smoothness condition invoked in either Theorem~\ref{finite_sample_theorem_conditional_trimmed_uni}, Example~\ref{ex_guzi} or Example~\ref{ex_bert}, then the analysis in this section extends to that type of event by setting $\alpha(r)= \widetilde{C} r^{d_{\mathbf{z}_1}} \cdot \mathbb{P}(\mathbf{z}_2 \in \mathcal{Z}_2)$ and noticing that  $\mathbb{Q}_{B(\mathbf{z}_1^*, r) \times \mathcal{Z}_2 \times \Xi_{\mathbf{y}}} \in \mathcal{R}_{1-\alpha(r)}(\mathbb{Q})$, $0 < r \leqslant r_0$, where $\mathbb{Q}_{B(\mathbf{z}_1^*, r) \times \mathcal{Z}_2 \times \Xi_{\mathbf{y}}}$ is the probability measure of $(\mathbf{z}_1,\mathbf{z}_2,\mathbf{y})$ conditional on $(\mathbf{z}_1,\mathbf{z}_2,\mathbf{y}) \in B(\mathbf{z}_1^*, r)\times \mathcal{Z}_2 \times\Xi_{\mathbf{y}}$.
\end{remark}

\section{Numerical Experiments}\label{numerics}
%

The following simulation experiments are designed to provide numerical evidence on the performance of the DRO framework with side information that we propose, with respect to other methods available in the technical literature. Here we only consider the case $\alpha = 0$, while additional numerical experiments for the case $\alpha > 0$ can be found in Appendix~\ref{NE_alpha_pos}.

%
To numerically illustrate the setting $\mathbb{Q}(\widetilde{\Xi}) = \alpha = 0$,  we consider two well-known problems, namely, the (single-item) newsvendor problem and the portfolio allocation problem, both posed in the form $\inf_{\mathbf{x}\in X}
\mathbb{E}_{\mathbb{Q}} \left[f(\mathbf{x},\boldsymbol{\xi})\mid \boldsymbol{\xi} \in \widetilde{\Xi}\right]$ to allow for side information. We compare four data-driven approaches to address the solution to these two problems: Our approach, i.e., problem $\text{P}_{(\alpha_N,\widetilde{\rho}_N)}$ with $\alpha_N=K_N/N$, which we denote ``DROTRIMM''; a Sample Average Approximation method based on a local predictive technique, in particular, the $K_N$ nearest neighbors, which we refer to as ``KNN'' (see \cite{Bertsimas2019b} for further details); this very same local predictive method followed by a standard Wasserstein-metric-based DRO approach to robustify it, as suggested in \cite[Section 5]{Bertsimas2019}, which we call ``KNNDRO''; and the robustified KNN method~\eqref{bertsimas_sturt}, also proposed in \cite{Bertsimas2019}, which we term ``KNNROBUST.'' 
%
We clarify that KNNDRO uses the $K$ nearest neighbors projected onto the set $\widetilde{\Xi}$ as the nominal ``empirical'' distribution that is used as the center of the Wasserstein ball in \cite{MohajerinEsfahani2018}.

We also note that the newsvendor problem and the portfolio optimization problem are structurally different if seen from the lens of  the standard Wasserstein-metric-based DRO approach.
Indeed, the newsvendor problem  features an objective function with a Lipschitz  constant with respect to the uncertainty that is independent of the decision $\mathbf{x}$. Consequently,
as per   \cite[Remark 6.7]{MohajerinEsfahani2018}, KNNDRO renders the same minimizer for this problem  as that of KNN whenever the  support set $\widetilde{\Xi}$ is equal  to the whole space. This is, in contrast, not true for the portfolio allocation problem, which has an objective function with a Lipschitz  constant with regard to the uncertainty that depends on the decision $\mathbf{x}$.

In all the numerical experiments, we take the $p$-norm with $p=1$ and, accordingly, we use the Wasserstein distance of order 1. Thus, all the optimization problems that we solve are linear programs. We consider a series of different values for the size $N$ of the sample data. Unless stated otherwise in the text, for each $N$, we choose  as the number of neighbors, $K_N$, the value $\lfloor N/\log (N+1) \rfloor$, where $\lfloor \cdot\rfloor$ stands for the floor function. Nevertheless, for the portfolio allocation problem, we also test the values $\lfloor N^{0.9}\rfloor$ and $\lfloor\sqrt{N} \rfloor$ 
 to assess the impact of the number of neighbors on the out-of-sample performance of  the four methods we compare.

We estimate $\mathbf{x}^* \in \argmin_{\mathbf{x}\in X}\;  \mathbb{E}_{\mathbb{Q}_{\widetilde{\Xi}}}
  \left[f(\mathbf{x},\boldsymbol{\xi})\right]$
  and $J^* = \mathbb{E}_{\mathbb{Q}_{\widetilde{\Xi}}} \left[f(\mathbf{x^*},\boldsymbol{\xi})\right]$ using a discrete proxy of the true conditional distribution $\mathbb{Q}_{\widetilde{\Xi}}$. In the newsvendor problem, this proxy is made up of 1085 data points, resulting from applying the KNN method (with the logarithmic rule) to 10\,000 samples from the true data-generating joint distribution. In the  portfolio optimization problem, we have an explicit form of $\mathbb{Q}_{\widetilde{\Xi}}$, which we utilize to directly construct a 10\,000-data-point approximation.
  To compare
  the four data-driven approaches under consideration, we use two performance metrics, specifically, the \emph{out-of-sample performance} of the data-driven solution  and its
  \emph{out-of-sample disappointment}. The former is given by $J=\mathbb{E}_{\mathbb{Q}_{\widetilde{\Xi}}} \left[f(\widehat{\mathbf{x}}_{N}^{m},\boldsymbol{\xi})\right]$, while the latter is calculated as $J -\widehat{J}_{N}^{m}$, where $m = \{\textrm{KNN, KNNROBUST, DROTRIMM, KNNDRO}\}$ and $\widehat{J}_{N}^{m}$ is the objective function value yielded by the data-driven optimization problem solved by method $m$. We note that a negative out-of-sample disappointment represents a favorable outcome.

Since $\mathbb{E}_{\mathbb{Q}_{\widetilde{\Xi}}} \left[f(\widehat{\mathbf{x}}_{N}^{m},\boldsymbol{\xi})\right]$ and $\widehat{J}_{N}^{m}$ are functions of the sample data, we conduct  a certain number of runs (400 for the newsvendor problem and  200 for the portfolio optimization problem)  for every $N$, each run with an independent sample of size $N$.
This way we can get (visual) estimates  of the out-of-sample performance and disappointment for several values of the sample size $N$ for different independent runs. These estimates are illustrated in the form
of box plots in a series of figures, where the dotted black horizontal line corresponds to either the optimal solution $\mathbf{x}^*$ (only in the newsvendor problem) or to its associated optimal cost $J^*$ with complete information.

As is customary in practice, we use a data-driven procedure to tune the robustness parameter of each method. In particular, for a desired value of reliability $1-\beta \in (0,1)$ (in our numerical experiments, we set $\beta$ to 0.15), and for each method $j$, where  $j = \{\textrm{KNNROBUST,\ KNNDRO,\ DROTRIMM}\}$, we aim for the value of the  robustness parameter for which the estimate of the objective value  $\widehat{J}_{N}^{j}$  given by method $j$
   provides an upper $(1-\beta)$-confidence bound  on the out-of-sample performance
  of its respective optimal solution (see Equation~\eqref{bounds_finite_sample_guarantee_1}),
  while delivering the best out-of-sample performance. As the optimal robustness parameter is unknown and depends on the available data sample, we need to  derive an estimator $param^{\beta,j}_N$  that is also a function of the training data. We construct $param^{\beta,j}_N$ and the corresponding reliability-driven solution  as follows:
    \begin{enumerate}
    \item  We generate $kboot$ resamples (with replacement) of size $N$, each playing the role of a different training set. In our experiments we set $kboot = 50$.
    Moreover, we build a validation dataset determining the $K_{N_{val}}$-neighbors of the $N_{val}$  data points of the original sample of size $N$ that have not been used to form the training set.

    \item For each resample $k=1,\ldots,kboot$ and each candidate value for $param$,  we compute a  solution by method $j$ with parameter  $param$ on the $k$-th resample.
    The resulting optimal decision  is denoted as $\widehat{x}^{j,k}_N(param)$ and its corresponding objective value as $\widehat{J}^{j,k}_N(param)$. Thereafter, we calculate the out-of-sample performance $J(\widehat{x}^{j,k}_N(param))$  of the data-driven solution $\widehat{x}^{j,k}_N(param)$ over the  validation set.
\item From among the candidate values for $param$ such that $\widehat{J}^{j,k}_N(param)$
exceeds the value $J(\widehat{x}^{j,k}_N(param))$
in at least
$(1-\beta)\times kboot$ different resamples, we take as $param^{\beta,j}_N$  the one yielding the best out-of-sample performance averaged over the $kboot$  validation datasets.
\item Finally, we compute the solution given by method $j$ with parameter
$param^{\beta,j}_N$, $\widehat{x}^j_N:=\widehat{x}^{j}_N(param^{\beta,j}_N)$
and the respective certificate
$\widehat{J}^j_N:=\widehat{J}^{j}_N(param^{\beta,j}_N)$.
\end{enumerate}

Recall that, in our approach, the robustness parameter $\widetilde{\rho}_{N}$ must be greater than or equal to the minimum transportation budget to the power of $p$, that is, $\underline{\varepsilon}^p_{N\alpha_N}$. Hence, if we decompose $\widetilde{\rho}_{N}$ as $\widetilde{\rho}_{N} = \underline{\varepsilon}^p_{N\alpha_N} + \Delta \widetilde{\rho}_N$, what one really needs to tune in DROTRIMM is the budget excess $\Delta \widetilde{\rho}_N$.
Furthermore,  for the same amount of budget $\Delta \widetilde{\rho}_N$, our approach will lead to more robust decisions $\mathbf{x}$ than KNNDRO,  because the worst-case distribution in KNNDRO is also feasible in DROTRIMM. Consequently, in practice, the tuning of one of these methods could guide the tuning of the other.



Lastly, all the  simulations have been run on a Linux-based server using up to 116 CPUs running in paralell, each clocking at 2.6 GHz with 4 GB of RAM. We have employed Gurobi 9.0 
under Pyomo 5.2 
to solve the associated linear programs.

\subsection{The single-item newsvendor problem}\label{num_exp:NV}
In this subsection, we deal with the popular single-item newsvendor problem,
which has received a lot of attention lately (see, for example,
\cite{Ban2019,Huber2019} and references therein). 
It is known that the solution to the single-item newsvendor problem is equivalent to that of a quantile regression problem, where the goal is to estimate the quantile  $b/(b+h)$ of the distribution of the uncertainty $y$, with $h$ and $b$ being the unit holding and backorder costs, respectively.

\begin{figure}[h!]
\centering
\subfloat[Heat map of the true joint distribution and kernel estimate of the true conditional density]{%
  \includegraphics[width=0.41\textwidth]{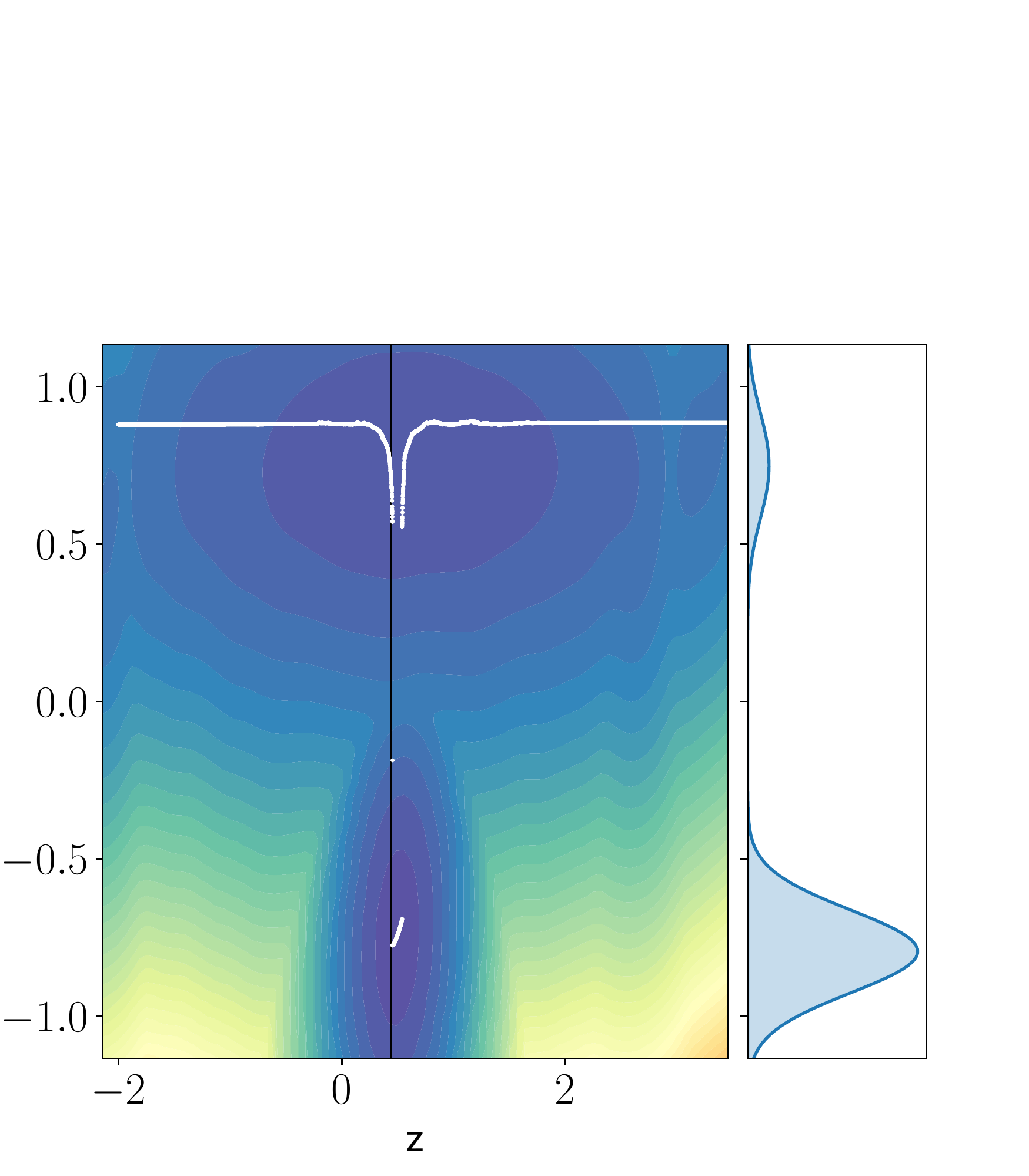}%
  \label{plot_distr_mixt_bivariate_asim}
}%
\subfloat[Optimal solution]{%
  \includegraphics[width=0.51\textwidth]{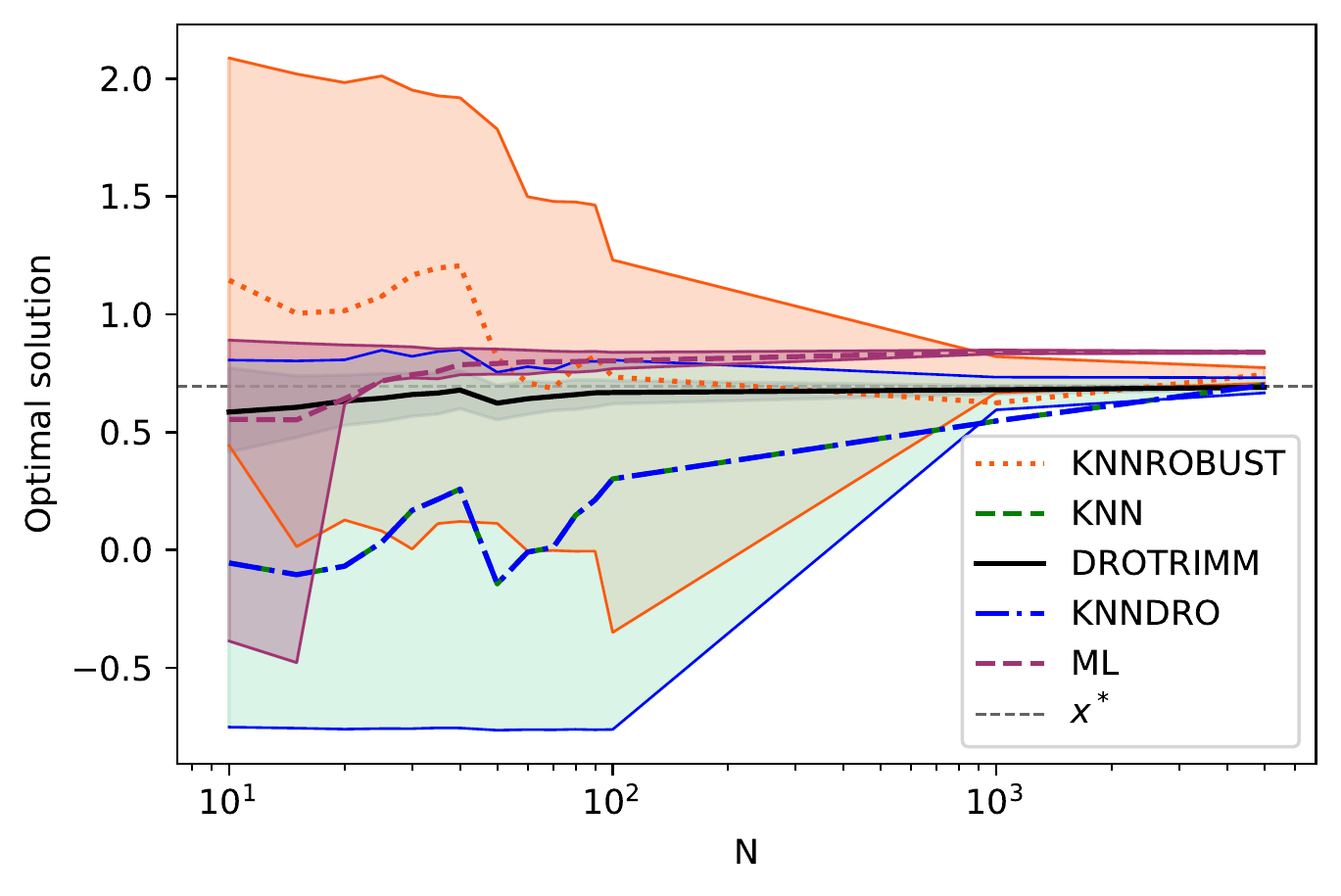}%
  \label{optimal_order_quantity1_univariate_ray_asim}
}

\subfloat[Out-of-sample disappointment]{%
  \includegraphics[width=0.47\textwidth]{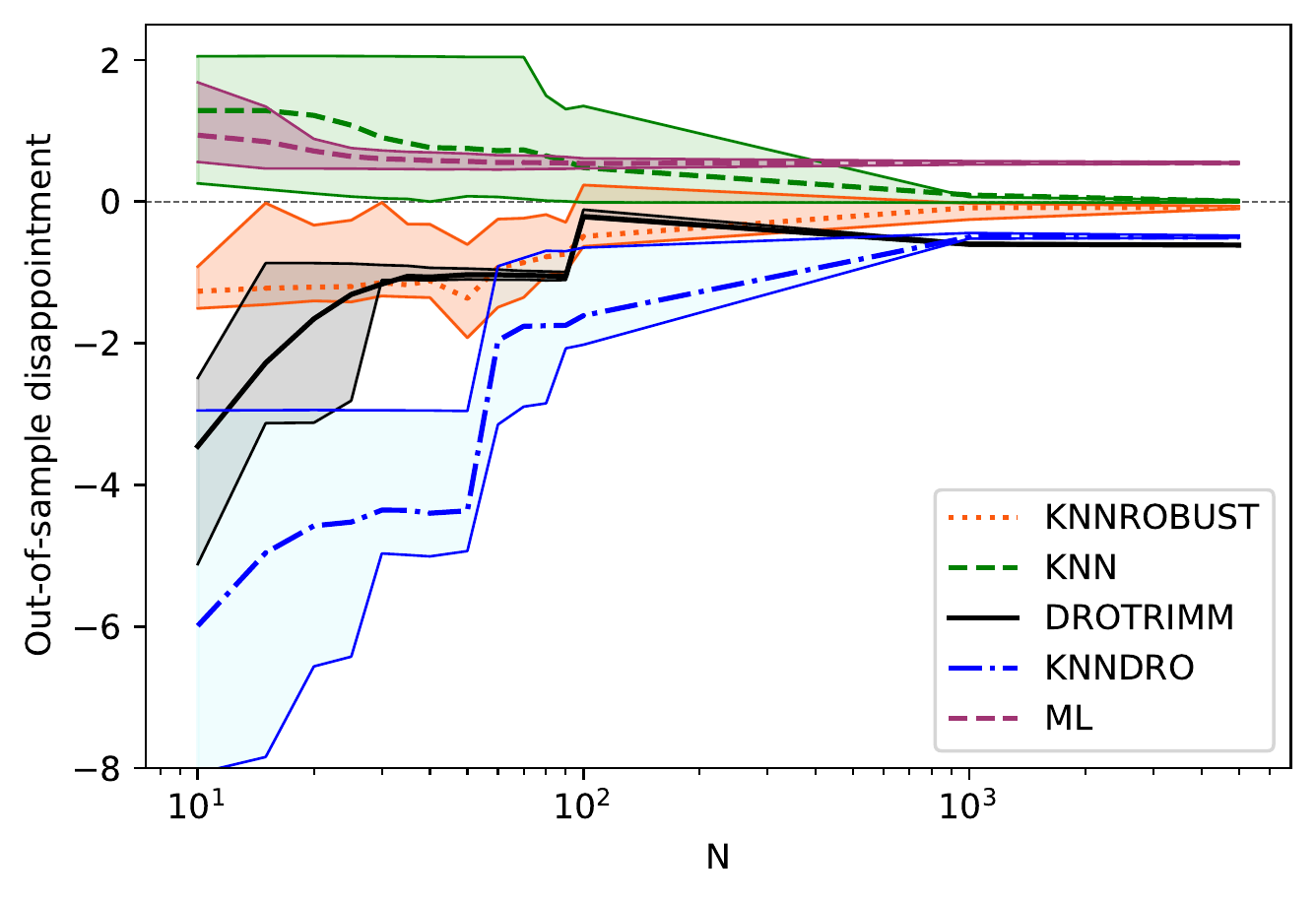}%
  \label{out-of-sample_univariate_ray_asim}
}%
\subfloat[Out-of-sample performance]{%
  \includegraphics[width=0.45\textwidth]{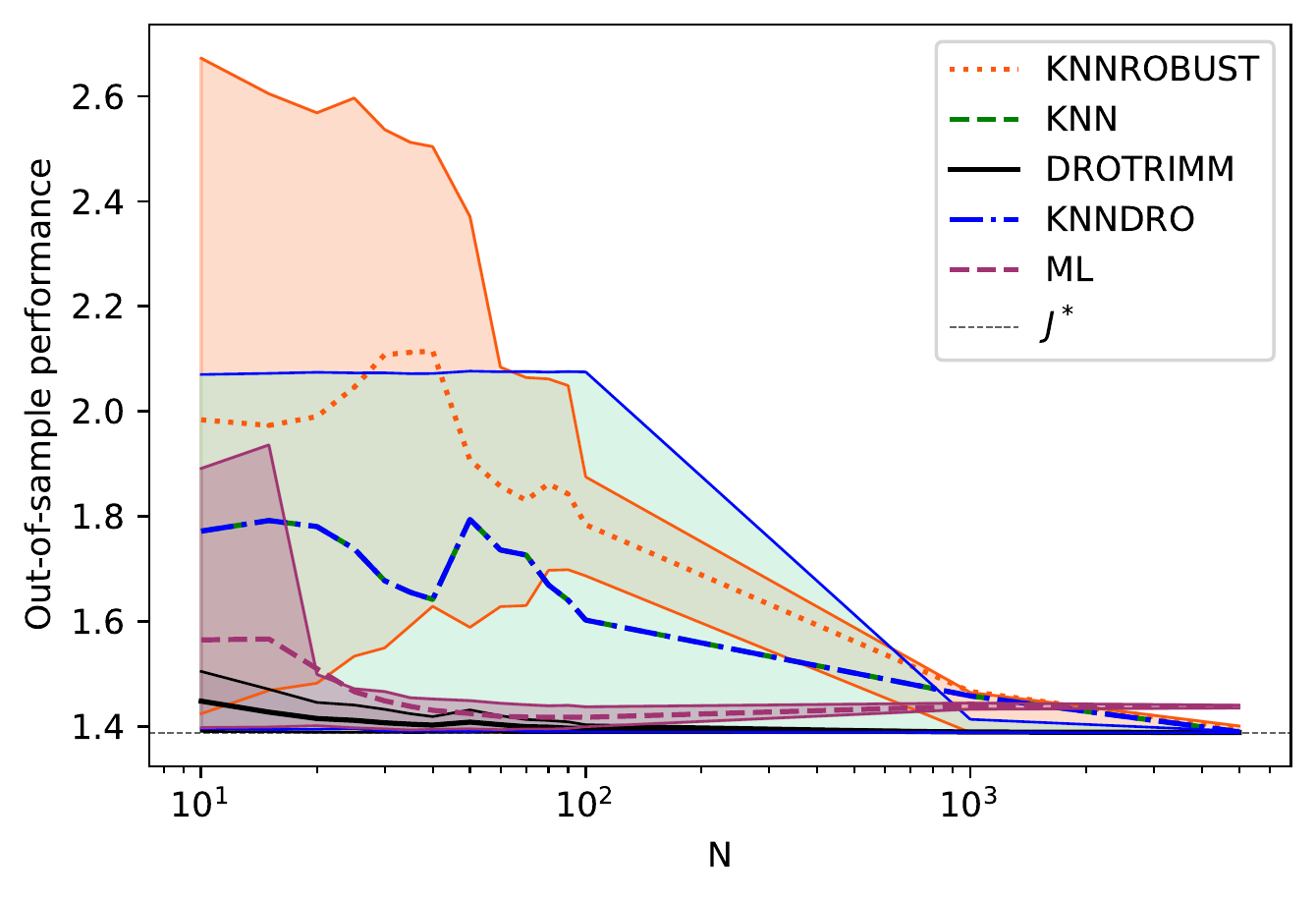}%
  \label{actual_expected_cost_univariate_ray_asim}
}

\vspace{5mm}

\caption{Newsvendor problem with features:   True distributions, quantile estimate and performance metrics}

\end{figure}

For the particular instance of this problem that we analyze next, we have considered $h=1$ and $b=10$. Furthemore, the true joint distribution of the data $\widehat{\boldsymbol{\xi}}_i:=(\widehat{z}_i,\widehat{y}_i)$,
$i=1,\ldots, N$ is assumed to follow a mixture (with equal weights) of two normal bivariate distributions with means $\mu_1=[0.6,0.75]^T, \; \mu_2=[0.5,-0.75]^T$ and covariance matrices $\Sigma_1=\begin{bmatrix}
0.5 & 0 \\
0 & 0.01
\end{bmatrix},$
$\Sigma_2=\begin{bmatrix}
0.0001 & 0 \\
0 & 0.1
\end{bmatrix}$,
respectively. Therefore, the support set of this distribution is the whole space $\mathbb{R}^{d_z+d_y}$, with $d_z =d_y=1$. In addition, we consider as $\mathcal{Z}$ the singleton  $\{z^* = 0.44\}$, with $\widetilde{\Xi}$ being the real line $\mathbb{R}$ as a result.  Figure~(\ref{plot_distr_mixt_bivariate_asim}) shows a heat  map  of the  true  joint  distribution, together with a kernel estimate of the probability density function of the random variable $y$ conditional on $z^*$.  Moreover, the white dotted curve in the figure corresponds to the optimal order quantity as a function of the feature $z$. Note that this curve is highly nonlinear around the context $z^*$. Also, the demand may be negative, which, in the context of the newsvendor problem, can be interpreted as items being returned to the stores due to, for example, some quality defect.
The set of candidate values from which the robustness parameters in methods KNNROBUST, KNNDRO and DROTRIMM have been selected is the discrete set composed of the thirty linearly spaced numbers between 0 and 2.  We also consider the machine learning algorithm proposed in \cite{Ban2019}, which was especially designed for the newsvendor problem with features. In this algorithm, a polynomial mapping between the optimal order quantity (i.e., the optimal quantile) and the covariates is presumed. The degree of the polynomial, up to the fourth degree, is tuned using the bootstrapping procedure described above. We denote this approach as ML from ``Machine Learning''.

Figures  \eqref{optimal_order_quantity1_univariate_ray_asim},  \eqref{out-of-sample_univariate_ray_asim},  and  \eqref{actual_expected_cost_univariate_ray_asim}  illustrate the  box  plots  corresponding  to  the  quantile estimators (i.e., the optimal solution of the problem),  the  out-of-sample disappointment and the out-of-sample performance delivered by each of the considered data-driven approaches for various sample sizes and runs, in that order. The  shaded color
areas have been obtained by joining the 15$th$ and 85$th$ percentiles of the box plots, while the associated bold colored lines link their means. The true optimal quantile (with complete information) and its out-of-sample performance are also depicted in Figures~\eqref{out-of-sample_univariate_ray_asim}  and  \eqref{optimal_order_quantity1_univariate_ray_asim}, respectively, using black dotted lines.

 Interestingly,  whereas the quantile estimators provided by DROTRIMM, KNNDRO and KNNROBUST  all lead to negative out-of-sample disappoinment in general,  KNNDRO and KNNROBUST exhibit substantially worse out-of-sample performance both in expectation and volatility. Recall that KNNDRO delivers the same solutions provided by KNN for this problem. Its behavior is, therefore, influenced by the bias introduced by the $K$-nearest neighbors estimation, which is particularly notorious for small-size samples in this case, given the shape of the true conditional density, see Figure~(\ref{plot_distr_mixt_bivariate_asim}). Actually, for some runs, the $K$-nearest neighbors, and hence KNNDRO, lead to negative quantile estimates, while the true one is positive and greater than 0.5. By construction, both KNNDRO and KNNROBUST are mainly affected by the estimation error of the conditional probability distribution incurred by the local predictive method. On the contrary, our approach DROTRIMM offers a natural protection against this error and a richer spectrum of data-driven solutions. Indeed, DROTRIMM is able to identify solutions that lead to a better out-of-sample performance with a negative out-of-sample disappointment.

Finally, both ML and DROTRIMM exhibit a notorious stable behavior against the randomness of the sample. The order quantity provided by the former, however, does not converge to the true optimal one, because the relationship between the true optimal order and the feature $z$ is far from being polynomial.  Note that ML is a \emph{global} method that seeks to learn the optimal order quantity  for \emph{all} possible contexts by using a polynomial up to the fourth degree. However, the (true) optimal order curve (that is, the white line in Figure~\ref{plot_distr_mixt_bivariate_asim}) is highly nonlinear within a neighborhood of the context $z^* = 0.44$, but practically constant outside of it.


%
\subsection{Portfolio optimization}\label{case:portf}
We consider next an instance of the portfolio optimization problem that is based on that used in \cite{BertsimasMCCORD2018} and
\cite{bertsimas2017bootstrap}. The instance corresponds to a single-stage  portfolio  optimization  problem  in  which  we
  wish  to  find  an allocation of a fixed budget to six different assets. Thus,
  $\mathbf{x}\in \mathbb{R}_+^6$ denotes the decision variable vector, that is, the asset allocations, and their uncertain  return  is represented by $\mathbf{y}\in \mathbb{R}^6$. In practice,
  these uncertain returns may be influenced by a set of features.
  First, the decision maker observes auxiliary covariates  and later,  selects the portfolio.
  We consider three different covariates that can potentially impact the returns and that we denote as $\mathbf{z} =(z_1,z_2,z_3)$.
The decision maker wishes to leverage this side information to improve  his/her decision-making process in which the goal is
 to maximize the expected value of the return
 while minimizing
 the  conditional  value  at  risk  (CVar)  of  the
  portfolio, that is,  the risk that the
loss $(-\langle \mathbf{x},\mathbf{y} \rangle)^+:=\max(-\langle \mathbf{x},\mathbf{y} \rangle,0)$ is large.
Using the reformulation of the CVar (see \cite{Rockafellar00optimizationof} and \cite{bertsimas2017bootstrap}) and introducing the auxiliary variable $\beta'$, the decision maker aims to solve the following  optimization problem given the value of the covariate $\mathbf{z}^*$($=(1000,0.01,5)$ in the numerical experiments):
  \begin{equation}
    \min_{(\mathbf{x},\beta') \in X }  \mathbb{E}\left[\beta'+\frac{1}{\delta}\left(-\langle \mathbf{x},\mathbf{y} \rangle -\beta' \right)^+-\lambda\langle \mathbf{x},\mathbf{y} \rangle
    \; \mid\; \mathbf{z}=\mathbf{z}^* \right]
  \end{equation}
  where  the feasible set of decision variables of the problem, that is, $X$ is equal to
$\{(\mathbf{x},\beta') \in \mathbb{R}_{+}^6\times \mathbb{R} :\; \sum_{j=1}^6 x_j=1   \}$. We set $\delta=0.5$ and $\lambda=0.1$ to simulate an investor with a moderate level of risk aversion. The parameter $\lambda \in \mathbb{R}_+$ serves to tradeoff between risk and return, and $\delta$ refers to the $(1-\delta)$-quantile of the loss distribution. We take the same marginal distributions for the covariates  as in Section 5.2 of \cite{bertsimas2017bootstrap}, i.e.,
$z_1 \leadsto  \mathcal{N}(1000,50)$, $z_2 \leadsto  \mathcal{N}(0.02,0.01)$ and $\log(z_3) \leadsto \mathcal{N}(0,1)$. Furthermore, we follow their approach to construct the joint true distribution of the covariates and the asset returns. In particular, we take
\begin{equation*}
    \mathbf{y}/(\mathbf{z}=(z_1,z_2,z_3)) \leadsto \mathcal{N}_6(\boldsymbol{\mu}+0.1\cdot(z_1-1000)\cdot \mathbf{v}_1+1000\cdot z_2\cdot \mathbf{v}_2+10\cdot \log (z_3+1)\cdot \mathbf{v}_3,\boldsymbol{\Sigma})
\end{equation*}
with $\mathbf{v}_1 = (1,1,1,1,1,1)^T$, $\mathbf{v}_2 = (4,1,1,1,1,1)^T$, $\mathbf{v}_3 = (1,1,1,1,1,1)^T$,  and with
$\boldsymbol{\mu}, \boldsymbol{\Sigma}^{1/2}$ given in \cite{bertsimas2017bootstrap,DROTRIMMINGS_github2021}.


Note that, unlike in \cite{bertsimas2017bootstrap}, not all the features affect equally all the asset returns.  Moreover, feature $z_3$ is log-normal and therefore, Assumption~\ref{light_tailed_assumption_joint} does not hold. Nonetheless, as we  show below,  DROTRIMM performs satisfactorily, which reveals that the conditions we derive in this paper to guarantee that our approach performs well are sufficient, but not necessary. Indeed, the condition $\mathbb{Q}_{\widetilde{\Xi}} \in \widehat{\mathcal{U}}_{N}(\alpha_N,\widetilde{\rho}_N)$ is not required to ensure performance guarantees \citep{Gao2020,Kuhn2019}.
For all the methods, we have standardized the covariates $\mathbf{z}$ and the asset returns $\mathbf{y}$ using their means and variances.  In all the simulations, the robustness parameter each method uses (i.e., $\varepsilon_N$ in KNNROBUST, the radius of the Wassertein ball, $\rho_N$, in KNNDRO, and the budget excess $\Delta \widetilde{\rho}_N$ in DROTRIMM) has been chosen from the discrete set  $\{b\cdot  10^c\;:\;  b \in \{0,\ldots, 9 \},\; c \in \{-2,-1,0\}\}$, following the above data-driven procedure.

\begin{figure}
\centering
\subfloat[$K_N=\lfloor N/(\log(N+1)) \rfloor$]{%
  \includegraphics[scale=0.5]{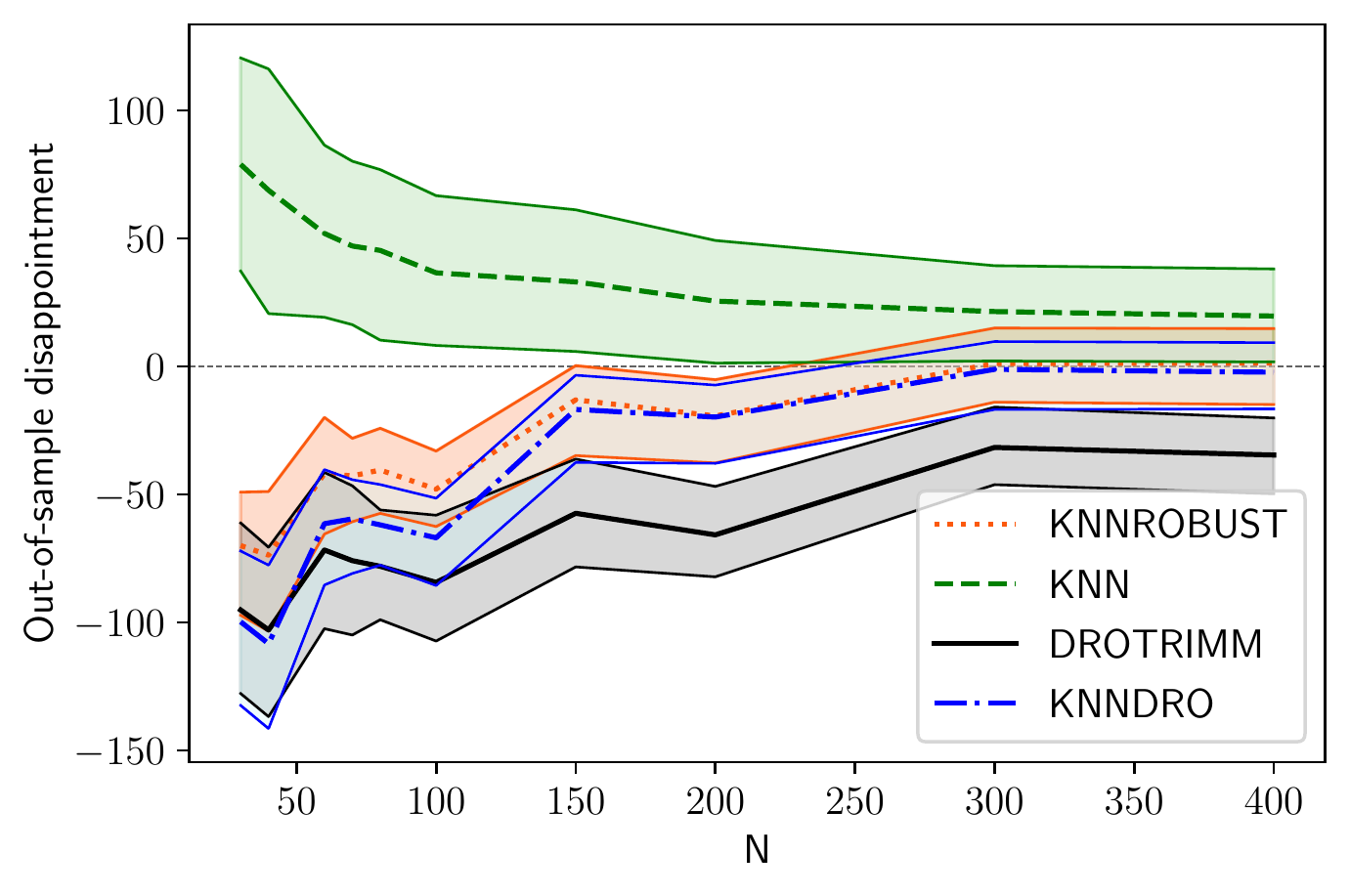}%
  \label{out-of-sample_portfolio_variante1.pdf}
}%
\subfloat[$K_N=\lfloor N/(\log(N+1)) \rfloor$]{%
  \includegraphics[scale=0.5]{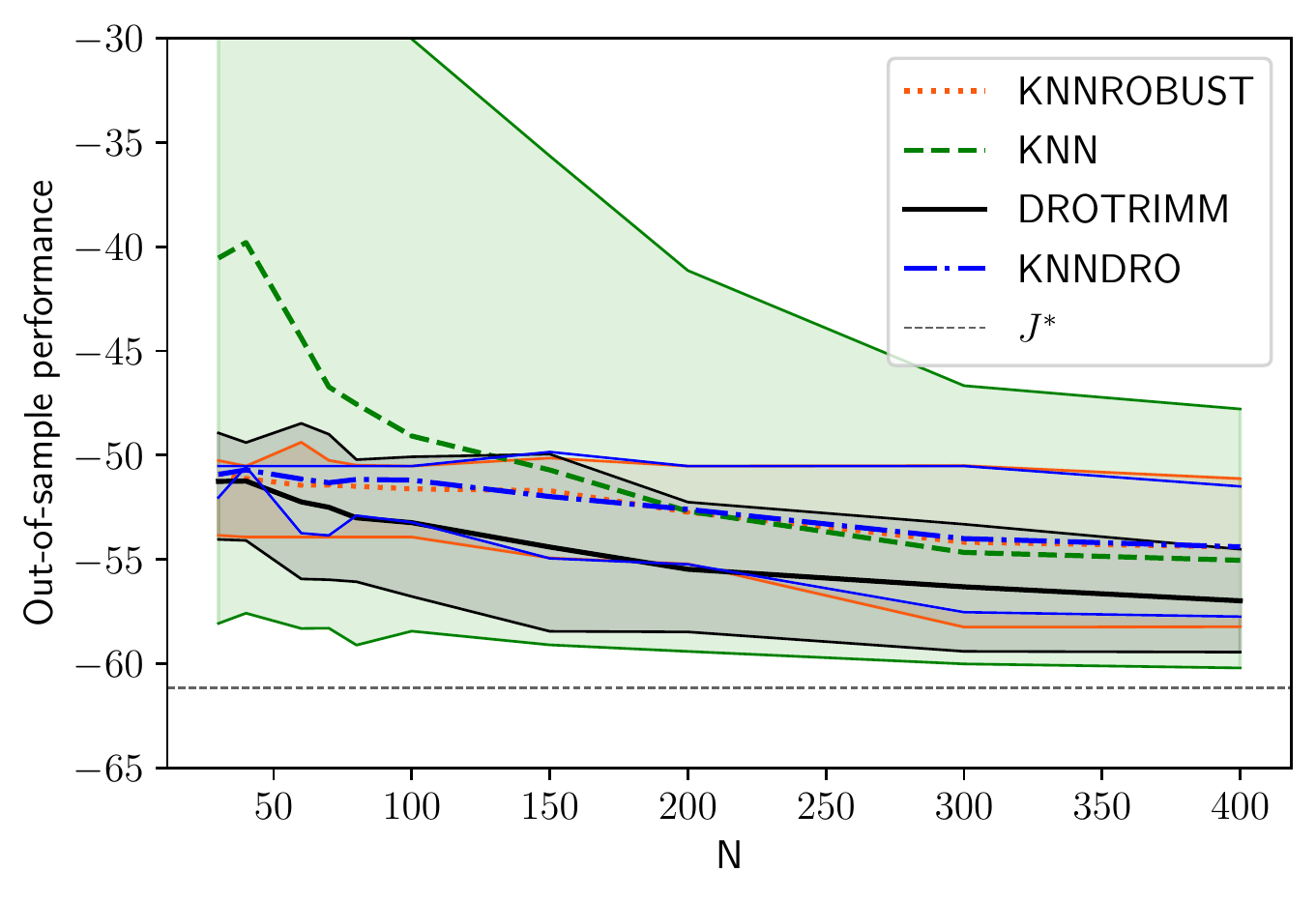}%
  \label{actual_expected_cost_portfolio_variante1.pdf}
}

\centering
\subfloat[$K_N=\lfloor N^{0.5}\rfloor$]{%
  \includegraphics[scale=0.5]{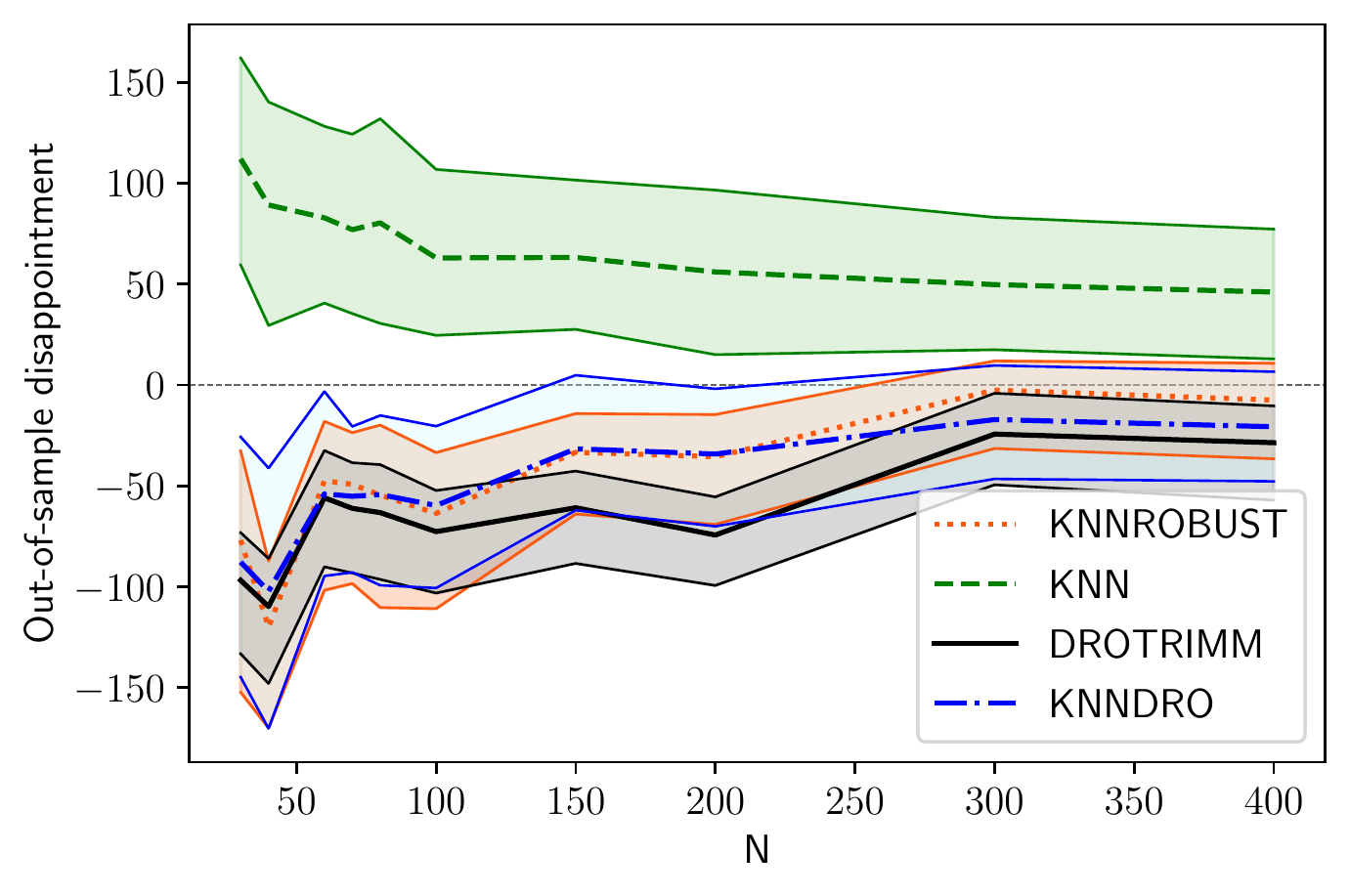}%
  \label{out-of-sample_portfolio_raizn_variante1.pdf}
}%
\subfloat[$K_N=\lfloor N^{0.5}\rfloor$]{%
  \includegraphics[scale=0.5]{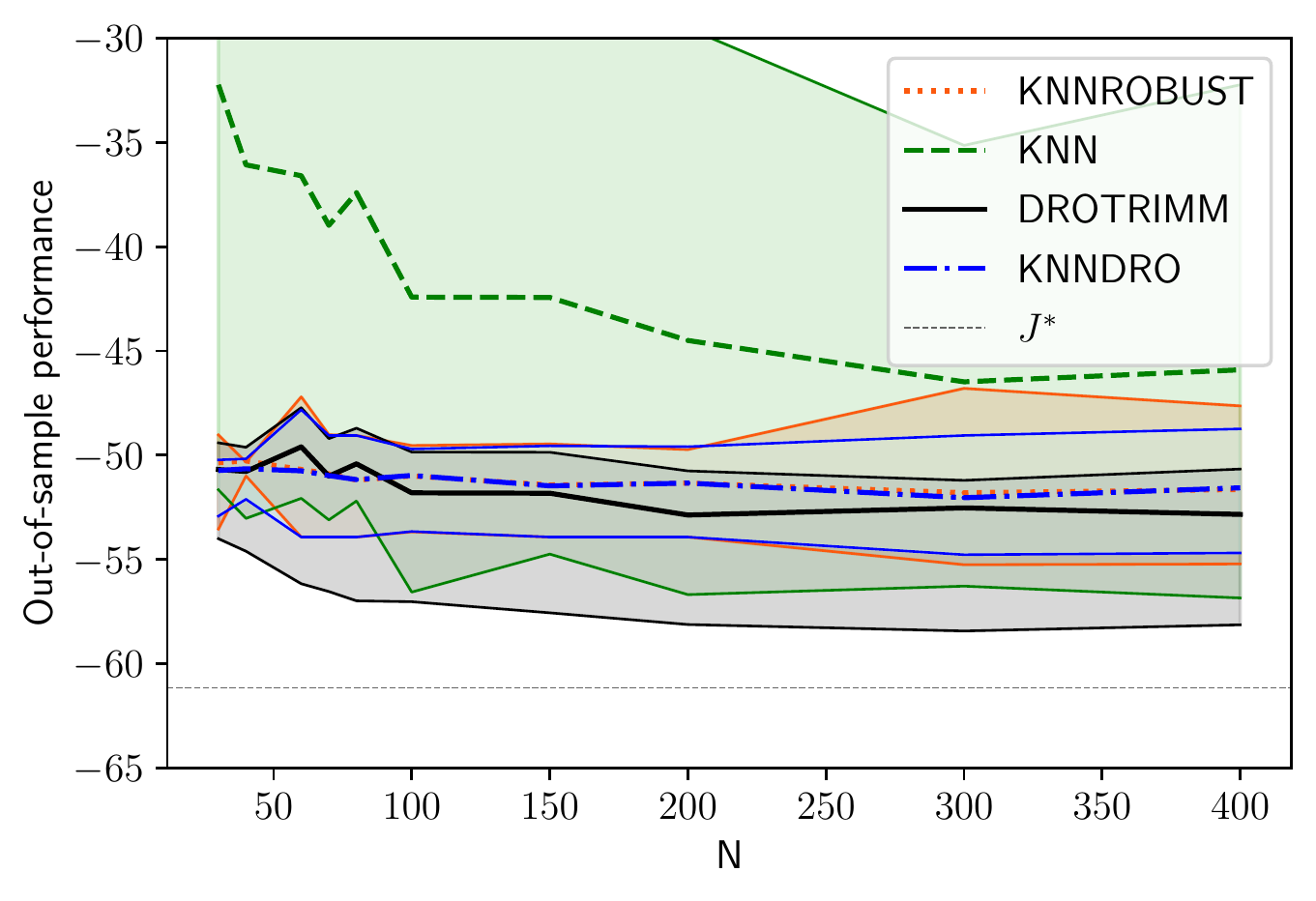}%
  \label{actual_expected_cost_portfolio_raizn_variante1.pdf}
}
\centering

\subfloat[$K_N=\lfloor N^{0.9}\rfloor$]{%
  \includegraphics[scale=0.5]{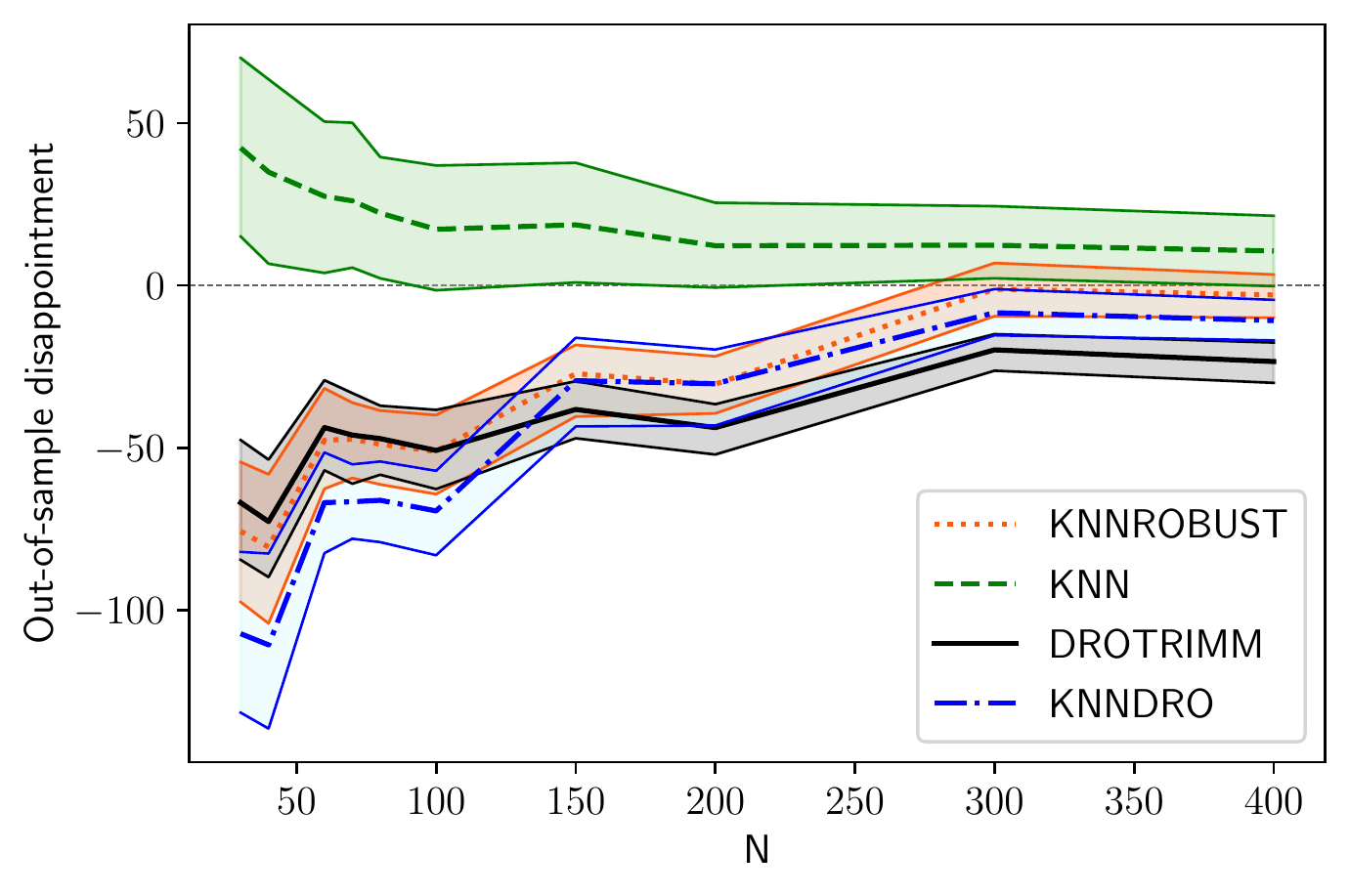}%
  \label{out-of-sample_portfolio_nnueve_variante1.pdf}
}%
\subfloat[$K_N=\lfloor N^{0.9}\rfloor$]{%
  \includegraphics[scale=0.5]{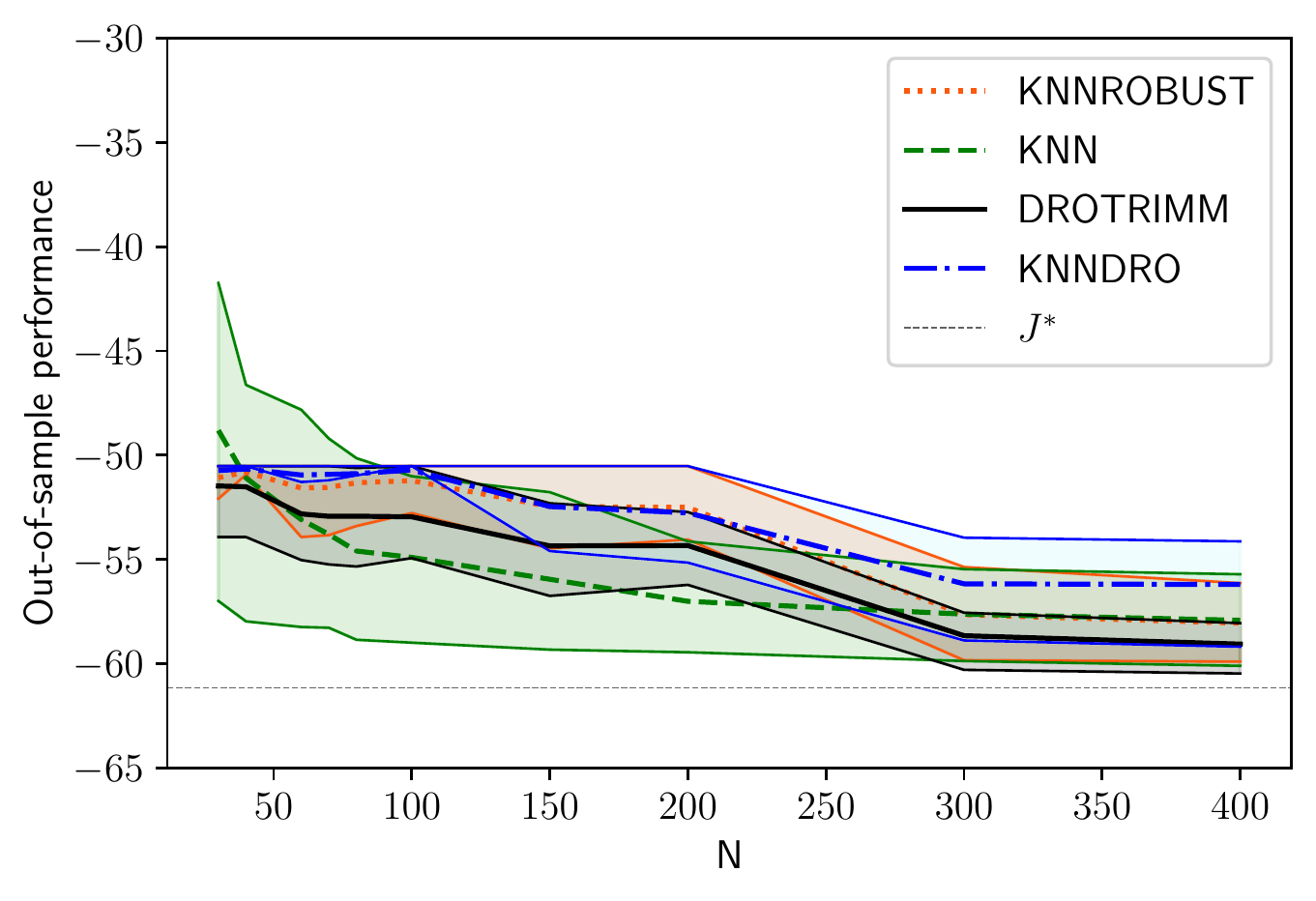}%
  \label{actual_expected_cost_portfolio_nnueve_variante1.pdf}
}

\centering

\subfloat[Tuned $K_N$]{%
  \includegraphics[scale=0.5]{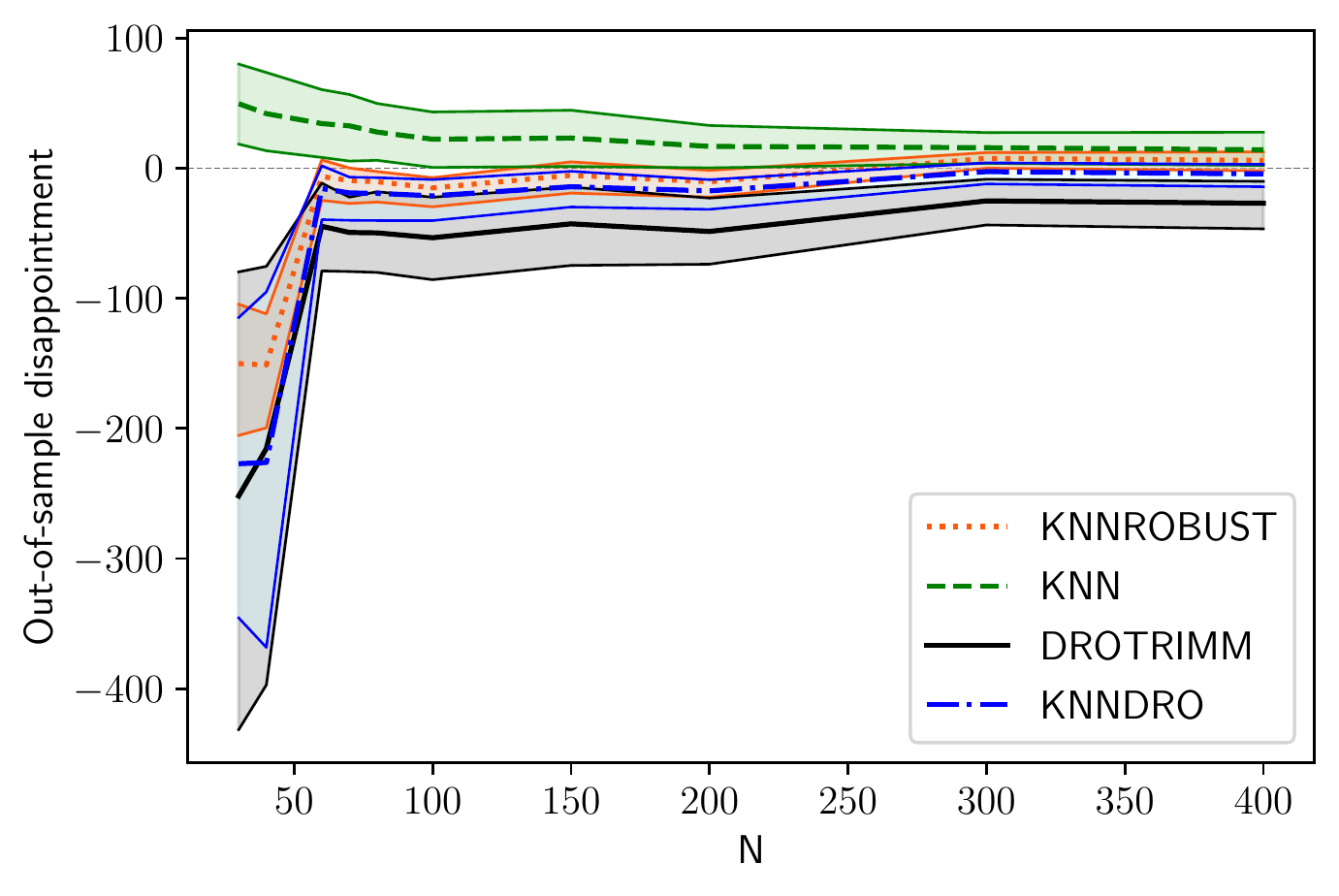}%
  \label{out-of-sample_portfolio_validalfarho.pdf}
}%
\subfloat[Tuned $K_N$]{%
  \includegraphics[scale=0.5]{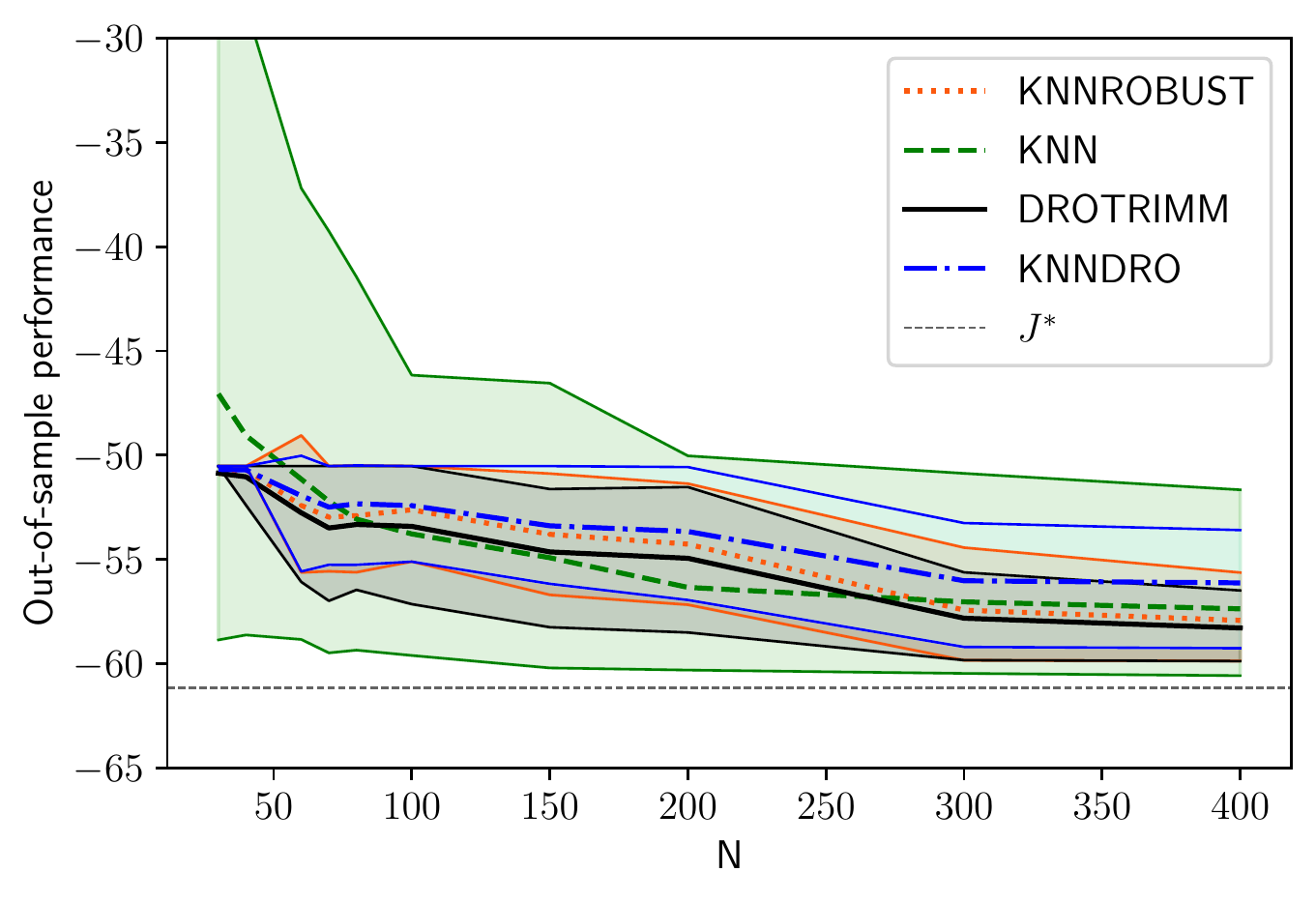}%
  \label{actual_expected_cost_portfolio_validalfarho.pdf}
}
\caption{  Portfolio problem  with features:  Performance metrics}\label{Results_alpha_0_performance}
\end{figure}

Similarly to the case of the single-item newsvendor problem, Figure~\ref{Results_alpha_0_performance} shows, for various sample sizes
and 200 runs, the box plots pertaining to  the out-of-sample disappointment and  performance associated with each of the considered data-driven approaches.
Each of the three pairs of subplots at the top of the figure has been obtained with a different rule to determine the number $K_N$ of nearest neighbors. Increasing this number seems to have a positive effect on the convergence speed of all the methods for this instance, although KNNROBUST (and KNNDRO to a lesser extent) has some trouble ensuring the desired reliability level, with the 85\% line above 0 for the largest values of $N$ we represent. In contrast, DROTRIMM manages to keep the disappointment negative. This is, in addition, accompanied by an important improvement of the the out-of-sample performance (in line with the criterion for selecting the best portfolio that we have established). In fact, DROTRIMM produces boxplots that appear to be shifted downward, i.e., in the direction of better objective function values. On the other hand, the KNN method substantially improves its performance by employing a larger number of neighbors. However, it is way too optimistic in any case.

The results shown in the pair of subplots at the bottom of Figure~\ref{Results_alpha_0_performance} correspond to a number $K_N$ of neighbors that has been tuned jointly with the robustness parameter and for each method independently. For this purpose, we have selected the best value of $K_N$ for each approach from the discrete set $\{N^{0.1},N^{0.2},\ldots,N^{0.9}\}$ following the bootstrapping-based procedure previously described. The data-driven tuning of the number $K_N$ of neighbors appears not to have a major effect on the performance of the different methods, especially in comparative terms. 
We do observe that the out-of-sample performance of KNNROBUST and KNNDRO is slightly improved on average. This improvement in cost performance is, however, accompanied by an increase in the number of sample sizes for which these methods do not satisfy the  reliability requirement, particularly in the case of KNNROBUST and small sample sizes.

\begin{figure}
\centering
\subfloat{%
  \includegraphics[width=0.5\textwidth]{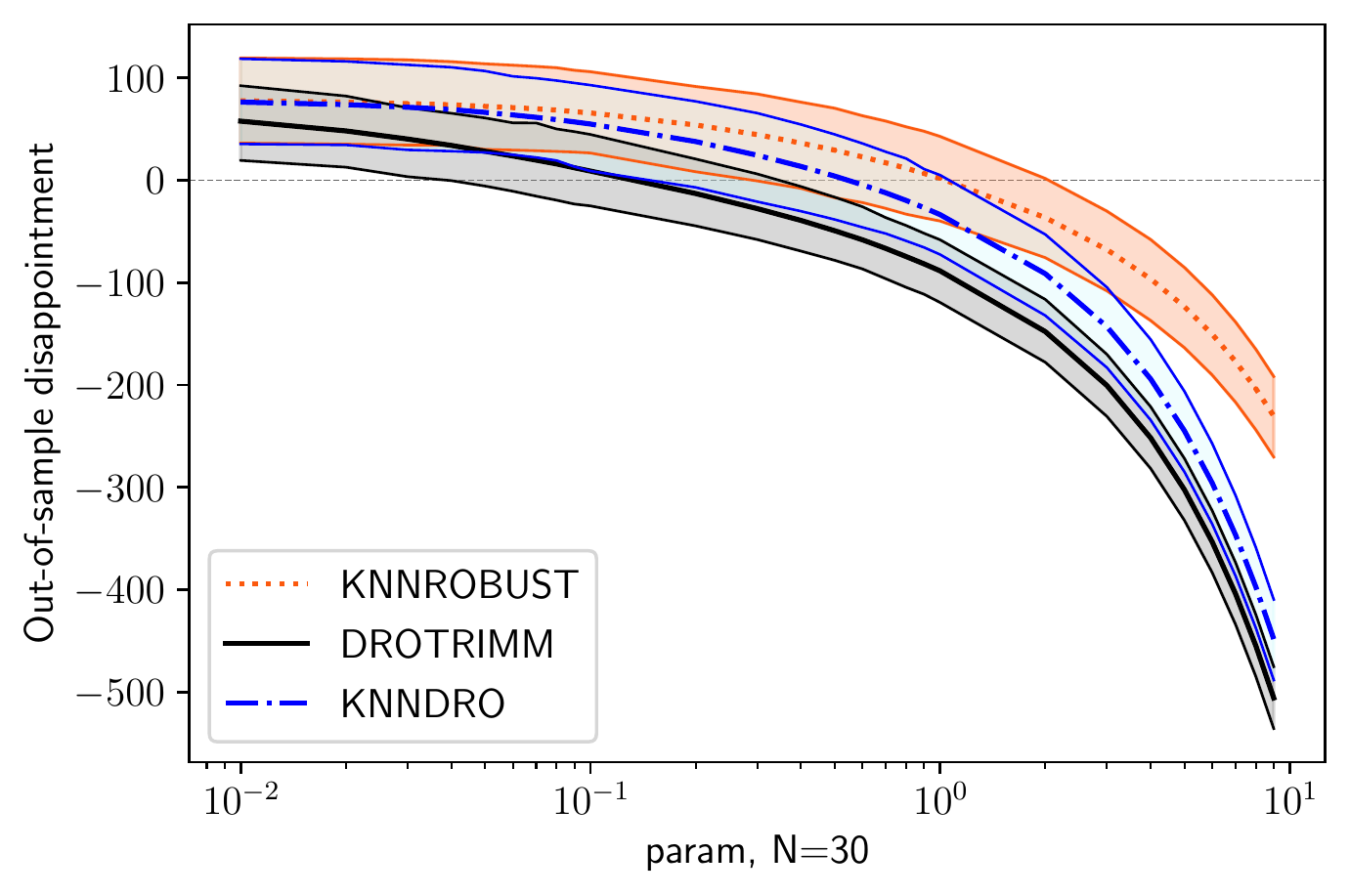}%
  \label{out-of-sample_portfolio_sens_30_DRO_variante1.pdf}
}%
\subfloat{%
  \includegraphics[width=0.5\textwidth]{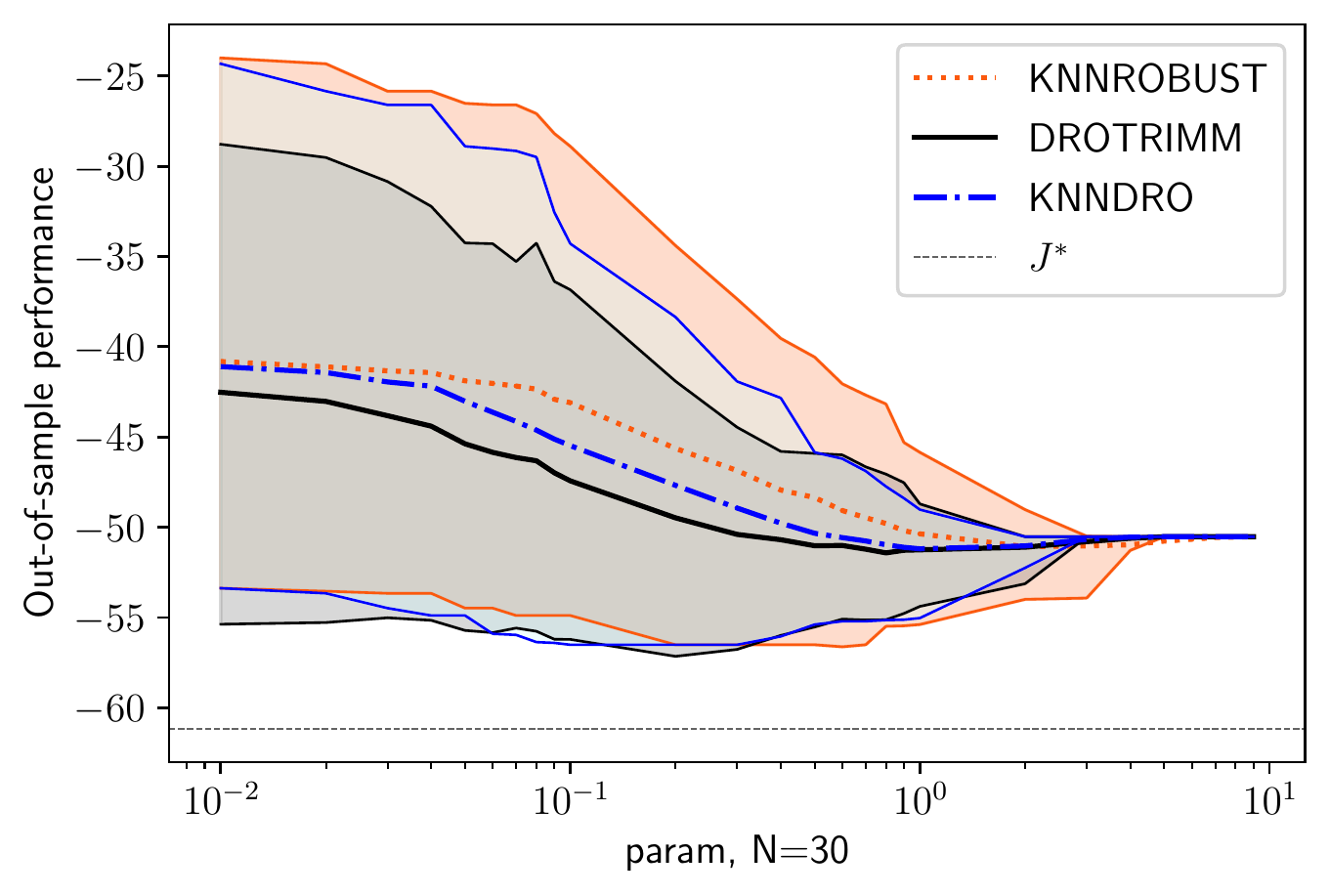}%
  \label{actualexpectedcost_portfolio_sens_30_DRO_variante1.pdf}
}

\vspace{5mm}
\centering
\subfloat{%
  \includegraphics[width=0.5\textwidth]{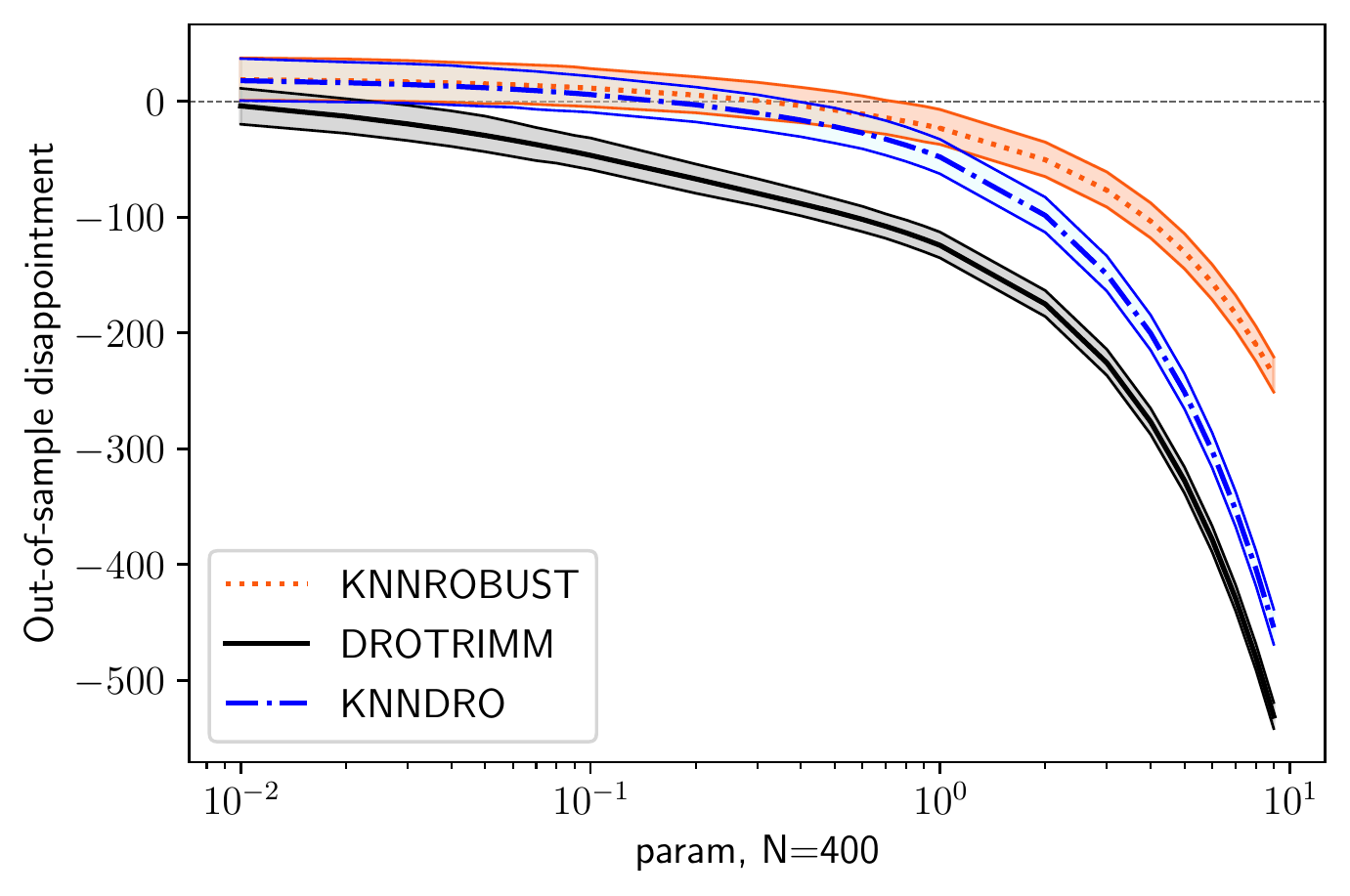}%
  \label{out-of-sample_portfolio_sens_400_DRO_variante1.pdf}
}%
\subfloat{%
  \includegraphics[width=0.5\textwidth]{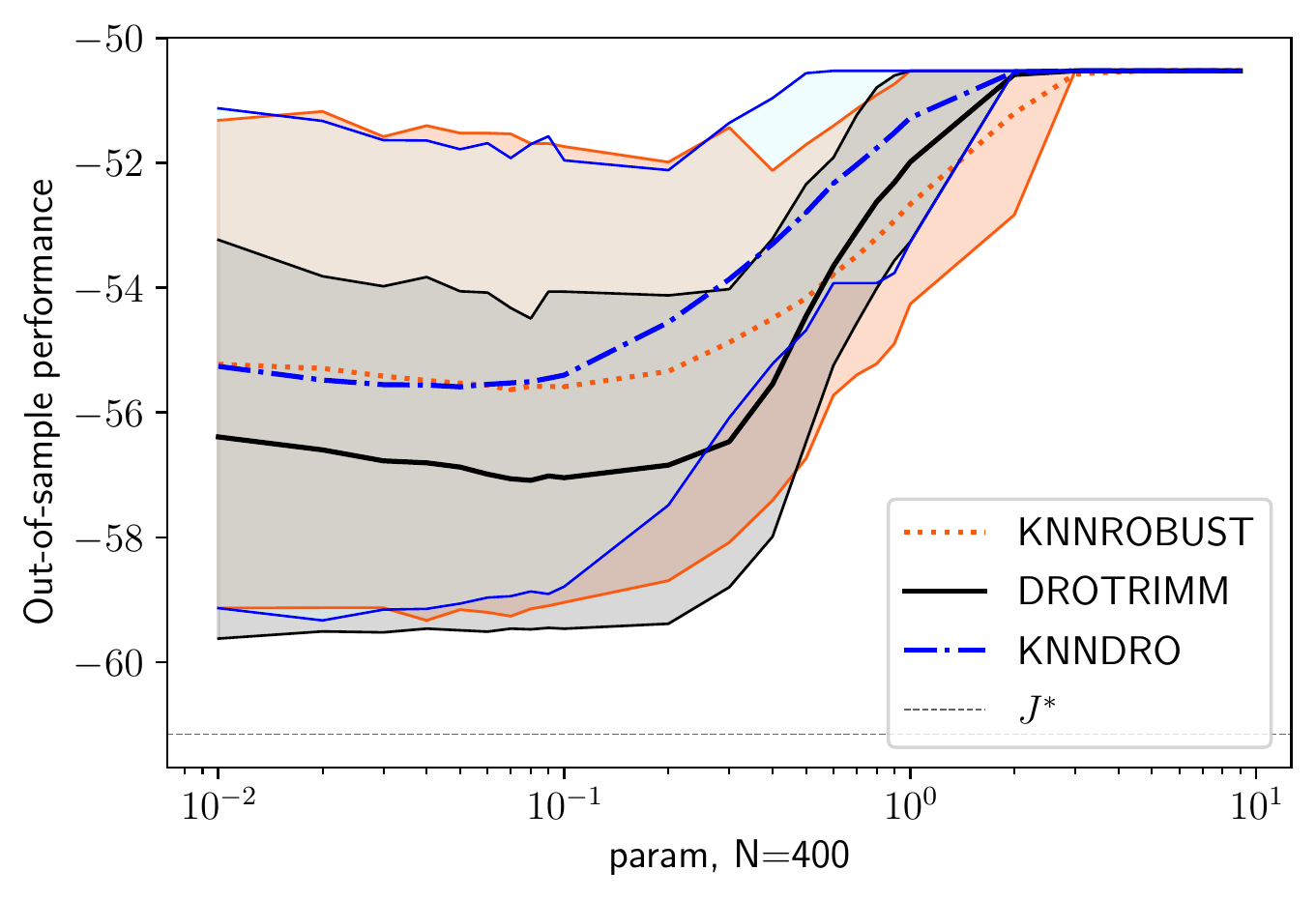}%
  \label{actualexpectedcost_portfolio_sens_400_DRO_variante1.pdf}
}

\vspace{5mm}

\caption{Impact of the robustness parameter with 200 training samples, $K_N=\lfloor N/(\log(N+1)) \rfloor$  and $\delta=0.5,\ \lambda=0.1$}\label{alpha_0_sens}
\end{figure}

To facilitate the analysis of the results shown in Figure~\ref{Results_alpha_0_performance}, we also provide Figure~\ref{alpha_0_sens}, which illustrates the (random) performance of the methods KNNROBUST, DROTRIMM and KNNDRO as a function of their respective robustness parameter, estimated over 200 independent runs. Again, the shaded areas cover the 15$th$ and 85$th$ percentiles, while the bold colored lines correspond to the average performance. The various plots are obtained for $N = 30$ and $N = 400$, with the number of neighbours given by the logarithmic rule. These plots are especially informative, because they are independent of the specific validation procedure used to tune the robustness parameters of the methods and thus, provide insight into the potential of each method to identify good solutions.
Note that the out-of-sample performance of all the three methods stabilizes around the same value as their respective robustness parameters grow large enough. This phenomenon is analogous to that discussed in   \cite[Section 7.1]{MohajerinEsfahani2018}. However, the value we observe here does not correspond to the ``equally weighted portfolio,'' because we have standardized the data on the asset returns. As a result, the ``robust portfolio'' that delivers this out-of-sample performance depends on and is solely driven by the standard deviations of the different assets. Very interestingly, DROTRIMM is able to uncover portfolios whose out-of-sample performance features a better mean-variance trade-off, in general. Furthermore, it requires a smaller value of the robustness parameter to guarantee reliability. All this is more evident (and useful) for the case $N = 400$, as we explain next. When $N =30$, all the  considered methods need large values of their robustness parameter to ensure reliability, so they all tend to operate close to the ``robust portfolio'' we mentioned above. DROTRIMM can certainly afford lower values of $\Delta \widetilde{\rho}$ in an attempt to improve performance, but this proves not to be that profitable for such a small sample size, for which the robust portfolio performs very well. As $N$ increases, the robust portfolio loses its appeal, since its performance gradually becomes comparatively worse. DROTRIMM is then able to identify portfolios that perform significantly better in expectation, while providing an estimate of their return such that the desired reliability is guaranteed. For their part, KNNDRO and KNNROBUST are also able to discover solutions with an actual average cost lower than that of the robust portfolio (albeit with a worse expectation and a higher variance than those given by DROTRIMM). However, they are more prone to overestimate their returns.

\begin{figure}[h!]

\centering

\subfloat{%
  \includegraphics[width=0.5\textwidth]{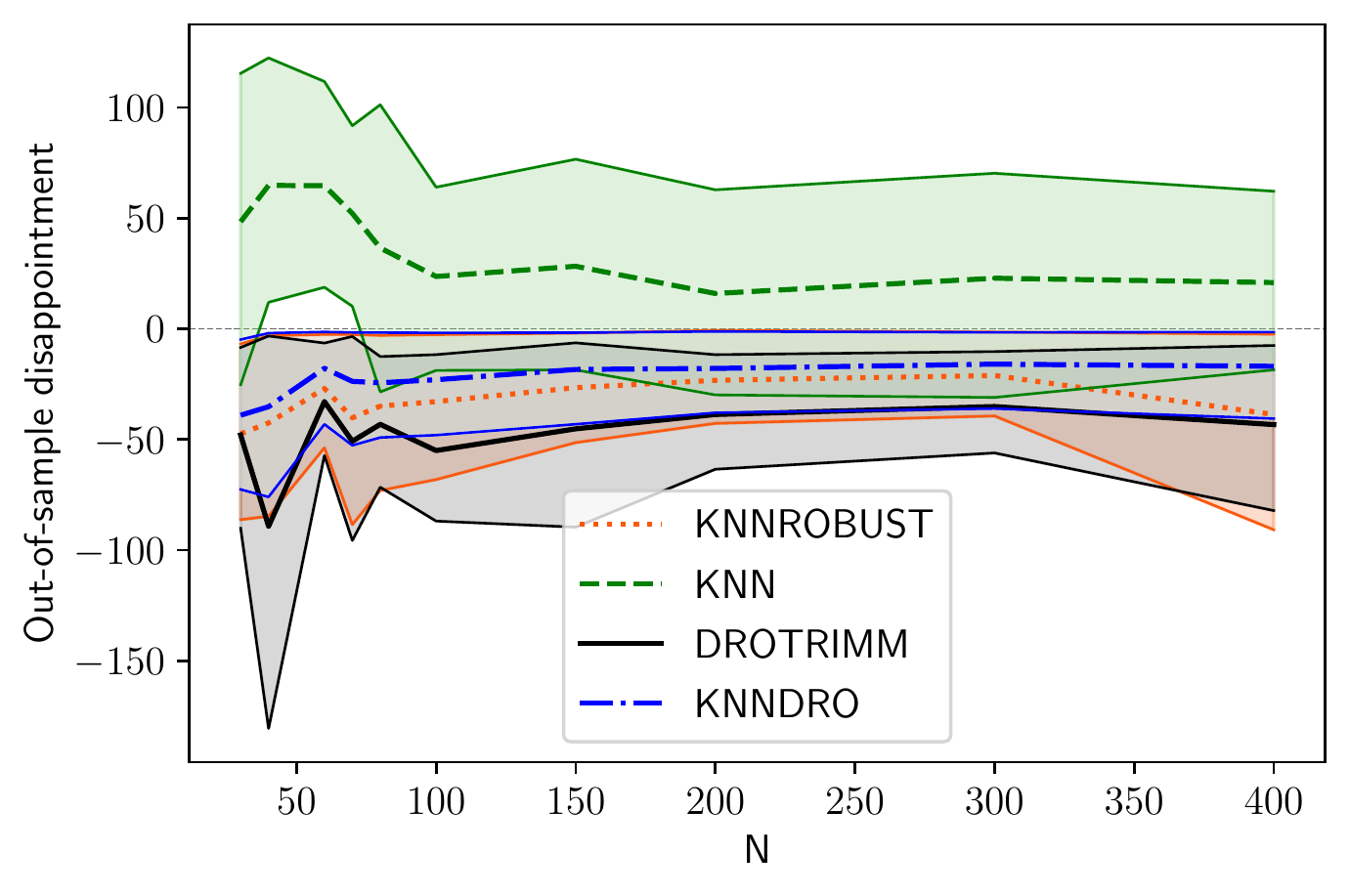}%
  \label{out-of-sample_portfolio_variandoz.pdf}
}%
\subfloat{%
  \includegraphics[width=0.5\textwidth]{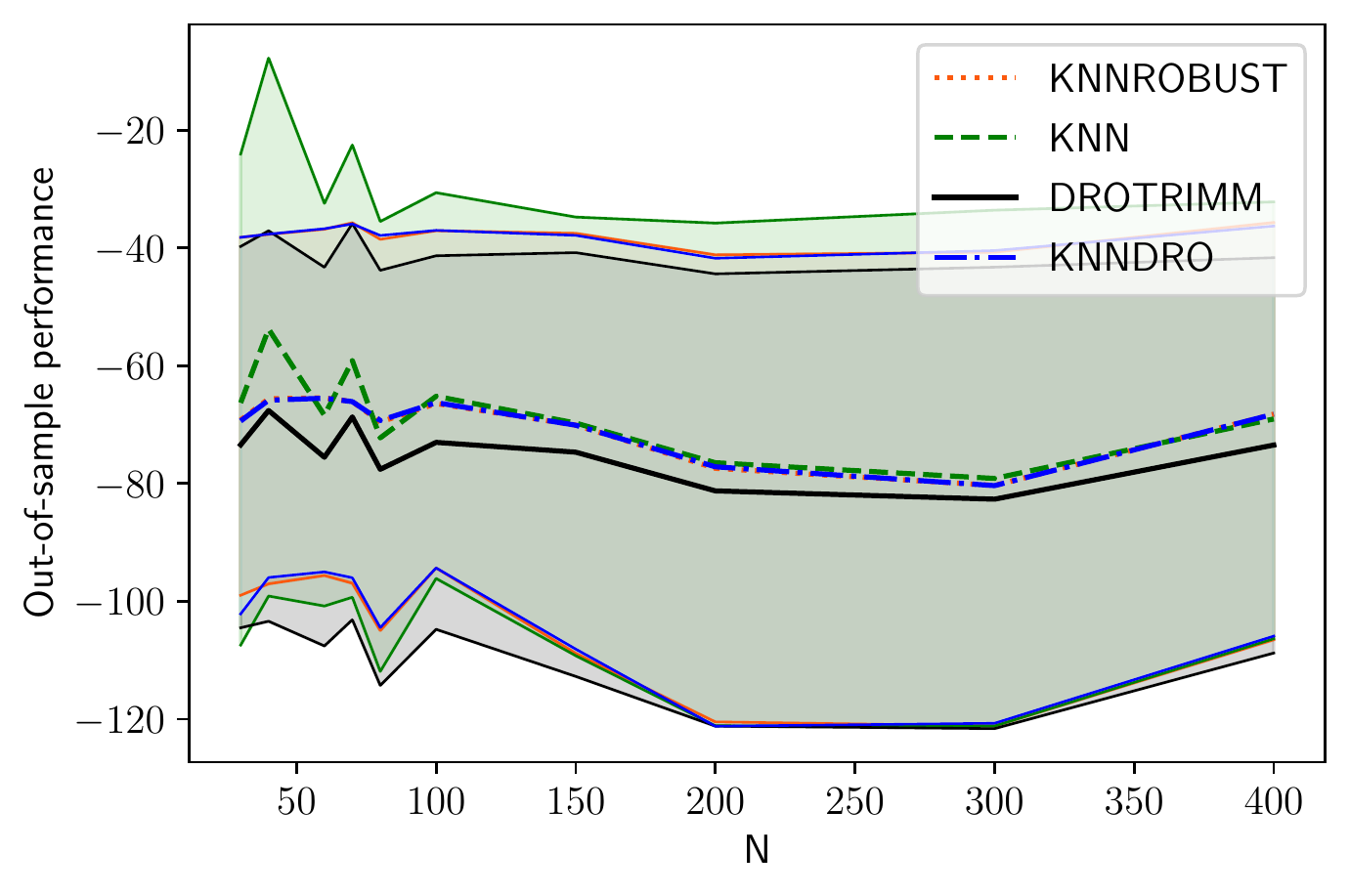}%
  \label{actual_expected_cost_portfolio_variandoz.pdf}
}

\vspace{5mm}

\caption{ Portfolio problem  with features: Varying context under an optimal selection of the robustness parameters,  $K_N=\lfloor N/(\log(N+1)) \rfloor$  and $\delta=0.5,\ \lambda=0.1$}\label{Results_alpha_pos_performance_perfectValidation}

\end{figure}

Finally, we study the behavior of the different methods under other contexts. For this, we consider several values of $N$, one random data sample for each $N$, and 200 different contexts $\mathbf{z}^*$ sampled from the marginal distributions of the features. The performance metrics (i.e., the out-of-sample disappointment and performance) are plotted in Figures~\ref{out-of-sample_portfolio_variandoz.pdf} and~\ref{actual_expected_cost_portfolio_variandoz.pdf}, respectively,  under an optimal selection of the robustness parameters (that is, for each method we use the value of the robustness parameter that, while ensuring a negative disappointment, delivers the best out-of-sample performance). We observe that DROTRIMM systematically performs better, with an actual cost averaged over the 200 contexts that is lower irrespective of the sample size. 

\section{Conclusions}\label{conclusions}
 In this paper, we have exploited the connection between probability trimmings and partial mass transportation to provide an easy, but powerful and novel way to extend the standard Wasserstein-metric-based DRO  to the case of \emph{conditional} stochastic programs. Our approach produces decisions that are distributionally robust against the uncertainty in the whole process of inferring the conditional probability measure of the random parameters from a finite sample coming from the true joint data-generating distribution. Through a series of numerical experiments built on the single-item newsvendor problem and a portfolio allocation problem, we have demonstrated that our method attains notably better out-of-sample performance than  some existing alternatives. We have supported these empirical findings with theoretical analysis, showing that our approach enjoys attractive performance guarantees.

%
%
%
 \begin{APPENDICES}
 \section{Proofs of theoretical results}\label{appen:proofs}

  This appendix compiles the proofs of some of the theoretical derivations that appear in the paper. The following technical results are needed to develop these proofs.
 \begin{definition}[Contamination of a distribution]\label{def:cont}
Given two probabilities $P,Q$ on $\mathbb{R}^d$, we say that $P$ is a $(1-\alpha)$-contaminated version of $Q$, if $P=\alpha Q+ (1-\alpha) R$, where $R$ is some probability. A $(1-\alpha)$-contamination neighbourhood of $Q$ is the set of all $(1-\alpha)$-contaminated versions of $Q$ and will be denoted as $\mathcal{F}_{1-\alpha}(Q)$.
\end{definition}

\begin{proposition}[Section 2.2. from \cite{Alvarez-esteban2012b} and p.18 in \cite{AgulloAntolin2018}]\label{Prop_trim_cont}
Let $P$, $Q$ be probabilities on $\mathbb{R}^d$ and $\alpha \in (0,1]$, then
\begin{equation}\label{equiv_trimmings}
Q \in \mathcal{R}_{1-\alpha}(P) \Longleftrightarrow P=\alpha Q+ (1-\alpha) R \Longleftrightarrow P \in \mathcal{F}_{1-\alpha}(Q)
\end{equation}
for some probability $R$. Moreover, if $D$ is a probability metric such that
$\mathcal{R}_{1-\alpha}(P)$ is closed for $D$ over an appropiate set of probability distributions, then \eqref{equiv_trimmings} is equivalent to
$D(Q,\mathcal{R}_{1-\alpha}(P))=0$.

\end{proposition}


\begin{remark}As  particular case, if we
 consider $D=\mathcal{W}_p$ over the set of probability distributions with finite $p$-th moment,  $\mathcal{P}_p$, we have that, if $P$, $Q \in \mathcal{P}_p$, then $Q\in \mathcal{R}_{1-\alpha}(P)$  if and only if $\mathcal{W}_p(Q,\mathcal{R}_{1-\alpha}(P)) = 0$.
\end{remark}

\begin{corollary}[Corollary 3.12 from \cite{AgulloAntolin2018}]\label{corollary312agullo}
Given two probabilities $P,Q \in \mathcal{P}_p (\mathbb{R}^d)$ and $\alpha \in (0,1)$, there exists $P_{1-\alpha}\in \mathcal{F}_{1-\alpha}(Q)$ such that
$P_{1-\alpha}=\alpha Q+ (1-\alpha) R_{1-\alpha}$ for some
$R_{1-\alpha} \in \mathcal{R}_{\alpha}(P)$ and
$\mathcal{W}_p(P,P_{1-\alpha})=\min_{R \in \mathcal{F}_{1-\alpha}(Q) } \mathcal{W}_p(P,R)$.

\end{corollary}

\begin{proposition}[Proposition 3.14 from \cite{AgulloAntolin2018}]\label{proposition314agullo}
Take $P,Q \in \mathcal{P}_p (\mathbb{R}^d)$. If $\alpha \in (0,1)$, then
$$\mathcal{W}_{p}^{p}\left(P, \mathcal{F}_{1-\alpha}(Q)\right)=\alpha \mathcal{W}_{p}^{p}\left(\mathcal{R}_{1-\alpha}(P), Q\right)$$
Moreover, if $\widehat{P}_{1-\alpha} \in \mathcal{R}_{1-\alpha}(P)$ is such that
$\mathcal{W}_{p}(\widehat{P}_{1-\alpha}, Q)=\mathcal{W}_{p}\left(\mathcal{R}_{1-\alpha}(P),
Q\right)$, then if we construct the probability measure $\widetilde{P}_{1-\alpha}=
\frac{1}{1-\alpha}\left(P-\alpha \widehat{P}_{1-\alpha}\right)$, we have that
$P_{1-\alpha}:=\alpha Q+(1-\alpha) \widetilde{P}_{1-\alpha} \in \mathcal{F}_{1-\alpha}(Q)$ and
$\mathcal{W}_{p}\left(P, P_{1-\alpha}\right)=\mathcal{W}_{p}\left(P, \mathcal{F}_{1-\alpha}(Q)\right)$.
\end{proposition}

 \subsection{Proof of Lemma~\ref{Lema_PMT}}\label{Lema_PMT_proof}

 We will prove the lemma by contradiction.
Suppose there are two different probability distributions $Q_{\widetilde{\Xi}}$ and $Q'_{\widetilde{\Xi}}$ such that $$D\left(\mathcal{R}_{1-\alpha}(Q), Q_{\widetilde{\Xi}}\right) = D(\mathcal{R}_{1-\alpha}(Q), Q'_{\widetilde{\Xi}}) = 0$$ and $Q_{\widetilde{\Xi}}(\widetilde{\Xi})= Q'_{\widetilde{\Xi}}(\widetilde{\Xi})=1$.

Because $D\left(\mathcal{R}_{1-\alpha}(Q), Q_{\widetilde{\Xi}}\right) = D(\mathcal{R}_{1-\alpha}(Q), Q'_{\widetilde{\Xi}}) = 0$, we know by Proposition~\ref{Prop_trim_cont} above
that $Q_{\widetilde{\Xi}}$, $Q'_{\widetilde{\Xi}} \in \mathcal{R}_{1-\alpha}(Q)$. Therefore, applying again Proposition~\ref{Prop_trim_cont}, we have
$$Q=\alpha Q_{\widetilde{\Xi}}+ (1-\alpha) R$$
$$Q=\alpha Q'_{\widetilde{\Xi}}+ (1-\alpha) R'$$
for some probabilities $R$ and $R'$ with $R(\widetilde{\Xi})=R'(\widetilde{\Xi}) = 0$.

Since, by hypothesis, $Q_{\widetilde{\Xi}}$ and $Q'_{\widetilde{\Xi}}$ are different, there must exist an event $A \subset \widetilde{\Xi}$ such that $Q_{\widetilde{\Xi}}(A) \neq Q'_{\widetilde{\Xi}}(A)$. We take that event and compute $Q(A)$ as follows:
$$Q(A)=\alpha Q_{\widetilde{\Xi}}(A)+ (1-\alpha) R(A) = \alpha Q'_{\widetilde{\Xi}}(A) + (1-\alpha) R'(A),$$
which renders a contradiction given that $R(A)=R'(A)=0.$
\qed

\subsection{Proof of Proposition~\ref{reform_potp}}\label{reform_potp_proof}
We begin by proving the first claim of Proposition~\ref{reform_potp}.

We show that every feasible solution of (SP1) can be mapped into a feasible solution of (SP2) with the same objective function value. To this end, take $Q$ as a feasible solution of (SP1) and let $Q_{\widetilde{\Xi}}$ be the $Q$-conditional probability measure given $\boldsymbol{\xi} \in \widetilde{\Xi}$. Take $\widehat{\mathbb{Q}}_N$ and $Q_{\widetilde{\Xi}}$ as the two probabilities in Corollary \ref{corollary312agullo} with $\alpha \in (0,1)$. There exists $Q_{1-\alpha} \in \mathcal{F}_{1-\alpha}(Q_{\widetilde{\Xi}})$ such that $Q_{1-\alpha}=\alpha Q_{\widetilde{\Xi}}+
(1-\alpha)\widetilde{Q}_{1-\alpha}$, with $\widetilde{Q}_{1-\alpha} \in \mathcal{R}_{\alpha}(\widehat{\mathbb{Q}}_N)$ and
$\mathcal{W}_p(\widehat{\mathbb{Q}}_N,Q_{1-\alpha} )=\mathcal{W}_p(\widehat{\mathbb{Q}}_N, \mathcal{F}_{1-\alpha}(Q_{\widetilde{\Xi}}))$. Furthermore,
it automatically follows from Proposition~\ref{proposition314agullo} that $\mathcal{W}_p^p(\widehat{\mathbb{Q}}_N, \mathcal{F}_{1-\alpha}(Q_{\widetilde{\Xi}}){)} =  \alpha\mathcal{W}_p^p(\mathcal{R}_{1-\alpha}( \widehat{\mathbb{Q}}_N),Q_{\widetilde{\Xi}})$.

Since $Q \in \mathcal{F}_{1-\alpha}(Q_{\widetilde{\Xi}}) $, we deduce that
$\mathcal{W}_p^p(\widehat{\mathbb{Q}}_N, \mathcal{F}_{1-\alpha}(Q_{\widetilde{\Xi}})) \leqslant \mathcal{W}_p^p(\widehat{\mathbb{Q}}_N,Q) \leqslant \widetilde{\rho}\cdot  \alpha$.
Hence, it holds that $\mathcal{W}_p^p(\mathcal{R}_{1-\alpha}( \widehat{\mathbb{Q}}_N),Q_{\widetilde{\Xi}})\leqslant \widetilde{\rho}$.
In other words, $Q_{\widetilde{\Xi}}$ is feasible in (SP2). Besides, since $Q_{\widetilde{\Xi}}$ is the $Q$-conditional probability
measure given $\boldsymbol{\xi} \in \widetilde{\Xi}$, we have that $\mathbb{E}_Q \left[f(\mathbf{x},\boldsymbol{\xi}) \;{\mid}\; \boldsymbol{\xi} \in \widetilde{\Xi} \right]$ = $\frac{1}{\alpha}\mathbb{E}_Q \left[f(\mathbf{x},\boldsymbol{\xi})\mathbb{I}_{\widetilde{\Xi}}(\boldsymbol{\xi})\right]$ $= \mathbb{E}_{Q_{\widetilde{\Xi}}} \left[f(\mathbf{x},\boldsymbol{\xi}) \right]$ a.s.
 \\\\
Next we prove the second claim of the proposition. For this purpose, first we show that, if $\widehat{\mathbb{Q}}_{N}(\widetilde{\Xi}) = 0$, then every feasible solution of (SP2) can also be mapped into a feasible solution of (SP1) with the same objective function value. To this end, take $Q_{\widetilde{\Xi}}$ feasible in (SP2) and consider $\widehat{Q}_{1-\alpha}\in \mathcal{R}_{1-\alpha}(\widehat{\mathbb{Q}}_N)$     such that $\mathcal{W}_p(\widehat{Q}_{1-\alpha},Q_{\widetilde{\Xi}} )$
$=
\mathcal{W}_p(\mathcal{R}_{1-\alpha}(\widehat{\mathbb{Q}}_N),Q_{\widetilde{\Xi}} )$. Fix $\widetilde{Q}_{1-\alpha}=\frac{1}{1-\alpha}(\widehat{\mathbb{Q}}_N-\alpha\widehat{Q}_{1-\alpha})$. By Proposition~\ref{proposition314agullo}, we have
$$Q_{1-\alpha}=\alpha Q_{\widetilde{\Xi}}+(1-\alpha)\widetilde{Q}_{1-\alpha}
=\alpha Q_{\widetilde{\Xi}}+\widehat{\mathbb{Q}}_N-\alpha\widehat{Q}_{1-\alpha}\in \mathcal{F}_{1-\alpha}(Q_{\widetilde{\Xi}})$$
Hence, $Q_{1-\alpha}(\widetilde{\Xi})=\alpha$, because $\widehat{\mathbb{Q}}_{N}(\widetilde{\Xi})$ gives zero measure to $\widetilde{\Xi}$ and so does any of its $(1-\alpha)$-trimmings. Besides, we have that
$$\mathcal{W}_p^p(\widehat{\mathbb{Q}}_N,Q_{1-\alpha} )=\mathcal{W}_p^p(\widehat{\mathbb{Q}}_N, \mathcal{F}_{1-\alpha}(Q_{\widetilde{\Xi}}))=\alpha\mathcal{W}_p^p(\mathcal{R}_{1-\alpha}( \widehat{\mathbb{Q}}_N),Q_{\widetilde{\Xi}}) \leqslant \alpha \widetilde{\rho}.$$
Therefore, $Q_{1-\alpha}$ is feasible in (SP1) and $Q_{\widetilde{\Xi}}$ is the $Q_{1-\alpha}$-conditional probability measure given $\boldsymbol{\xi} \in \widetilde{\Xi}$.

Finally, if $\alpha = 1$, then $\mathcal{R}_{1-\alpha}( \widehat{\mathbb{Q}}_N) = \widehat{\mathbb{Q}}_N$,  $ \mathbb{E}_Q \left[f(\mathbf{x},\boldsymbol{\xi}) \;{\mid}\; \boldsymbol{\xi} \in \widetilde{\Xi} \right] =  \mathbb{E}_Q \left[f(\mathbf{x},\boldsymbol{\xi}) \right]$ and the mapping is direct, namely, $Q=Q_{\widetilde{\Xi}}.$ \qed

\subsection{Proof of Theorem~\ref{Theorem_partial_mass_reformulation_duality}}\label{Theorem_partial_mass_reformulation_duality_proof}
Thanks to Lemma~\ref{lemma_trimming}, the subproblem (SP2) can be written equivalently as follows:
\begin{align*}
 \text{(SP2)}\;\;  &\sup_{Q_{\widetilde{\Xi}};\;{\mathbf{b}\in \Delta(\alpha_N)}} \;\;
     \mathbb{E}_{Q_{\widetilde{\Xi}}} \left[f(\mathbf{x},\boldsymbol{\xi}) \right]\\
&\text{s.t.} \enskip Q_{\widetilde{\Xi}}(\widetilde{\Xi})=1\\
& \phantom{s.t.} \enskip  \mathcal{W}_p\left(\sum_{i=1}^N b_i \delta_{\widehat{\boldsymbol{\xi}}_i},Q_{\widetilde{\Xi}}\right)\leqslant \widetilde{\rho}^{1/p}
\end{align*}
 where $\Delta(\alpha_N)$ stands for the set of constraints $\{ 0 \leqslant b_i \leqslant \frac{1}{N \alpha_N}, \forall i \leqslant N, \sum_{i=1}^N b_i=1\}$.

which, in turn, can be reformulated as
\begin{align}
     & \left\{ \begin{array}{cl}\displaystyle \sup_{Q_{\widetilde{\Xi}};\  \Pi;\ {\mathbf{b}\in \Delta(\alpha_N)}} &{}\displaystyle \int_{\widetilde{\Xi}}f(\mathbf{x},(\mathbf{z},\mathbf{y}))  Q_{\widetilde{\Xi}}(d\mathbf{z},d\mathbf{y}) \\ \text {s.t.}\\
&\displaystyle \int_{\widetilde{\Xi}} Q_{\widetilde{\Xi}}(d\mathbf{z},d\mathbf{y})=1
\\&{}\displaystyle
\left(\int_{\widetilde{\Xi}\times\Xi}
    \left\|(
\mathbf{z},\mathbf{y})- (
\mathbf{z},\mathbf{y})'\right\|^p \Pi
    (d(
\mathbf{z},\mathbf{y}), d(
\mathbf{z},\mathbf{y})')\right)^{1/p} \leqslant
    \widetilde{\rho}^{1/p}  \\[1ex] &{} \left\{ \begin{array}{l} \Pi
    \text{
is } \text{ a } \text{ joint } \text{ distribution } \text{ of }
(
\mathbf{z},\mathbf{y}) \text{ and } (
\mathbf{z},\mathbf{y})'\\ \text{ with }
\text{ marginals } Q_{\widetilde{\Xi}} \text{ and } \sum_{i=1}^N b_i \delta_{\widehat{\boldsymbol{\xi}}_i}\text{, } \text{
respectively } \end{array}\right.  \\\end{array} \right. \\=&
\left\{ \begin{array}{cl} \displaystyle\sup_{Q^i_{\widetilde{\Xi}}; \ {\mathbf{b}\in \Delta(\alpha_N)} } &{}
\displaystyle \sum_{i=1}^N b_i\displaystyle\int_{\widetilde{\Xi}}
f(\mathbf{x},(\mathbf{z},\mathbf{y}))
Q^i_{\widetilde{\Xi}} (d\mathbf{z},d\mathbf{y}) \\ \text
{s.t.}      &\displaystyle \int_{\widetilde{\Xi}}Q^i_{\widetilde{\Xi}} (d
\mathbf{z},d\mathbf{y})=1,\;
\forall i \leqslant N\\
&{}\displaystyle\sum_{i=1}^{N}b_i\int_{\widetilde{\Xi}}
    \left\|(\mathbf{z},\mathbf{y})-(\widehat{\mathbf{z}}_{i},\widehat{\mathbf{y}}_{i})\right\|^p
 Q^i_{\widetilde{\Xi}} (d
\mathbf{z},d\mathbf{y})
     \leqslant \widetilde{\rho}
    \label{rho_constraint_1} \\
\end{array} \right.
\end{align}
where reformulation \eqref{rho_constraint_1} follows from the fact that the marginal
distribution of $(\mathbf{z},\mathbf{y})'$ is the discrete
distribution supported on points $(\widehat{\mathbf{z}}_{i},\widehat{\mathbf{y}}_{i})$, with probability masses $b_i$, $i = 1, \ldots, N$.
Thus, $\Pi$ is completely determined by the conditional
distributions  $Q^i_{\widetilde{\Xi}}$  of $(\mathbf{z},\mathbf{y})$ given
$(\mathbf{z},\mathbf{y})'=(\widehat{\mathbf{z}}_{i},\widehat{\mathbf{y}}_{i})$, $i = 1, \ldots, N$, that is,
$$\Pi
    (d(
\mathbf{z},\mathbf{y}), d(
\mathbf{z},\mathbf{y})')=
\sum_{i=1}^Nb_i
\delta_{(\widehat{\mathbf{z}}_{i},\widehat{\mathbf{y}}_{i})}(d(
\mathbf{z},\mathbf{y})')  Q^i_{\widetilde{\Xi}}(d(
\mathbf{z},\mathbf{y}))$$

Now we split up the supremum into two:
\begin{subequations}
\begin{align}
     \sup_{ {\mathbf{b}\in \Delta(\alpha_N)} }\quad &\sup_{ Q^i_{\widetilde{\Xi}}, \forall i \leqslant N}\sum_{i=1}^N b_i\int_{\widetilde{\Xi}}
f(\mathbf{x},(\mathbf{z},\mathbf{y}))
Q^i_{\widetilde{\Xi}} (d\mathbf{z},d\mathbf{y})\\
 &   \text{s.t} \;\displaystyle \int_{\widetilde{\Xi}}Q^i_{\widetilde{\Xi}} (d
\mathbf{z},d\mathbf{y})=1,\enskip
\forall i \leqslant N\\
&   \phantom{s.t} \; \sum_{i=1}^{N}b_i\int_{\widetilde{\Xi}}
    \left\|(\mathbf{z},\mathbf{y})-(\widehat{\mathbf{z}}_{i},\widehat{\mathbf{y}}_{i})\right\|^p
 Q^i_{\widetilde{\Xi}} (d
\mathbf{z},d\mathbf{y})
     \leqslant \widetilde{\rho} \label{constraint_dual_lamb}
%
\end{align}
\end{subequations}
  If we set $\lambda$ as the dual variable of  constraint~\eqref{constraint_dual_lamb}, then
using standard duality arguments, we can equivalently rewrite the inner supremun as
    \begin{align}
            &  \sup_{  {\mathbf{b}\in \Delta(\alpha_N)}}\inf_{\lambda \geqslant 0} \sup_{ Q^i_{\widetilde{\Xi}}, \forall i \leqslant N}\!\!\!\!\lambda \widetilde{\rho}+\sum_{i=1}^N b_i\!\!\int_{\widetilde{\Xi}}\left(
f(\mathbf{x},(\mathbf{z},\mathbf{y}))-\lambda  \left\|(\mathbf{z},\mathbf{y})-(\widehat{\mathbf{z}}_{i},\widehat{\mathbf{y}}_{i})\right\|^p\right)
Q^i_{\widetilde{\Xi}} (d\mathbf{z},d\mathbf{y})\\
 &  \hspace{4cm}       \text{s.t} \;\;\displaystyle \int_{\widetilde{\Xi}}Q^i_{\widetilde{\Xi}} (d
\mathbf{z},d\mathbf{y})=1,\enskip
\forall i \leqslant N &\\
=&     \sup_{ {\mathbf{b}\in \Delta(\alpha_N)}}\;\inf_{\lambda \geqslant 0}
\lambda \widetilde{\rho}+\sum_{i=1}^N b_i\sup_{ (\mathbf{z},\mathbf{y}) \in  \widetilde{\Xi}}\left(
f(\mathbf{x},(\mathbf{z},\mathbf{y}))-\lambda  \left\|(\mathbf{z},\mathbf{y})-(\widehat{\mathbf{z}}_{i},\widehat{\mathbf{y}}_{i})\right\|^p  \right)\label{swapping}
\\
=&  \inf_{\lambda \geqslant 0} \enskip   \sup_{  {\mathbf{b}\in \Delta(\alpha_N)}}  \lambda \widetilde{\rho}+\sum_{i=1}^N b_i \sup_{ (\mathbf{z},\mathbf{y}) \in  \widetilde{\Xi}}\left(
f(\mathbf{x},(\mathbf{z},\mathbf{y}))-\lambda  \left\|(\mathbf{z},\mathbf{y})-(\widehat{\mathbf{z}}_{i},\widehat{\mathbf{y}}_{i})\right\|^p  \right)
\\
=&  \inf_{\lambda \geqslant 0; \overline{\mu}_i, \forall i\leqslant N; \theta \in \mathbb{R}} \quad   \lambda \widetilde{\rho}+\theta +\dfrac{1}{N\alpha}\sum_{i=1}^N \overline{\mu}_i  \label{supremum_first} \\
& \hspace{2cm}  \text{s.t.} \  \overline{\mu}_i+\theta\geqslant \sup_{ (\mathbf{z},\mathbf{y}) \in  \widetilde{\Xi}}\left(
f(\mathbf{x},(\mathbf{z},\mathbf{y}))-\lambda  \left\|(\mathbf{z},\mathbf{y})-(\widehat{\mathbf{z}}_{i},\widehat{\mathbf{y}}_{i})\right\|^p  \right),\enskip \forall i \leqslant N \label{supremum}\\
& \hspace{2cm}  \phantom{s.t.} \ \overline{\mu}_i \geqslant 0,\enskip \forall i \leqslant N\label{supremum_last}
    \end{align}
where  we have swapped  the supremum and the infimum in \eqref{swapping} by appealing to Sion's min-max theorem \citep{Sion1958},  given that the objective function in \eqref{swapping} is linear in the $b_i, i=1,\ldots, N$, over a compact convex set, and a positively weighted sum of convex functions in $\lambda$. \qed

\begin{remark}[Limiting case $\alpha = 0$]\label{remark:lim_case}
 If $\alpha = 0$, $\mathcal{R}_{1}(\widehat{\mathbb{Q}}_{N}) = \{\sum_{i=1}^{N}{b_i\delta_{\widehat{\boldsymbol{\xi}}_i}}$ such that $b_i \geqslant 0$, $\forall i =1, \ldots, N$, and $\sum_{i=1}^{N}{b_i}=1\}$. Therefore, dual variables $\overline{\mu}_i, \forall i \leqslant N$, do not appear in~\eqref{supremum_first}--\eqref{supremum_last} in this case. Similarly, if $\frac{1}{N\alpha} \geqslant 1$, the constraints $b_i \leqslant \frac{1}{N \alpha}, \forall i \leqslant N$, become redundant and hence we can set $\overline{\mu}_i=0, \forall i \leqslant N$.
\end{remark}

\end{APPENDICES}

\theendnotes

\ACKNOWLEDGMENT{This research has received funding from the European Research Council (ERC) under the European Union's Horizon 2020 research and innovation programme (grant agreement no. 755705). This work was also supported in part by the Spanish Ministry of Economy, Industry and Competitiveness and the European Regional Development Fund (ERDF) through project ENE2017-83775-P.}


\bibliographystyle{informs2014} 
\bibliography{References} 




\ECSwitch


\ECHead{Electronic Companion}

This electronic companion contains some additional material of interest related to the DRO framework we propose to handle conditional stochastic programs.  First, we state some complementary theoretical results. Second, we use tools from nearest neighbors to show that our DRO approach is asymptotically consistent under assumptions slightly different than those made in the main text (some of which are less restrictive). Finally, numerical experiments for the case $\mathbb{Q}(\widetilde{\Xi}) = \alpha >0$ are presented and discussed.

\textbf{Notation}.
Given any norm $\left\| \cdot
\right\|$ in the Euclidean space (of a certain
dimension $d$), the dual norm is defined as
$\left\| \mathbf{u}
\right\|_{*}=\sup_{\left\| \mathbf{v}
\right\| \leqslant 1} \langle
\mathbf{u},\mathbf{v} \rangle$. Let $g$ be a function from $\mathbb{R}^d$ to $\overline{\mathbb{R}}$, we will say that $g$ is a proper function if $g(\mathbf{x}) < +\infty$ for at least one $\mathbf{x}$ and $g(\mathbf{x}) > -\infty$ for all $\mathbf{x} \in \mathbb{R}^d$. In addition, the convex conjugate function of $g$, $g^*$, is given by
$g^*(\mathbf{y}):= \sup_{\mathbf{x} \in \mathbb{R}^d} \langle \mathbf{y},\mathbf{x} \rangle -g(\mathbf{x})$. It is well known that if $g$ is a proper function, then $g^*$ is a proper function as well.
The support function of set $A$, $S_A$, is
defined as
$S_A(\mathbf{b}):= \sup_{\mathbf{a} \in A}
\langle \mathbf{b}, \mathbf{a} \rangle$.
The recession cone of a non-empty set  $A \subseteq \mathbb{R}^d$ is given by $\{ \mathbf{b} \in \mathbb{R}^d\;/\; \mathbf{a}+\lambda \mathbf{b} \in A,\; \forall a \in A, \ \forall \lambda \geqslant 0 \}$.
\section{Complementary theoretical results}
This section contains some theoretical results which are complementary to the theory developed in the manuscript.  First, we introduce a few preliminary concepts and definitions. Second,  we state the topological properties of the ambiguity set $\widehat{\mathcal{U}}_{{N}}(\alpha,\widetilde{\rho})$ in problem (P). Finally, we introduce a tractable reformulation of our DRO approach, which is similar to that in \citet[Theorem 8]{Kuhn2019}.
\subsection{ Auxiliary measure theoretic concepts and Wasserstein metric}\label{auxiliary concepts}
This subsection compiles some definitions and results from the measure theory that underpins our research. It starts with  concepts related to the weak convergence of measures and compactness. Subsequently, some known facts in connection with the topology generated by the Wasserstein metric
$\mathcal{W}_p$ are presented. We denote  the set of all Borel probability measures supported on $\mathcal{X}$ as $\mathcal{P}(\mathcal{X})$. Although some of the following concepts and results are still true in the more general setting of Polish spaces, we restrict ourselves here to $\mathcal{X} \subseteq\mathbb{R}^d$. Similarly, we denote the $p$-Wasserstein space as $\mathcal{P}_p(\mathcal{X})$, that is, the set of all Borel probability measures supported on $\mathcal{X}$ with a finite $p$-th moment.  It is well known that the $p$-Wassertein metric defines a metric in  $\mathcal{P}_p(\mathcal{X})$ \cite[Theorem 7.3]{Villani2003}.

\begin{definition}[Weak convergence of probability measures]
Given a sequence of probability
  measures $\{Q_N\}_N\subseteq\mathcal{P}(\mathcal{X})$, we say that it converges weakly
  to $Q$
 if
 \begin{equation}
     \lim_{N\rightarrow \infty}\int_{\mathcal{X}} \ell(\boldsymbol{\xi})Q_N(d\boldsymbol{\xi})= \int_{\mathcal{X}} \ell(\boldsymbol{\xi})Q(d\boldsymbol{\xi})
 \end{equation}
for all bounded and continuous
function $\ell$ on $\mathcal{X}$.
\end{definition}
 \begin{definition}[Tightness]
 A given   set  $\mathcal{K} \subseteq \mathcal{P}(\mathcal{X})$ is \emph{tight} if for all $\varepsilon>0$, there is
a compact set $X_{\varepsilon} \subset \mathcal{X}$ such that
$\inf_{Q \in \mathcal{K}}Q(X_{\varepsilon})> 1-\varepsilon$. If $\mathcal{K}$ reduces to a singleton, then we refer to the ``tightness
of a probability measure''.
 \end{definition}

  \begin{definition}[Closed sets]
A given   set  $\mathcal{K} \subseteq \mathcal{P}(\mathcal{X})$ is \emph{closed} (under the topology of weak convergence) if for all sequence
 $\{Q_N\}_N \subset \mathcal{K}$ such that $Q_N $ converges weakly to
 $Q$,
 we have $Q \in \mathcal{K}$.
 \end{definition}
 The following theorem, which is known as  Prokhorov's Theorem, connects the notions of weak compactness and tightness.
 \begin{theorem}[Prokhorov's Theorem]
 A  set $\mathcal{K} \subseteq \mathcal{P}(\mathcal{X})$   is tight  if and only if the closure of $\mathcal{K}$  is weakly compact in $\mathcal{P}(\mathcal{X})$.
 \end{theorem}
 \begin{definition}[Weak compactness]
A set $\mathcal{K} \subseteq \mathcal{P}(\mathcal{X})$  is \emph{weakly compact}  if for all sequence of probability measures $\{Q_N\}_N \subset \mathcal{K}$, there exists a subsequence
$\{Q_{N'} \}_{N'}$ that converges weakly to $Q \in \mathcal{K}$.
 \end{definition}

 \begin{definition}[$p$-uniform integrability]
 A set $\mathcal{K}\subseteq \mathcal{P}(\mathcal{X})$  is said to have $p$-uniformly integrable moments if
\begin{equation}
    \lim_{t\rightarrow\infty}\int_{\{\boldsymbol{\xi}/ \| \boldsymbol{\xi}\|>t \}}\|\boldsymbol{\xi}\|^p Q(d\boldsymbol{\xi})=0\;  \text{uniformly w.r.t.}\;  Q\in \mathcal{K}
\end{equation}
 \end{definition}

 Finally, we introduce a proposition that connects some of the aforementioned concepts  with the Wasserstein metric. More concretely, this proposition establishes the topological properties of the Wasserstein space.

\begin{proposition}\label{topology_wasserstein}
Given $p\geqslant1$ and $\mathcal{X}\subseteq \mathbb{R}^d$ a closed set, we have:
$\mathcal{P}_p(\mathcal{X})$ endowed with $\mathcal{W}_p$ is a Polish space.
A closed set $\mathcal{K} \subseteq \mathcal{P}_p(\mathcal{X})$  is weakly compact if and only if it has $p$-uniformly integrable moments (and hence tight).
Specifically, given a sequence of probability measures
 $\{Q_N\}_N\subseteq\mathcal{P}_p(\mathcal{X})$,
 the following statements
 are equivalent:
 \begin{enumerate}
     \item $\mathcal{W}_p(Q_N,Q)\rightarrow 0$.
     \item $Q_N$ converges weakly to $Q$ and
     $\{Q_N\}_N$ has $p$-uniformly integrable moments.
     \item $Q_N$ converges weakly to $Q$ and the following holds
      $$\int_{\mathcal{X}} \|\boldsymbol{\xi}\|^pQ_N(d\boldsymbol{\xi})\stackrel{N\rightarrow\infty}{\longrightarrow} \int_{\mathcal{X}} \|\boldsymbol{\xi}\|^pQ(d\boldsymbol{\xi}).$$

     \item For any $L>0$ and any continuous function $\ell:\mathcal{X}\rightarrow \mathbb{R}$  such that verifies $|\ell(\boldsymbol{\xi})| \leqslant L(1+\|\boldsymbol{\xi}\|^p)$ for all $\boldsymbol{\xi}$, the following holds
     $$\int_{\mathcal{X}} \ell(\boldsymbol{\xi})Q_N(d\boldsymbol{\xi})\stackrel{N\rightarrow\infty}{\longrightarrow} \int_{\mathcal{X}} \ell(\boldsymbol{\xi})Q(d\boldsymbol{\xi}).$$

 \end{enumerate}

\end{proposition}

\begin{remark}
Proposition \ref{topology_wasserstein} compiles results from Prop. 7.1.5 in
\citet{ambrosio2005} and Th. 7.12 in \cite{Villani2003}. It implies  that
 the topology generated by $\mathcal{W}_p$ and the weak topology do coincide on any subset $\mathcal{K}$ which has $p$-uniformly integrable moments. We note that assertion 2 in Proposition \ref{topology_wasserstein} is reduced to weak convergence if  $\mathcal{X}$ is a compact set (see, for example,  \citet[Corollary 2.2.2 ]{Panaretos2020}).
\end{remark}

 \subsection{Topological properties of the ambiguity set}\label{proof_th_compactness}
The following proposition formally establishes   the topological properties of our ambiguity set:
\begin{proposition}\label{th_compactness}
Given $\mathbb{Q}\in \mathcal{P}_p(\mathbb{R}^d)$,  $\alpha>0$, and
$\widetilde{\rho}\geqslant\underline{\epsilon}^p_{N\alpha}$, the ambiguity set of problem \emph{(P)}, $\widehat{\mathcal{U}}_N(\alpha,\widetilde{\rho})$,
is non-empty, tight, weakly compact, and $p$-uniformly integrable.
\end{proposition}

\proof{Proof}
 The set  $\widehat{\mathcal{U}}_N(\alpha,\widetilde{\rho})$ is non-empty, because
$\widetilde{\rho}\geqslant\underline{\epsilon}^p_{N\alpha}$. We can equivalently rewrite  $\widehat{\mathcal{U}}_N(\alpha,\widetilde{\rho})$  as
$$\left\{ { Q_{\widetilde{\Xi}} \in \mathcal{P}(\widetilde{\Xi}) :\mathcal{W}^p_p(R, Q_{\widetilde{\Xi}})\leqslant \widetilde{\rho} \text { for some } R\in \mathcal{R}_{1-\alpha}(\widehat{\mathbb{Q}}_N)}\right\}.$$

If $\alpha>0$, then the trimming set $\mathcal{R}_{1-\alpha}(\widehat{\mathbb{Q}}_N) $ is  tight and weakly compact, see   \citet[Lemmas 2 and 3]{cascos2008}.
Furthermore, $\widehat{\mathcal{U}}_N(\alpha,\widetilde{\rho})$
 is a subset of
$$\mathcal{K}:=\left\{ { Q_{\widetilde{\Xi}} \in \mathcal{P}(\mathbb{R}^d) :\mathcal{W}^p_p(R, Q_{\widetilde{\Xi}})\leqslant \widetilde{\rho} \text { for some } R\in \mathcal{R}_{1-\alpha}(\widehat{\mathbb{Q}}_N)}\right\}$$
which is tight and weakly compact by
 \citet[Proposition 3]{Pichler2018}. The tightness of $\widehat{\mathcal{U}}_N(\alpha,\widetilde{\rho})$ is trivially guaranteed, since any subset of a tight set is also tight. Hence, by Prokhorov's theorem, to demonstrate that  $\widehat{\mathcal{U}}_N(\alpha,\widetilde{\rho})$ is also weakly compact, it suffices to show that it is closed. For this purpose, let  $\{Q^N_{\widetilde{\Xi}}\}_N$ be a sequence of probability measures in $\widehat{\mathcal{U}}_N(\alpha,\widetilde{\rho})$ that converges weakly to $Q$.
We need to show that $Q$ is  in $\widehat{\mathcal{U}}_N(\alpha,\widetilde{\rho}) $ too. In turn,  since $\widehat{\mathcal{U}}_N(\alpha,\widetilde{\rho})$ is a subset of  $\mathcal{K}$, which is closed, this boils down to proving that the weak limit  satisfies the condition
$Q\in \mathcal{P}(\widetilde{\Xi})$, that is,  $Q(\widetilde{\Xi})=1$.
Given that the sequence $\{Q^N_{\widetilde{\Xi}}\}_N$ converges weakly to $Q$ and the support set $\widetilde{\Xi}$ is closed,  Portmanteau's theorem (see \citet[Theorem 2.1]{Billingsley})  tells us that
  $\lim \sup_{N\rightarrow \infty} Q^N_{\widetilde{\Xi}}(\widetilde{\Xi}) =1\leqslant Q(\widetilde{\Xi})$.
This implies that
$Q(\widetilde{\Xi})=1$.

Finally,
the $p$-uniform integrability of our ambiguity set  follows from Proposition \ref{topology_wasserstein}. To apply this proposition, we only need to check whether  any distribution of $\widehat{\mathcal{U}}_N(\alpha,\widetilde{\rho})$ has a finite $p$-th moment.  From  \cite{Alvarez-Esteban2011} (see p. 363 for the  the
case $p = 2$, although the proof works similarly for any $p\geqslant 1$), we know that $\mathcal{R}_{1-\alpha}(\widehat{\mathbb{Q}}_N)\subset \mathcal{P}_p(\mathbb{R}^d)$ if $\mathbb{Q}\in \mathcal{P}_p(\mathbb{R}^d)$. Now,  assume that there is
a distribution $Q_{\widetilde{\Xi}}$ in $\widehat{\mathcal{U}}_N(\alpha,\widetilde{\rho})$ that does not have a finite $p$-th moment. If this were the case, we would have
$\mathcal{W}_p(Q_{\widetilde{\Xi}},R)=\infty$ for some $ R\in \mathcal{R}_{1-\alpha}(\widehat{\mathbb{Q}}_N)$, which is in contradiction with the fact that $\mathcal{W}_p(Q_{\widetilde{\Xi}},R)$ must be less or equal to a finite $\widetilde{\rho}^{1/p}$.
\qed
\endproof
\subsection{Tractable reformulation and maximizer of problem (SP2)}\label{sup_reformulation}
%


Next we provide a more manageable reformulation of  problem (SP2), which can be  used directly to address the decision-making problems considered in our numerical experiments. However, we omit its proof, as it runs in parallel with that of   \citet[Theorem 4.2]{MohajerinEsfahani2018} and  \citet[Theorem 8]{Kuhn2019}. See also \cite{JianzheZhen2021}. Said reformulation relies on the following assumption.
\begin{assumption}\label{assumption_convex}
The region $\widetilde{\Xi}$ is a closed convex set, and
$f(\mathbf{x}, \boldsymbol{\xi}):=\max_{k \leqslant K} g_k(\mathbf{x},
\boldsymbol{\xi})$, with $g_k$, for each $k\leqslant K$, being a proper, concave and upper semicontinuous function with respect to $\boldsymbol{\xi}$ (for any fixed value of $\mathbf{x}\in X$)  and not identically $\infty$ on $\widetilde{\Xi}$.
\end{assumption}
\begin{theorem}\label{Theorem_partial_mass_reformulation}
Let $p,q \geqslant 1$ such that $\frac{1}{p}+\frac{1}{q}=1$.
If Assumption \ref{assumption_convex} holds,  then for any value of $\widetilde{\rho}\geqslant\underline{\epsilon}^p_{N\alpha}$, subproblem \emph{(SP2)} is equivalent to the following finite convex problem:

      \begin{align*}
 {\rm (SP2'')}\;\;     \inf_{\lambda, \overline{\mu}_i,\theta ,\mathbf{v}_{ik},\mathbf{v}'_{ik},\mathbf{w}_{ik}, \mathbf{w}'_{ik} } &
    \lambda\widetilde{\rho} + \theta+  \frac{1}{N\alpha}\sum_{i=1}^{N}\overline{\mu}_i\\
  \text{s.t.} \   \overline{\mu}_i \geqslant    & [-g_k]^*((\mathbf{v}_{ik}, \mathbf{w}_{ik})-(\mathbf{v}'_{ik}, \mathbf{w}'_{ik}))\\ & +S_{\widetilde{\Xi}}((\mathbf{v}'_{ik}, \mathbf{w}'_{ik}))-\left\langle (\mathbf{v}_{ik}, \mathbf{w}_{ik}), (\widehat{\mathbf{z}}_i,\widehat{\mathbf{y}}_i)\right\rangle \\
  & + \varphi(q) \lambda \left\|\dfrac{(\mathbf{v}_{ik}, \mathbf{w}_{ik})}{\lambda}\right\|_{*}^q  -\theta ,\ \forall i \leqslant N, \forall k \leqslant K \label{cons_reformulationsup2}  \\
  & \lambda \geqslant 0\\
& \overline{\mu}_i\geqslant 0, \ \forall i \leqslant N
\end{align*}
where $[-g_k]^*((\mathbf{v}_{ik}, \mathbf{w}_{ik})-(\mathbf{v}'_{ik}, \mathbf{w}'_{ik}))$ is the conjugate function of $-g_k$ evaluated at $(\mathbf{v}_{ik}, \mathbf{w}_{ik})-(\mathbf{v}'_{ik}, \mathbf{w}'_{ik})$ and $S_{\widetilde{\Xi}}$ is the support function of $\widetilde{\Xi}$. Moreover, $\varphi(q)=(q-1)^{q-1}/q^q$ if $q>1$, and $\varphi(1)=1$. If
$\lambda=0$, then $0\left\|\dfrac{(\mathbf{v}_{ik}, \mathbf{w}_{ik})}{0}\right\|_{*}^q :=\lim_{\lambda \downarrow 0} \lambda \left\|\dfrac{(\mathbf{v}_{ik}, \mathbf{w}_{ik})}{\lambda}\right\|_{*}^q  $.
\end{theorem}
In problem ${\rm (SP2'')}$, we have suppressed the dependence of functions $g_k$ on $\mathbf{x}$ for ease of notation.

The following theorem serves  to construct a maximizer (i.e., a worst-case distribution) of problem (SP2). Again, we omit its proof, as it is analogous to the proof of  \citet[Theorem 4.4]{MohajerinEsfahani2018} and   \citet[Theorem~9]{Kuhn2019}.

\begin{theorem}[Worst-case distributions]\label{th_worst_case_partial}
Under the assumptions of Theorem~\ref{Theorem_partial_mass_reformulation}, the worst-case expectation in \emph{(SP2)} is equal to the optimal objective value  of the following finite convex optimization problem
\begin{align*} \left\{ \begin{array}{clll} \mathop {\sup }\limits _{\gamma _{ik}, \mathbf{q}_{ik}} &{} \sum \limits _{i = 1}^{N} \sum \limits _{k = 1}^{K} \gamma _{ik}g _k \big ( \widehat{\boldsymbol{\xi}}_{i} - {\mathbf{q}_{ik} \over \gamma _{ik}}\big ) \\ \text {s.t.}&{} \sum \limits _{i =1}^{N} \sum \limits _{k =1}^{K}\gamma_{ik}
\left\| {\mathbf{q}_{ik} \over \gamma_{ik}}\right\|^p \le \widetilde{\rho} \\
&{} \sum \limits_{i=1}^{N} \sum \limits_{k = 1}^{K} \gamma_{ik}=1\\
&{} \sum \limits_{k = 1}^{K} \gamma_{ik} \leqslant \frac{1}{N\alpha}&{}\forall i \le N\\
&{} \gamma _{ik} \ge 0 &{}\forall i \le N, \quad \forall k\le K \\
&{} \widehat{\boldsymbol{\xi}}_{i} - {\mathbf{q}_{ik} \over \gamma _{ik}} \in \widetilde{\Xi} \;  &{}\forall i \le N, \quad \forall k\le K \end{array} \right. \end{align*}
where
$0g_k(\widehat{\boldsymbol{\xi}}_{i} - {\mathbf{q}_{ik} \over 0})$
is interpreted as the value
which makes the function
$\gamma_{ik}g_k(\widehat{\boldsymbol{\xi}}_{i} - {\mathbf{q}_{ik} \over \gamma_{ik}})$ upper semicontinuous at $(\mathbf{q}_{ik} , \gamma_{ik})=(\mathbf{q}_{ik} ,0)$. Also, the
constraint $\widehat{\boldsymbol{\xi}}_{i} - {\mathbf{q}_{ik}/0} \in \widetilde{\Xi}$
means that $\mathbf{q}_{ik}$ is in the recession cone of $\widetilde{\Xi}$, and
$0\left\|\mathbf{q}_{ik}/0 \right\|^p$ is understood
as $\lim_{\gamma_{ik}\downarrow 0}\gamma_{ik}\left\|\mathbf{q}_{ik}/\gamma_{ik} \right\|^p$.

Moreover,
if we assume that $p>1$ or that $\widetilde{\Xi}$ is bounded (with $p \geqslant 1$), then
if
$(\gamma^* _{ik}, \mathbf{q}^*_{ik})$ maximizes the problem above, we have that the discrete probability distribution $Q_{\widetilde{\Xi}}$ defined as
$$Q_{\widetilde{\Xi}}=\sum \limits _{i = 1}^{N} \sum \limits _{k = 1}^{K} \gamma^* _{ik}\delta_{\boldsymbol{\xi}^*_{ik}}$$
where $\boldsymbol{\xi}^*_{ik}:=\widehat{\boldsymbol{\xi}}_{i} - {\mathbf{q}^*_{ik} \over \gamma^* _{ik}} \in \widetilde{\Xi}$, represents a maximizer of the worst-case expectation  problem.
\end{theorem}

\section{Asymptotic consistency under a nearest neighbors lens}\label{case_alpha_0_supplementary}



In this section, we show that the asymptotic consistency of our DRO framework for the case $\mathbb{Q} \ll \lambda^{d}$ with $\mathbb{Q}(\widetilde{\Xi}) = \alpha = 0$  can also be proved using a nearest-neighbors approach.

If the density of $\mathbb{Q}$ is sufficiently smooth, it is known that $\mathbb{Q}_{\widetilde{\Xi}}$ can be inferred from information on $\mathbb{Q}$ within a neighborhood of $\mathbf{z}= \mathbf{z^*}$. This essentially means that the portion of mass from the empirical distribution $\widehat{\mathbb{Q}}_{N}$ that is the closest to $\widetilde{\Xi}$ is statistically representative of the conditional distribution $\mathbb{Q}_{\widetilde{\Xi}}$.
Inspired by popular data-driven local predictive methods such as  $K$ nearest neighbours and kernel regression, we can solve problem (P) for a series of pairs $(\alpha_N, \widetilde{\rho}_N)$, both of which tend to zero appropriately as $N$ increases.
%
Indeed, we will demonstrate that, in doing so, problem $\left({\rm P}_{(\alpha_N,\widetilde{\rho}_N)}\right)$ naturally produces distributionally robustified versions of those popular methods when applied to solve problem~\eqref{conditional_expectation_problem}.
Next, we formalize these ideas.

\begin{remark}
Throughout this section, we will assume that $\textrm{dist}
    (\widehat{\boldsymbol{\xi}}_i,
    \widetilde{\Xi}) = \|\mathbf{\widehat{z}}_i-\mathbf{z}^*\|$. This assumption is standard in the technical literature.  The geometry of the joint support set $\Xi$ is expected to have a negligible impact on the asymptotic performance of problem $\left({\rm P}_{(\alpha_N,\widetilde{\rho}_N)}\right)$ (i.e., for large samples), because, under a smoothness condition on $\mathbb{Q}$ and $K/N \rightarrow 0$, it holds that $\textrm{dist}
    (\widehat{\boldsymbol{\xi}}_{K:N},
    \widetilde{\Xi}) \rightarrow \|\mathbf{\widehat{z}}_{K:N}-\mathbf{z}^*\| \rightarrow 0$ almost surely (see   \citet[Lemmas 2.2 and 2.3]{Biau2015}), where $\mathbf{\widehat{z}}_{K:N}$ is the $\mathbf{z}$-component of the $K$-\emph{th} nearest neighbor to $\mathbf{z}^*$ after reordering the data sample $\{\widehat{\boldsymbol{\xi}}_i:=(\mathbf{\widehat{z}}_i, \mathbf{\widehat{y}}_i)\}_{i=1}^{N}$ in terms of $\|\mathbf{\widehat{z}}_{i}-\mathbf{z}^*\|$ only.
\end{remark}

Here, we show that the solutions of the distributionally robust optimization problem $\left({\rm P}_{(\alpha_N,\widetilde{\rho}_N)}\right)$ converge to the solution of the targeted conditional stochastic program \eqref{conditional_expectation_problem} as $N$ increases, for a careful choice of parameters $\alpha_N$ and $\widetilde{\rho}_N$. This result is underpinned by the fact that, under that selection of parameters $\alpha_N$ and $\widetilde{\rho}_N$, any distribution in
$\widehat{\mathcal{U}}_N(\alpha_N,\widetilde{\rho}_N)$
converges to the true conditional distribution $\mathbb{Q}_{\widetilde{\Xi}}$.
\begin{assumption}[Lipschitz-regularity]\label{assumption_lipschitz}
 We assume that there exists an integrable function
$\ell:\mathbb{R}^{d_{\mathbf{y}}} \rightarrow \mathbb{R}_{+}$ such that for all $\mathbf{y} \in \mathbb{R}^{d_{\mathbf{y}}} $
\begin{equation}\label{condition_function_density_lip}
\left|\phi_{\mathbf{y}/\mathbf{z}=\mathbf{z}'}(\mathbf{y})-\phi_{\mathbf{y}/\mathbf{z}=\mathbf{z}^*}(\mathbf{y})  \right| \leqslant \ell(\mathbf{y}) \|\mathbf{z}'-\mathbf{z}^* \|,\enskip \forall \mathbf{z}' \textrm{ such that } \|\mathbf{z}'-\mathbf{z}^* \| \leqslant r_0
\end{equation}
where $\phi_{\mathbf{y}/\mathbf{z}=\mathbf{z}'}(\cdot)$ stands for the density function of $\mathbf{y}$ conditional on $\mathbf{z}=\mathbf{z}'$.

\end{assumption}
\begin{lemma}[Convergence of transported trimmed distributions]\label{lemma_convergence_trimmed_distributions}
Suppose that Assumptions \ref{assumption_knn} and \ref{assumption_lipschitz} hold. Take $(\alpha_N, \widetilde{\rho}_N)$ such that  $\alpha_N
\rightarrow 0$, $\frac{N\alpha_N}{\log(N)}\rightarrow \infty$,
and $\widetilde{\rho}_N \rightarrow 0$ a.s., with $\widetilde{\rho}_N \geqslant \underline{\epsilon}^p_{N\alpha_N}$, where $\underline{\epsilon}_{N\alpha_N}$ is the minimum transportation budget as in Definition~\ref{MTB}.
%
Then, we have that
$$\mathcal{W}_p( Q^N_{\widetilde{\Xi}}, \mathbb{Q}_{\widetilde{\Xi}}) \rightarrow 0 \enskip a.s.$$
where $Q^N_{\widetilde{\Xi}}:=
\sum_{i=1}^N b^N_i \delta_{(\mathbf{z}^*,
\mathbf{\widehat{y}}_i)} \in \widehat{\mathcal{U}}_N(\alpha_N, \widetilde{\rho}_N)$
is the distribution that results from transporting
the distribution $\sum_{i=1}^N b^N_i \delta_{(\mathbf{\widehat{z}}_i,
\mathbf{\widehat{y}}_i)}$ in the trimming set $\mathcal{R}_{1-\alpha_N}(\widehat{\mathbb{Q}}_N)$ onto $\widetilde{\Xi}$.
\end{lemma}
\proof{Proof}
%
%
Since $\mathbf{y}$ is bounded, we only need to prove that $Q^N_{\widetilde{\Xi}}$ converges weakly to $\mathbb{Q}_{\widetilde{\Xi}}$. For this purpose, take a continuous and bounded function $h$ and let $m(\mathbf{z}^*)=\mathbb{E}[h(\mathbf{y})\mid\mathbf{z}=\mathbf{z}^*]$. We have

\begin{equation*}
 \left| \sum_{i=1}^N b^N_i h(\widehat{\mathbf{y}}_i)
    -m(\mathbf{z}^*)\right| \leqslant \left| \sum_{i=1}^N b^N_i h(\widehat{\mathbf{y}}_i)
    -\sum_{i=1}^N b^N_i m(\widehat{\mathbf{z}}_i)
    \right|
    + \left|\sum_{i=1}^N b^N_i m(\widehat{\mathbf{z}}_i)
    -m(\mathbf{z}^*)
    \right|
\end{equation*}

We deal with each of the terms in the inequality above one by one. First, we use \citet[Lemma 6]{devroye1982necessary} to get
%
%
\begin{equation*}
    \mathbb{P}\left( \left| \sum_{i=1}^N b^N_i h(\widehat{\mathbf{y}}_i)
    -\sum_{i=1}^N b^N_i m(\widehat{\mathbf{z}}_i)
    \right| > \varepsilon\  \mid \ \widehat{\mathbf{z}}_1, \ldots, \widehat{\mathbf{z}}_N  \right) \leqslant 2\exp{
    \left(\dfrac{-(N\alpha_N) \varepsilon^2}{4\left\|h\right\|_{\infty}(2
    \left\|h\right\|_{\infty}+\varepsilon)}\right)}
\end{equation*}
Given that
\begin{equation*}
    \mathbb{P}\left( \left| \sum_{i=1}^N b^N_i \left(h(\widehat{\mathbf{y}}_i)
    - m(\widehat{\mathbf{z}}_i)\right)
    \right| > \varepsilon \right) = \mathbb{E}\left[\mathbb{P}\left( \left| \sum_{i=1}^N b^N_i \left(h(\widehat{\mathbf{y}}_i)
    - m(\mathbf{\widehat{z}}_i)\right)
    \right| > \varepsilon\  \mid \ \widehat{\mathbf{z}}_1, \ldots, \widehat{\mathbf{z}}_N  \right)\right]
\end{equation*}
we have
\begin{equation}\label{Conv_transp_trimmed_dist_1term}
    \mathbb{P}\left( \left| \sum_{i=1}^N b^N_i h(\widehat{\mathbf{y}}_i)
    -\sum_{i=1}^N b^N_i m(\widehat{\mathbf{z}}_i)
    \right| > \varepsilon \right) \leqslant 2\exp{
    \left(\dfrac{-(N\alpha_N) \varepsilon^2}{4\left\|h\right\|_{\infty}(2
    \left\|h\right\|_{\infty}+\varepsilon)}\right)}
\end{equation}

Now let $\ell$
 be an integrable function satisfying condition \eqref{condition_function_density_lip}.
 Hence, for any $\mathbf{z}'$ such that $\|\mathbf{z}'-\mathbf{z}^* \| \leqslant r_0$
 \begin{equation}\label{Lip_cond_dens}
     | m(\mathbf{z}') -m(\mathbf{z}^*) | \leqslant \left\|h\right\|_{\infty} \left\|\ell\right\|_{1} \|\mathbf{z}'-\mathbf{z}^* \|=: L\  \|\mathbf{z}'-\mathbf{z}^* \|
 \end{equation}

In addition,
\begin{equation*}
    \left|\sum_{i=1}^N b^N_i m(\widehat{\mathbf{z}}_i)-m(\mathbf{z}^*)  \right|=
     \left|\sum_{i=1}^N b^N_i( m(\widehat{\mathbf{z}}_i)-m(\mathbf{z}^*)) \right|\leqslant
     \sum_{i=1}^N b^N_i \left|m(\widehat{\mathbf{z}}_i)-m(\mathbf{z}^*)\right|
\end{equation*}

Let $J$ be the number of samples such that their distance from the set $\widetilde{\Xi}$ is smaller than or equal to $r_0$. We can write
\begin{align*}
    \left|\sum_{i=1}^N b^N_i m(\widehat{\mathbf{z}}_i)-m(\mathbf{z}^*)  \right| &\leqslant
    \sum_{i=1}^J b^N_{i:N} \left|m(\widehat{\mathbf{z}}_{i:N})-m(\mathbf{z}^*)\right| + \sum_{i=J+1}^N b^N_{i:N} \left|m(\widehat{\mathbf{z}}_{i:N})-m(\mathbf{z}^*)\right|\\
     &\leqslant L\sum_{i=1}^{J}b_{i:N}^N
     \left\|\widehat{\mathbf{z}}_{i:N}-
     \mathbf{z}^*\right\| + 2\left\|h\right\|_{\infty}\sum_{i=J+1}^{N}{b_{i:N}^N}\\
     &\leqslant L\  \mathcal{W}_{1}\left(Q^N_{\widetilde{\Xi}}, \sum_{i=1}^N{b_{i}^N\delta_{(\widehat{\mathbf{z}}_i, \widehat{\mathbf{y}}_i)}}\right) +  2\left\|h\right\|_{\infty}\sum_{i=J+1}^{N}{b_{i:N}^N}\\
     &\leqslant L\  \mathcal{W}_{p}\left(Q^N_{\widetilde{\Xi}}, \sum_{i=1}^N{b_{i}^N\delta_{(\widehat{\mathbf{z}}_i, \widehat{\mathbf{y}}_i)}}\right) +  2\left\|h\right\|_{\infty}\sum_{i=J+1}^{N}{b_{i:N}^N}\\
     &\leqslant L \left(\widetilde{\rho}_N\right)^{\frac{1}{p}} +  2\left\|h\right\|_{\infty}\sum_{i=J+1}^{N}{b_{i:N}^N}
\end{align*}

Next we upper bound the second term in right-hand side of the last inequality.
\begin{align*}
& \sum_{i=J+1}^{N}{b_{i:N}^N}\leqslant \sup\left\{\sum_{i=J+1}^{N}\!\!\!{b_{i:N}^N},  0\leqslant b^N_{i:N} \leqslant \frac{1}{N\alpha_N},\forall i;\ \sum_{i=1}^N b^N_{i:N} = 1; \  \sum_{i=1}^N b^N_{i:N} \left\|\widehat{\mathbf{z}}_{i:N}-\mathbf{z}^*\right\|^p \leqslant \widetilde{\rho}_N \right\}\\
 &=\inf\Bigg\{\frac{1}{N\alpha_N}\sum_{i=1}^{N}\mu_{i:N}+\theta + \lambda \widetilde{\rho}_N, \enskip \mu_{i:N}+\theta+\lambda\left\|\widehat{\mathbf{z}}_{i:N}-\mathbf{z}^*\right\|^p-\gamma_{i:N} =0,\  \forall i=1, \ldots, J;\\
&\phantom{=\inf\Bigg\{} \mu_{i:N}+\theta+\lambda\left\|\widehat{\mathbf{z}}_{i:N}-\mathbf{z}^*\right\|^p-\gamma_{i:N} =1, \ \forall i=J+1, \ldots, N; \ \lambda \geqslant 0; \ \gamma_{i:N}, \mu_{i:N} \geqslant 0, \forall i\Bigg\}
\end{align*}
It suffices to take a feasible solution. In particular, we consider $\mu_{i:N} =0,\  \forall i$, $\theta = 0$, and $\lambda =1/r_{0}^p$, which renders
\begin{equation*}
 \sum_{i=J+1}^N b^N_{i:N} \leqslant \frac{ \widetilde{\rho}_N}{r_0^p}
\end{equation*}
Hence,
\begin{equation*}
 \left|\sum_{i=1}^N b^N_i m(\widehat{\mathbf{z}}_i)-m(\mathbf{z}^*)  \right| \leqslant L \left(\widetilde{\rho}_N\right)^{\frac{1}{p}} +  \frac{2\left\|h\right\|_{\infty}}{r_0^p}\ \widetilde{\rho}_N
\end{equation*}
Consequently, we essentially need that $\lim_{N\rightarrow \infty}{\widetilde{\rho}_N = 0}$ with probability one. To show this, as $\widetilde{\rho}_N \geqslant \underline{\epsilon}^p_{N\alpha_N}$, we decompose $\widetilde{\rho}_N$ into $\underline{\epsilon}^p_{N\alpha_N}$ plus $\Delta\widetilde{\rho}_N$ and use $(\underline{\epsilon}^p_{N\alpha_N} +\Delta\widetilde{\rho}_N)^{1/p} \leqslant \underline{\epsilon}_{N\alpha_N} + \left(\Delta\widetilde{\rho}_N\right)^{1/p}$ to recast the  expression above as
\begin{equation*}
 \left|\sum_{i=1}^N b^N_i m(\widehat{\mathbf{z}}_i)-m(\mathbf{z}^*)  \right| \leqslant L\  \underline{\epsilon}_{N\alpha_N} + \frac{2\left\|h\right\|_{\infty}}{r_0^p}\ \underline{\epsilon}^p_{N\alpha_N} + L \ \left(\Delta\widetilde{\rho}_N\right)^{\frac{1}{p}} +  \frac{2\left\|h\right\|_{\infty}}{r_0^p}\ \Delta\widetilde{\rho}_N
\end{equation*}

Importantly, the budget $\Delta \widetilde{\rho}_N$ is under the decision-maker's control, who simply needs to guarantee that $\Delta \widetilde{\rho}_N \rightarrow 0$ so that the last two terms on the right-hand side of the previous inequality vanishes. Group these two terms into $a_N(\Delta\widetilde{\rho}_N)$, set $K:=\lceil N\alpha_N \rceil$ and note that $\underline{\epsilon}_{N\alpha_N} \leqslant \left\|\widehat{\mathbf{z}}_{K:N}-\mathbf{z}^*\right\|$.

Thus, for any arbitrary $\varepsilon >0$,
\begin{align*}
     \mathbb{P}\left(\left|\sum_{i=1}^N b^N_i m(\widehat{\mathbf{z}}_i)-m(\mathbf{z}^*)
     \right|-a_N(\Delta\widetilde{\rho}_N) > \varepsilon  \right)& \leqslant
     \mathbb{P}\left(
     L\left\|\widehat{\mathbf{z}}_{K:N}-\mathbf{z}^*\right\|  > \frac{\varepsilon}{2} \right)\\
     &+ \mathbb{P}\left(\frac{2\left\|h\right\|_{\infty}}{r_0^p} \left\|\widehat{\mathbf{z}}_{K:N}-\mathbf{z}^*\right\|^p > \frac{\varepsilon}{2}\right)
\end{align*}
In turn,
\begin{align*}
  \mathbb{P}\left(
     L\left\|\widehat{\mathbf{z}}_{K:N}-\mathbf{z}^*\right\|   > \frac{\varepsilon}{2} \right) &= \mathbb{P}\left(
     \left\|\widehat{\mathbf{z}}_{K:N}-\mathbf{z}^*\right\|  > \frac{\varepsilon}{2L}\right)\\
    \mathbb{P}\left(\frac{2\left\|h\right\|_{\infty}}{r_0^p}\ \left\|\widehat{\mathbf{z}}_{K:N}-\mathbf{z}^*\right\|^p > \frac{\varepsilon}{2}\right) &= \mathbb{P}\left(\left\|\widehat{\mathbf{z}}_{K:N}-\mathbf{z}^*\right\| > r_0\left(\frac{\varepsilon}{4\left\|h\right\|_{\infty}}\right)^{\frac{1}{p}}\right)
\end{align*}
Furthermore, due to the first point in Assumption~\ref{assumption_knn}, it holds that
$$\mathbb{P}\left(\left\|\widehat{\mathbf{z}}_{K:N}-\mathbf{z}^*\right\|> \eta\right) \leqslant \text{exp}\left(-\frac{\widetilde{C}}{8}N\eta^{d_\mathbf{z}}\right)$$
for any $0<\eta \leqslant r_0$ and provided that $\frac{K}{N} \leqslant \frac{\widetilde{C}}{2}\eta^{d_\mathbf{z}}$ (see   \citet[formula (34)]{loubes2017}, which is an application of the lower-tail of Chernoff's bound).

Therefore, in that case,
\begin{align*}
     \mathbb{P}\left(\left|\sum_{i=1}^N b^N_i m(\widehat{\mathbf{z}}_i)-m(\mathbf{z}^*)
     \right|-a_N(\Delta\widetilde{\rho}_N) > \varepsilon  \right)& \leqslant \text{exp}\left(-\frac{\widetilde{C}}{8}N \left(\frac{\varepsilon}{2L}\right)^{d_\mathbf{z}}\right)\\
     &+ \text{exp}\left(-\frac{\widetilde{C}}{8}N r_0^{d_\mathbf{z}}\left(\frac{\varepsilon}{4\left\|h\right\|_{\infty}}\right)^{\frac{d_\mathbf{z}}{p}}\right)
\end{align*}
whenever
$$\frac{K}{N} \leqslant \min\left\{\frac{\widetilde{C}}{2}\left(\frac{\varepsilon}{2L}\right)^{d_\mathbf{z}},\  \frac{\widetilde{C}\, r_0^{d_{\mathbf{z}}}}{2}\left(\frac{\varepsilon}{4\left\|h\right\|_{\infty}}\right)^{\frac{d_\mathbf{z}}{p}}\right\} $$
which we guarantee, for $N$ large enough, by enforcing $\alpha_N \rightarrow 0$.

This way, for any arbitrarily small $\varepsilon>0$,   we finally have
     \begin{align}
            \mathbb{P}\left(\left| \sum_{i=1}^N b^N_i h(\widehat{\mathbf{y}}_i)- m(\mathbf{z}^*)
     \right| -a_N(\Delta\widetilde{\rho}_N)> \varepsilon  \right) &\leqslant
      2\exp{\left(\dfrac{-(N\alpha_N) \left(\varepsilon/3\right)^2}{4\left\|h\right\|_{\infty}(2\left\|h\right\|_{\infty}+\varepsilon/3)}\right)}\notag\\
      &+
      \text{exp}\left(-\frac{\widetilde{C}}{8}N \left(\frac{\varepsilon}{3L}\right)^{d_\mathbf{z}}\right)\notag\\
      &+ \text{exp}\left(-\frac{\widetilde{C}}{8}N r_0^{d_\mathbf{z}}\left(\frac{\varepsilon}{6\left\|h\right\|_{\infty}}\right)^{\frac{d_\mathbf{z}}{p}}\right)\label{convergence_summable}
     \end{align}

The last two terms on the right-hand side of~\eqref{convergence_summable} are summable over $N$, while the first one is summable if $\frac{N\alpha_N}{\log(N)}\rightarrow \infty$. Consequently, the Borel-Cantelli Lemma allows us to conclude that
\begin{align*}
 \mathbb{P}\left(\lim_{N\rightarrow \infty}\left|\sum_{i=1}^N b^N_i h(\widehat{\mathbf{y}}_i)- m(\mathbf{z}^*)
     \right| -a_N(\Delta\widetilde{\rho}_N)=0\right)& = \mathbb{P}\left(\lim_{N\rightarrow \infty}\left|\sum_{i=1}^N b^N_i h(\widehat{\mathbf{y}}_i)- m(\mathbf{z}^*)
     \right| =0\right)\\
     &= 1
\end{align*}

given that $a_N \rightarrow 0$ when $\Delta\widetilde{\rho}_N \rightarrow 0$. Thus, $Q^N_{\widetilde{\Xi}}$ converges weakly to $\mathbb{Q}_{\widetilde{\Xi}}$ almost surely.\qed
\endproof

The following corollary extends the convergence to any distribution in the proposed ambiguity set (apart from the transported trimmings of the empirical distribution).

\begin{corollary}[Convergence of conditional distributions]\label{Convergence_conditional}
Suppose that the conditions in Lemma~\ref{lemma_convergence_trimmed_distributions} hold. Then, it follows that
$$\mathcal{W}_p( Q^N_{\widetilde{\Xi}}, \mathbb{Q}_{\widetilde{\Xi}}) \rightarrow 0 \enskip a.s.$$
where $Q^N_{\widetilde{\Xi}}$ is any distribution from the ambiguity set $ \widehat{\mathcal{U}}_N(\alpha_N, \widetilde{\rho}_N)$.
\end{corollary}
\proof{Proof} This corollary is an immediate result of the previous lemma. With some abuse of notation, let $\sum_{i=1}^N b^N_i \delta_{(\mathbf{\widehat{z}}_i,
\mathbf{\widehat{y}}_i)}$ be the distribution in the trimming set $\mathcal{R}_{1-\alpha_N}(\widehat{\mathbb{Q}}_N)$ such that $\mathcal{W}_p\left(\mathcal{R}_{1-\alpha_N}( \widehat{\mathbb{Q}}_N),Q^N_{\widetilde{\Xi}}\right) = \mathcal{W}_p\left(\sum_{i=1}^N b^N_i \delta_{(\mathbf{\widehat{z}}_i,
\mathbf{\widehat{y}}_i)},Q^N_{\widetilde{\Xi}}\right)$.

By the triangle inequality, we have
\begin{equation}\label{ineq_corollary_conv}
 \mathcal{W}_p( Q^N_{\widetilde{\Xi}}, \mathbb{Q}_{\widetilde{\Xi}}) \leqslant
 \mathcal{W}_p\left( Q^N_{\widetilde{\Xi}}, \sum_{i=1}^N b^N_i \delta_{(\mathbf{\widehat{z}}_i,
\mathbf{\widehat{y}}_i)}\right)+ \mathcal{W}_p\left(\sum_{i=1}^N b^N_i \delta_{(\mathbf{\widehat{z}}_i,
\mathbf{\widehat{y}}_i)}, \mathbb{Q}_{\widetilde{\Xi}}\right)
\end{equation}
where
$\mathcal{W}^p_p\left( Q^N_{\widetilde{\Xi}},\sum_{i=1}^N b^N_i \delta_{(\mathbf{\widehat{z}}_i,
\mathbf{\widehat{y}}_i)}\right) \leqslant \widetilde{\rho}_N$, because $Q^N_{\widetilde{\Xi}} \in \widehat{\mathcal{U}}_N(\alpha_N, \widetilde{\rho}_N)$.
We  again use the triangle inequality to upper bound the second term on the right-hand side of \eqref{ineq_corollary_conv}.
$$\mathcal{W}_p\left(\sum_{i=1}^N b^N_i \delta_{(\mathbf{\widehat{z}}_i,
\mathbf{\widehat{y}}_i)}, \mathbb{Q}_{\widetilde{\Xi}}\right) \leqslant \mathcal{W}_p\left(\sum_{i=1}^N b^N_i \delta_{(\mathbf{\widehat{z}}_i,
\mathbf{\widehat{y}}_i)}, \sum_{i=1}^N b^N_i \delta_{(\mathbf{z}^*,
\mathbf{\widehat{y}}_i)}\right) + \mathcal{W}_p\left(\sum_{i=1}^N b^N_i \delta_{(\mathbf{z}^*,
\mathbf{\widehat{y}}_i)}, \mathbb{Q}_{\widetilde{\Xi}}\right)$$
where $\sum_{i=1}^N b^N_i \delta_{(\mathbf{z}^*,
\mathbf{\widehat{y}}_i)}$ is the distribution with support on $\widetilde{\Xi}$ that is the closest (in $p$-Wasserstein distance) to $\sum_{i=1}^N b^N_i \delta_{(\mathbf{\widehat{z}}_i,
\mathbf{\widehat{y}}_i)}$. Therefore,
$$\mathcal{W}_p^p\left(\sum_{i=1}^N b^N_i \delta_{(\mathbf{\widehat{z}}_i,
\mathbf{\widehat{y}}_i)}, \sum_{i=1}^N b^N_i \delta_{(\mathbf{z}^*,
\mathbf{\widehat{y}}_i)}\right)\leqslant \mathcal{W}^p_p\left(\sum_{i=1}^N b^N_i \delta_{(\mathbf{\widehat{z}}_i,
\mathbf{\widehat{y}}_i)}, Q^N_{\widetilde{\Xi}}\right) \leqslant \widetilde{\rho}_N$$
That is, $\sum_{i=1}^N b^N_i \delta_{(\mathbf{z}^*,
\mathbf{\widehat{y}}_i)}$ is in $\widehat{\mathcal{U}}_N(\alpha_N, \widetilde{\rho}_N)$ and is precisely one of the transported trimmed distributions to which Lemma~\ref{lemma_convergence_trimmed_distributions} refers.

Hence,
$$\mathcal{W}_p( Q^N_{\widetilde{\Xi}}, \mathbb{Q}_{\widetilde{\Xi}}) \leqslant 2(\widetilde{\rho}_{N})^{\frac{1}{p}}+\mathcal{W}_p\left(\sum_{i=1}^N b^N_i \delta_{(\mathbf{z}^*,
\mathbf{\widehat{y}}_i)}, \mathbb{Q}_{\widetilde{\Xi}}\right)$$

Since both $\widetilde{\rho}_N \rightarrow 0$ and $\mathcal{W}_p\left(\sum_{i=1}^N b^N_i \delta_{(\mathbf{z}^*,
\mathbf{\widehat{y}}_i)}, \mathbb{Q}_{\widetilde{\Xi}}\right) \rightarrow 0$ a.s. by Lemma~\ref{lemma_convergence_trimmed_distributions}, the claim of the corollary follows. \qed
\endproof
Finally, the following theorem formally states the asymptoptic consistency guarantee of our model.

\begin{theorem}[Asymptotic consistency]\label{theorem_consistency_trimmed}
Suppose that the assumptions in Corollary~\ref{Convergence_conditional}
 hold.  Then, we have
 \begin{enumerate}
 \item[(i)] If for any fixed $\boldsymbol{\xi}\in \widetilde{\Xi}$, $f(\cdot,\boldsymbol{\xi})$ is   continuous on $X$, and for any fixed value $\mathbf{x}\in X$,  $f(\mathbf{x},\boldsymbol{\xi})$ is continuous in $\boldsymbol{\xi}$ and there is  $L\geqslant 0$ such that
     $|f(\mathbf{x},\boldsymbol{\xi})|\leqslant L(1+ \left\|\boldsymbol{\xi} \right\|^p) $  for all $\mathbf{x}\in X$ and $\boldsymbol{\xi} \in \widetilde{\Xi}$, then we have that $\widehat{J}_N \rightarrow J^*$ almost surely when $N$ grows to infinity.
     \item[(ii)] Let $X_N, X^*$ be the  set of optimal solutions of problems $\left({\text P}_{(\alpha_N,\widetilde{\rho}_N)}\right)$ and \eqref{conditional_expectation_problem_Case1}, respectively. If the assumptions in (i) are satisfied, the feasible set $X$ is closed and $X_N,X^*$ are non-empty,  then we have that any accumulation point of the sequence $\{ \widehat{\mathbf{x}}_N \}_N$ is almost surely an optimal solution of problem~\eqref{conditional_expectation_problem_Case1}.

 \end{enumerate}
\end{theorem}\
\proof{Proof}
 Set $v_N(\mathbf{x})= \sup_{Q_{\widetilde{\Xi}} \in \widehat{\mathcal{U}}_N(\alpha_N, \widetilde{\rho}_N) }\mathbb{E}_{Q_{\widetilde{\Xi}}}[f(\mathbf{x},\boldsymbol{\xi})]$ and $v(\mathbf{x})= \mathbb{E}_{\mathbb{Q}_{\widetilde{\Xi}}}[f(\mathbf{x},\boldsymbol{\xi})]$. Let  $\mathcal{F}$ be  the class of random functions defined as follows
  \begin{equation}
\begin{split}
\mathcal{F}:= & \left\{f(\boldsymbol{\xi}):=f(\mathbf{x},\boldsymbol{\xi})\; \text{continuous}\;\text{such that}\; \mathbf{x}\in X \right.\\
    &   \text{ and } \exists L\geqslant 0 \
\text{with}\ |f(\mathbf{x},\boldsymbol{\xi})|\leqslant L(1+ \left\|\boldsymbol{\xi} \right\|^p),\ \forall \mathbf{x} \in X, \forall  \boldsymbol{\xi} \in \widetilde{\Xi}\}
\end{split}
\end{equation}
%
%
and let
$\mathcal{D}$ be the pseudometric between two probability measures $P$ and $Q$ given by
$$\mathcal{D}(P,Q):=\sup_{f\in \mathcal{F}}\left|\mathbb{E}_{P}[f]-\mathbb{E}_{Q}[f]\right|$$
For two sets of probability measures $\mathcal{U}_1 \text{ and } \mathcal{U}_2$, define the \emph{excess} of $\mathcal{U}_1$ over $\mathcal{U}_2$ as
$$\mathcal{D}(\mathcal{U}_1,\mathcal{U}_2):=\sup_{P\in \mathcal{U}_1}\inf_{Q \in \mathcal{U}_2}\mathcal{D}(P,Q)$$
First, we show that $v_N(\mathbf{x})<\infty$ for all $\mathbf{x}\in X$. Fix $\mathbf{x}\in X$ and define
$$ \mathcal{V}:=\{\mathbb{E}_{\mathbb{Q}_{\widetilde{\Xi}}}[f(\mathbf{x}, \boldsymbol{\xi})] \} \;\text{and}\;  \mathcal{V}_N:=\{\mathbb{E}_{Q_{\widetilde{\Xi}}}[f(\mathbf{x}, \boldsymbol{\xi})]\; : \;  Q_{\widetilde{\Xi}} \in \widehat{\mathcal{U}}_N(\alpha_N, \widetilde{\rho}_N)\}. $$
The function $f$  satisfies the following uniform-integrability-type condition for all $\mathbf{x}$,
$$\lim_{t\rightarrow \infty}\sup_{Q_{\widetilde{\Xi}} \in \widehat{\mathcal{U}}_N(\alpha_N,\widetilde{\rho_N})}\int_{\{\widetilde{\Xi}: |f(\mathbf{x},\boldsymbol{\xi})|\geqslant t\}}|f(\mathbf{x},\boldsymbol{\xi})|Q_{\widetilde{\Xi}}(d\boldsymbol{\xi})=0$$ due to the limitation on the maximum growth of $f$ established in point (i) and the $p$-uniform integrability of $\widehat{\mathcal{U}}_N(\alpha_N,\widetilde{\rho_N})$. Furthermore, the set $\widehat{\mathcal{U}}_N(\alpha_N,\widetilde{\rho_N})$  is also tight. Consequently,
using  \citet[Proposition 1]{Sun2016}, we deduce that the set $\mathcal{V}_N$ is compact (and hence bounded). Thus, $v_N(\mathbf{x})<\infty$.

  Let $a_N:=\inf_{v \in \mathcal{V}_N}v,\ b_N:=\sup_{v \in \mathcal{V}_N} v$ and $c:=\inf_{v \in \mathcal{V}}v=\sup_{v \in \mathcal{V}}v$.
 Now, denote the Hausdorff distance between the respective convex hulls of the sets $\mathcal{V}$ and $\mathcal{V}_N$ as $\mathbb{H}(\text{conv} \mathcal{V},\text{conv} \mathcal{V}_N)$. We have
 $$\mathbb{H}(\text{conv} \mathcal{V},\text{conv} \mathcal{V}_N)=\mathbb{H}( \mathcal{V},\text{conv} \mathcal{V}_N)=\max \{|b_N-c|,|c-a_N| \}$$
 where
 $$b_N-c=\max_{Q_{\widetilde{\Xi}}\in \widehat{\mathcal{U}}_N(\alpha_N, \widetilde{\rho}_N)} \mathbb{E}_{Q_{\widetilde{\Xi}}}[f(\mathbf{x}, \boldsymbol{\xi})]-\mathbb{E}_{\mathbb{Q}_{\widetilde{\Xi}}}[f(\mathbf{x}, \boldsymbol{\xi})] $$
  $$c-a_N=\mathbb{E}_{\mathbb{Q}_{\widetilde{\Xi}}}[f(\mathbf{x}, \boldsymbol{\xi})]-\min_{Q_{\widetilde{\Xi}}\in \widehat{\mathcal{U}}_N(\alpha_N, \widetilde{\rho}_N)} \mathbb{E}_{Q_{\widetilde{\Xi}}}[f(\mathbf{x}, \boldsymbol{\xi})] $$



 On the other hand, by  \citet[Proposition 2.1 (c)]{Hess1999} and the definition  of the Hausdorff distance, the following holds
$$\mathbb{H}( \mathcal{V},\text{conv} \mathcal{V}_N) \leqslant \mathbb{H}( \mathcal{V},\mathcal{V}_N)=\max (\mathbb{D}(\mathcal{V},\mathcal{V}_N), \mathbb{D}(\mathcal{V}_N,\mathcal{V})) = \mathbb{D}(\mathcal{V}_N,\mathcal{V})$$
where
\begin{align*}
 \mathbb{D}(\mathcal{V}_N,\mathcal{V})&=\max_{v'\in \mathcal{V}_N}d(v',\mathcal{V})=\max_{v'\in \mathcal{V}_N} \min_{v\in \mathcal{V}}|v'-v|\\
 &= \max_{Q_{\widetilde{\Xi}}\in \widehat{\mathcal{U}}_N(\alpha_N, \widetilde{\rho}_N)}
\left|\mathbb{E}_{Q_{\widetilde{\Xi}}}[f(\mathbf{x}, \boldsymbol{\xi})]-\mathbb{E}_{\mathbb{Q}_{\widetilde{\Xi}}}[f(\mathbf{x}, \boldsymbol{\xi})]\right|  \\
&\leqslant
\max_{Q_{\widetilde{\Xi}}\in \widehat{\mathcal{U}}_N(\alpha_N, \widetilde{\rho}_N)} \sup_{\mathbf{x}\in X}
\left|\mathbb{E}_{Q_{\widetilde{\Xi}}}[f(\mathbf{x}, \boldsymbol{\xi})]-\mathbb{E}_{\mathbb{Q}_{\widetilde{\Xi}}}[f(\mathbf{x}, \boldsymbol{\xi})]\right|\\
&=\max_{Q_{\widetilde{\Xi}}\in \widehat{\mathcal{U}}_N(\alpha_N, \widetilde{\rho}_N)} \mathcal{D}(Q_{\widetilde{\Xi}},\mathbb{Q}_{\widetilde{\Xi}})\\
&= \mathcal{D}(\widehat{\mathcal{U}}_N(\alpha_N, \widetilde{\rho}_N),\mathbb{Q}_{\widetilde{\Xi}})
\end{align*}

Note that $ \mathcal{D}( \widehat{\mathcal{U}}_N(\alpha_N, \widetilde{\rho}_N), \mathbb{Q}_{\widetilde{\Xi}}) \stackrel{N\rightarrow\infty}{\longrightarrow}0$, because, for any $f\in \mathcal{F}$, we have that
$\mathbb{E}_{Q_{\widetilde{\Xi}}}[f]\stackrel{N\rightarrow\infty}{\longrightarrow} \mathbb{E}_{\mathbb{Q}_{\widetilde{\Xi}}}[f]
$ under Corollary~\ref{Convergence_conditional} and Proposition~\ref{topology_wasserstein}.
Thus,
$$
 \mathbb{H}(\mathcal{V},\text{conv} \mathcal{V}_N) \leqslant \mathbb{H}( \mathcal{V},\mathcal{V}_N)= \mathbb{D}(\mathcal{V}_N,\mathcal{V})
  \leqslant \mathcal{D}(\widehat{\mathcal{U}}_N(\alpha_N, \widetilde{\rho}_N),\mathbb{Q}_{\widetilde{\Xi}})
$$
Therefore,
$$|v_N(\mathbf{x})-v(\mathbf{x})|\leqslant  \mathbb{H}(\mathcal{V},\text{conv} \mathcal{V}_N) \leqslant \mathcal{D}(\widehat{\mathcal{U}}_N(\alpha_N, \widetilde{\rho}_N),\mathbb{Q}_{\widetilde{\Xi}})\stackrel{N\rightarrow\infty}{\longrightarrow}0$$

Hence, since the inequality above is independent of the value
of $\mathbf{x}$, we have
$\lim_{N\rightarrow \infty}\sup_{\mathbf{x}\in X}|v_N(\mathbf{x})-v(\mathbf{x})|=0$ a.s.

Now, we show that the functions
$v_N(\mathbf{x})$ and $v(\mathbf{x})$ are continuous in $\mathbf{x}\in X$:
Fix an arbitrary $\mathbf{x}\in X$ and consider a sequence $(\mathbf{x}_N)_N$ such that $\mathbf{x}_N \rightarrow \mathbf{x}$ as $N$ grows to infinity. We want to prove that
$v_N(\mathbf{x}_N)\rightarrow v_N(\mathbf{x})$ and $v(\mathbf{x}_N)\rightarrow v(\mathbf{x})$. First, there exist $Q_{\mathbf{x}_N},Q_{\mathbf{x}} \in \widehat{\mathcal{U}}_N(\alpha_N, \widetilde{\rho}_N)$ such that $v_N(\mathbf{x}_N)=\mathbb{E}_{Q_{\mathbf{x}_N} }f(\mathbf{x}_N,\boldsymbol{\xi})$ and
$v_N(\mathbf{x})=\mathbb{E}_{Q_{\mathbf{x}} }f(\mathbf{x},\boldsymbol{\xi})$.
For any $\varepsilon>0$, there exists $N'>0$ sufficiently large such that for $N \geqslant N'$ the following holds:
\begin{align*}
 |v_N(\mathbf{x}_N)-v_N(\mathbf{x})|&= |\mathbb{E}_{Q_{\mathbf{x}_N} }f(\mathbf{x}_N,\boldsymbol{\xi})-\mathbb{E}_{Q_{\mathbf{x}} }f(\mathbf{x},\boldsymbol{\xi})|  \\
 &\leqslant |\mathbb{E}_{Q_{\mathbf{x}_N} }f(\mathbf{x}_N,\boldsymbol{\xi})-\mathbb{E}_{Q_{\mathbf{x}_N} }f(\mathbf{x},\boldsymbol{\xi})|+
|\mathbb{E}_{Q_{\mathbf{x}_N} }f(\mathbf{x},\boldsymbol{\xi})-\mathbb{E}_{Q_{\mathbf{x}} }f(\mathbf{x},\boldsymbol{\xi})|\\
&\leqslant \varepsilon/2+\varepsilon/2=\varepsilon
\end{align*}
since $|\mathbb{E}_{Q_{\mathbf{x}_N} }f(\mathbf{x}_N,\boldsymbol{\xi})-\mathbb{E}_{Q_{\mathbf{x}_N} }f(\mathbf{x},\boldsymbol{\xi})|<\varepsilon/2$ because $f$ is continuous in  $\mathbf{x}$ and $$|\mathbb{E}_{Q_{\mathbf{x}_N} }f(\mathbf{x},\boldsymbol{\xi})-\mathbb{E}_{Q_{\mathbf{x}} }f(\mathbf{x},\boldsymbol{\xi})|\leqslant\mathcal{D}(Q_{\mathbf{x}_N},Q_{\mathbf{x}})\leqslant
\mathcal{D}(Q_{\mathbf{x}_N},\mathbb{Q}_{\widetilde{\Xi}})+\mathcal{D}(\mathbb{Q}_{\widetilde{\Xi}},Q_{\mathbf{x}})\stackrel{N\rightarrow\infty}{\longrightarrow}0,$$ because $\mathcal{D}(\widehat{\mathcal{U}}_N(\alpha_N, \widetilde{\rho}_N),\mathbb{Q}_{\widetilde{\Xi}} )\stackrel{N\rightarrow\infty}{\longrightarrow}0$. As $\varepsilon>0$ is arbitrary, this implies that the function $v_N(\mathbf{x})$ is continuous in $\mathbf{x}\in X$. Similarly, since $f$ is continuous in  $\mathbf{x}$, we have that  the function $v(\mathbf{x})$ is continuous in $\mathbf{x}\in X$.
Finally, as $v_N(\mathbf{x})$ and $v(\mathbf{x})$ are continuous in $\mathbf{x}\in X$ and $\lim_{N\rightarrow \infty}\sup_{\mathbf{x}\in X}|v_N(\mathbf{x})-v(\mathbf{x})|=0$ a.s., we deduce from \citet[Lemma 3.4]{Xu2007} that $\widehat{J}_N\rightarrow J^*$ a.s. and the proof of (i) is complete.

The proof of (ii) is given by the application of  \citet[Lemma 3.8]{Liu2013}. \qed
\endproof

\begin{remark}
The theoretical framework underpinned  by Lemma~\ref{lemma_convergence_trimmed_distributions}, Corollary~\ref{Convergence_conditional}  and~\ref{theorem_consistency_trimmed} leaves the decision maker with considerable freedom to choose the values for $\alpha_N$ and $\widetilde{\rho}_N$. In the following two corollaries, we show that our framework naturally produces distributionally robust variants of popular non-parametric regression techniques such as the $K$-nearest neighbors and the Nadaraya-Watson kernel regression. This could serve to guide the selection of $\alpha_N$ and $\widetilde{\rho}_N$.

\begin{corollary}[Distributionally robust $K$-nearest neighbors]
   Let $K_N$ be the number of nearest neighbors, chosen such that  $K_N \rightarrow \infty,\; K_N/N \rightarrow 0$ and $ \frac{K_N}{\log N}\rightarrow \infty$
   when the sample size $N$ grows to infinity. This defines a standard KNN regression method.

   Take problem $\left({\text P}_{(\alpha_N,\widetilde{\rho}_N)}\right)$, set $\alpha_N:=K_N/N$ and compute the minimum transportation budget $\underline{\epsilon}_{K_N}$ as in Definition~\ref{MTB}. Problem $\left({\text P}_{(\alpha_N,\widetilde{\rho}_N)}\right)$ for any sequence of $\widetilde{\rho}_N$, $N \in \mathbb{N}$, such that $\widetilde{\rho}_N = \underline{\epsilon}_{K_N}^p + \Delta{\widetilde{\rho}_N}$ with $\Delta{\widetilde{\rho}_N}\downarrow 0$ is a distributionally robust variant of that KNN method.
    \end{corollary}
    \proof{Proof}  The proof of this claim directly follows from the fact that all the conditions in Lemma~\ref{lemma_convergence_trimmed_distributions} are satisfied if we choose $\alpha_N=K_N/N$. Actually, if we set $\widetilde{\rho}_N = \underline{\epsilon}_{K_N}^p$, the ambiguity set consisting of all distributions $Q^N_{\widetilde{\Xi}}$ such that $Q^N_{\widetilde{\Xi}}\in \widehat{\mathcal{U}}_N(\alpha_N, \widetilde{\rho}_N)$ is reduced, for each $N \in \mathbb{N},$ to the singleton $Q^N_{\widetilde{\Xi}}:=
\sum_{i=1}^{K_{N}} \frac{1}{K_N} \delta_{
(\mathbf{z}^*,\mathbf{\widehat{y}}_{i:N})}$, where $\mathbf{\widehat{y}}_{i:N}$ represents the $\mathbf{y}$-coordinate of the data point in the sample that is the $i$-th nearest neighbor. The decision maker can thus use the extra budget $\Delta\widetilde{\rho}_N$ to control the degree of robustness of the KNN solution.   \qed
\endproof

    \begin{corollary}[Distributionally robust Nadaraya-Watson kernel regression]
   Consider a Nadaraya-Watson (NW) kernel regression method with bandwidth $h_N$ such that  $h_N \rightarrow 0$  and $Nh_N^{d_{\mathbf{z}}}/\log(N)\rightarrow \infty$
   when $N$ grows to infinity. Also, assume that the non-negative Kernel $\mathcal{K}$ of the NW method satisfies that there exist positive numbers $c_1$, $c_2$ and $r$ such that
    $c_1 \mathbb{I}_{\{ \left\|\mathbf{v} \right\| \leqslant r  \}} \leqslant \mathcal{K}(\mathbf{v}) \leqslant c_2 \mathbb{I}_{\{ \left\|\mathbf{v} \right\| \leqslant r  \}} $.

    Let $w_i$, $i = 1, \ldots, N$ be the weights given by the NW method to the data points in a certain sample of size $N$ and let $w^{max}:= \max_{i}{w_i}$. Compute $$\widetilde{\rho}_N^{NW} = \sum_{i=1}^{N}w_i \, dist\left((\widehat{\mathbf{z}}_i,\widehat{\mathbf{y}}_i),\widetilde{\Xi}\right)^p.$$  The
    choices $\alpha_N:=1/(N w^{max})$ and $\widetilde{\rho}_N:= \widetilde{\rho}_N^{NW}+\Delta\widetilde{\rho}_N$ with $\Delta\widetilde{\rho}_N \downarrow 0$ produce an asymptotically consistent and distributionally robust Nadaraya-Watson kernel regression method.
    \end{corollary}
   \proof{Proof} To prove this corollary, we will use the following lemma, which appears in \cite{Devroye1981}.

\begin{lemma}[Lemma 4.1 from \cite{Devroye1981}]\label{lema_devroye}
If $n$ is a binomial random variable with parameters $N$ and $\hat{p}$, then
$$\sum_{N=1}^{\infty} \mathbb{E}\left[\exp{(-s\,n)}\right] < \infty, \enskip for\  all \ s>0$$
whenever $N\hat{p}/\log{N} \rightarrow \infty$.
\end{lemma}

    Define $A_i$ as the event $(\left\|  \widehat{\mathbf{z}}_i-\mathbf{z}^* \right\| \leqslant rh_N)$. Then, $n=\sum_{i=1}^{N}\mathbb{I}_{A_i}$ is a binomial random variable with parameters $N$ and $\hat{p} = \mathbb{P}(\left\|  \widehat{\mathbf{z}}_i-\mathbf{z}^* \right\| \leqslant rh_N)$ that represents the number of samples that are given a weight different from zero by the NW method. By Assumption~\ref{assumption_knn}, it follows that $\hat{p} \geqslant \widetilde{C} r^{d_{\mathbf{z}}} h_{N}^{d_{\mathbf{z}}}$, when $rh_N < r_0$. Furthermore, by the way the weights are constructed in this method and the choice of $\alpha_{N}$, we have that $\widetilde{\rho}_{N}^{NW} \geqslant \underline{\epsilon}^{p}_{N\alpha_N}$, provided that $n \geqslant 1$. In that case, it also holds $1/n \leqslant w^{max} \leqslant c_2/(c_1\, n)$ and thus, $(c_1 \,n)/c_2\leqslant N\alpha_N \leqslant n$. Note that the event $(n = 0)$ can happen only in a finite number of instances as $N$ increases. Indeed, for $N$ sufficiently large, $\mathbb{P}(n = 0) = \left(1-\hat{p}\right)^N \leqslant \exp{(-N\hat{p})} \leqslant \exp{\left(-N\widetilde{C} r^{d_{\mathbf{z}}} h_{N}^{d_{\mathbf{z}}}\right)}$, which is summable over $N$, because $Nh_N^{d_{\mathbf{z}}}/\log(N)\rightarrow \infty$. Therefore, in practice, the bandwidth of the NW method could be occasionally augmented in those specific instances so that $n \geq 1$, without affecting the convergence of the method.

    Thus, we have
    \begin{equation*}
        \widetilde{\rho}^{NW}_{N} \leqslant \frac{c_2}{c_1\, n}\sum_{i=1}^N \left\|  \widehat{\mathbf{z}}_i-\mathbf{z}^* \right\|^p \mathbb{I}_{A_i} \leqslant \frac{c_2 r^p h_N^p }{c_1} \rightarrow 0
    \end{equation*}
   because $h_N$ tends to $0$ as  $N$ grows to infinity.

   Now, we need to revisit Equation~\eqref{Conv_transp_trimmed_dist_1term}, since $N\alpha_N$ is  random here (contingent on the training sample). In particular, we have
   \begin{align*}
    \mathbb{P}\left( \left| \sum_{i=1}^N b^N_i h(\widehat{\mathbf{y}}_i)
    -\sum_{i=1}^N b^N_i m(\widehat{\mathbf{z}}_i)
    \right| > \varepsilon\  \mid \ \widehat{\mathbf{z}}_1, \ldots, \widehat{\mathbf{z}}_N  \right) &\leqslant 2\exp{
    \left(\dfrac{-(N\alpha_N) \varepsilon^2}{4\left\|h\right\|_{\infty}(2
    \left\|h\right\|_{\infty}+\varepsilon)}\right)}\\
    &\leqslant 2\exp{
    \left(\dfrac{-(c_1\,n/c_2) \varepsilon^2}{4\left\|h\right\|_{\infty}(2
    \left\|h\right\|_{\infty}+\varepsilon)}\right)}
\end{align*}
for any arbitrary $\varepsilon > 0$.

Hence,
\begin{equation*}
    \mathbb{P}\left( \left| \sum_{i=1}^N b^N_i h(\widehat{\mathbf{y}}_i)
    -\sum_{i=1}^N b^N_i m(\widehat{\mathbf{z}}_i)
    \right| > \varepsilon \right) \leqslant \mathbb{E}\left[ 2\exp{
    \left(\dfrac{-(c_1\,n/c_2) \varepsilon^2}{4\left\|h\right\|_{\infty}(2
    \left\|h\right\|_{\infty}+\varepsilon)}\right)}\right]
\end{equation*}
The summability with respect to $N$ of the expectation on the right-hand side of the inequality above is ensured by Lemma~\ref{lema_devroye}, given that, for $N$ large enough, $N\hat{p}/\log{N} \geqslant \widetilde{C} r^{d_{\mathbf{z}}}Nh_N^{d_{\mathbf{z}}}/\log(N)\rightarrow \infty$. The Borel-Cantelli lemma does the rest to conclude the proof.

While not explicitly required in this proof, it is easy to check that $\alpha_N \rightarrow 0$ almost surely as well. Note that $\frac{c_1\,n}{c_2\,N}\leqslant \alpha_N \leqslant \frac{n}{N}$, with $\mathbb{E}\left[\frac{n}{N}\right] = \hat{p} \rightarrow 0$, since $h_N \rightarrow 0$. Using \citet[Lemma 6]{devroye1982necessary}, we get, for any $\varepsilon >0$,
$$\mathbb{P}\left( \left|\frac{n}{N} - \hat{p}\right| > \varepsilon\right) = \mathbb{P}\left( \left|\sum_{i=1}^{N} \frac{1}{N} (\mathbb{I}_{A_i}-\hat{p})\right| > \varepsilon\right) \leqslant 2\exp{\left(-\frac{N\varepsilon^2}{2(1+\varepsilon)}\right)} $$
which is summable with respect to $N$. Thus, $\lim_{N\rightarrow \infty} \frac{n}{N} = \hat{p} = 0$ with probability one (as expected) and consequently, $\alpha_{N} \rightarrow 0$ a.s.  \qed
\endproof

Similarly as before, the extra budget $\Delta\widetilde{\rho}_N$ can be used by the decision-maker to robustify the NW solution. Nevertheless, in this case, as $\widetilde{\rho}_{N}^{NW} \geqslant \underline{\epsilon}_{N\alpha_N}^p$, the ambiguity set is not necessarily a singleton, meaning that our DRO approach already confers some degree of robustness on the decision vector $\mathbf{x}$ even if we set  $\widetilde{\rho}_{N} = \widetilde{\rho}_{N}^{NW}$.

\end{remark}

We conclude this section with a corollary that extends Lemma~\ref{lemma_convergence_trimmed_distributions}  to the case of unbounded uncertainty $\mathbf{y}$ under certain conditions. This extension guarantees that the solution to problem $\left({\rm P}_{(\alpha_N,\widetilde{\rho}_N)}\right)$ is asymptotically consistent also for this case.

\begin{corollary}[Extension of Lemma~\ref{lemma_convergence_trimmed_distributions} to unbounded $\mathbf{y}$]\label{Extension_unbounded_y}
Suppose that Assumptions~\ref{assumption_knn}.1 and \ref{assumption_lipschitz}  hold. Consider the true data-generating distribution $\mathbb{Q}$ of the random vector $\boldsymbol{\xi}:= (\mathbf{z},\mathbf{y})$ with support $\Xi:= \Xi_{\mathbf{z}}  \times \mathbb{R}^{d_{\mathbf{y}}}$ and define $m(\mathbf{z}^*) = \mathbb{E}[\|\mathbf{y}\|^a\mid\mathbf{z}=\mathbf{z}^*]$, for some $a \geqslant p$.

Assume that there exists a constant $\overline{m} >0$ such that $m(\mathbf{z}) <  \overline{m}$ for almost all $\mathbf{z} \in \Xi_{\mathbf{z}} $, and that there are non-negative numbers $(\sigma, \nu)$ such that
$$\log\mathbb{E}\left[ \exp\left\{t(\|   \mathbf{y}\|^{a}-m(\mathbf{z}))\right\}\mid\mathbf{z}=\mathbf{z^*}\right]\leqslant  \sigma^2 t^2/2,\; |t|\leq 1/\nu,$$  for almost all $\mathbf{z}^* \in \Xi_{\mathbf{z}} $.
Then, if the sequence $(\alpha_N, \widetilde{\rho}_N)$, $N \in \mathbb{N}$, meets the conditions stated in Lemma~\ref{lemma_convergence_trimmed_distributions}, we have that the convergence result stated in that lemma, also applies in the following two cases: i) $a=p$ and function $\ell:\mathbb{R}^{d_{\mathbf{y}}} \rightarrow \mathbb{R}_{+}$ in Assumption~\ref{assumption_lipschitz} is such that $\int{\|\mathbf{y}\|^p}\ell(\mathbf{y}) d\mathbf{y} < R < \infty$; and ii) $a>p$.
\end{corollary}
\proof{Proof} Since the weak convergence of distributions is guaranteed by way of Lemma~\ref{lemma_convergence_trimmed_distributions}, we just need to prove that $\int_{\widetilde{\Xi}}{\|\mathbf{y}\|^p dQ_{\widetilde{\Xi}}^N} \rightarrow \int_{\widetilde{\Xi}}{\|\mathbf{y}\|^p d\mathbb{Q}_{\widetilde{\Xi}}}$ (i.e., convergence of the $p$-th moment, see Proposition~\ref{topology_wasserstein}). For this purpose, we will use different strategies in cases i) and ii).

\paragraph{Case i):} Here we  follow a similar strategy  to that used to prove Lemma~\ref{lemma_convergence_trimmed_distributions}.
%


We have
\begin{equation*}
 \left| \sum_{i=1}^N b^N_i \|\widehat{\mathbf{y}}_i\|^p
    -m(\mathbf{z}^*)\right| \leqslant \left| \sum_{i=1}^N b^N_i \|\widehat{\mathbf{y}}_i\|^p
    -\sum_{i=1}^N b^N_i m(\widehat{\mathbf{z}}_i)
    \right|
    + \left|\sum_{i=1}^N b^N_i m(\widehat{\mathbf{z}}_i)
    -m(\mathbf{z}^*)
    \right|
\end{equation*}
To upper bound the first term on the right-hand side of the above inequality, we  exploit the subexponential character of  $\|\widehat{\mathbf{y}}_i\|^p$, $i=1, \ldots, N$ (understood as random variables). To this end, we  employ the following technical result, which corresponds to Theorem 2.51 in \cite{Bercu2015}.
\begin{theorem}[Theorem 2.51 from \cite{Bercu2015}]\label{Bercu}
Let $Z_1,\ldots, Z_n$ be a finite sequence of independent and centered random variables such that, for all $1\leqslant k\leqslant  n$, the random variable $Z_k$ satisfies $\log\mathbb{E}[\exp(tZ_k)] \leqslant l(t)$ for any $t\geqslant 0$, with $l(t)$ being a function from $[0, \infty)$ to $[0,\infty]$ with a concave derivative such that $l(0)=l'(0)=0$.

Denote $S_N = b_1Z_1 +\ldots+b_NZ_N$ for some positive real numbers $b_1,\ldots, b_N$. For any positive $\varepsilon$,
$$\mathbb{P}(S_N\geqslant \varepsilon) \leqslant \exp\left(-\frac{\|b\|_1^2}{\|b\|_2^2}l^*\left(\frac{\varepsilon}{\|b\|_1}\right)\right)$$
where $l^*$ stands for the  convex conjugate of $l$.
\end{theorem}

By assumption, we have
$$\log \mathbb{E}\left[\exp\{t(\|\mathbf{y}\|^p-m(\mathbf{z}^*))\}/\mathbf{z}=\mathbf{z}^*\right] \leqslant \frac{\sigma^2 t^2}{2} \text{ for } 0 \leqslant t \leqslant 1/\nu \text{ and for almost all } \mathbf{z}^* \in \Xi_{\mathbf{z}} $$

We take then $l(t):= \frac{\sigma^2 t^2}{2}$, if $0 \leqslant t \leqslant 1/\nu$, and  $l(t):= \infty$, if $t>1/\nu$. Therefore, $l^*(s)=\frac{s^2}{2\sigma^2}$, if $0<s\leqslant \sigma^2/\nu$
and $l^*(s)=\frac{s}{\nu}-\frac{\sigma^2}{2\nu^2}$, if $s>\sigma^2/\nu$.

Thus, for any arbitrary $\varepsilon >0$,
$$\mathbb{P}\left( \left| \sum_{i=1}^N b^N_i \|\widehat{\mathbf{y}}_i\|^p
    -\sum_{i=1}^N b^N_i m(\widehat{\mathbf{z}}_i)
    \right| \geqslant \varepsilon\, \mid\, \mathbf{\widehat{z}}_1, \ldots, \mathbf{\widehat{z}}_N \right) \leqslant 2\, \exp\left(-\frac{\|b\|_1^2}{\|b\|_2^2}l^*\left(\frac{\varepsilon}{\|b\|_1}\right)\right)$$
It holds $\|b\|_1 = 1$,  $\|b\|_2^2 \leqslant 1/N\alpha_N$, and $\frac{s}{\nu}-\frac{\sigma^2}{2\nu^2} > \frac{s}{2\nu}$, if $s>\sigma^2/\nu$. Hence,
\begin{align*}
\mathbb{P}\left( \left| \sum_{i=1}^N b^N_i \|\widehat{\mathbf{y}}_i\|^p
    -\sum_{i=1}^N b^N_i m(\widehat{\mathbf{z}}_i)
    \right| \geqslant \varepsilon\ \right)
    &\leqslant 2\, \exp\left(-\frac{N\alpha_N\varepsilon^2}{2\sigma^2}\right) \mathbb{I}_{(\varepsilon \leqslant \sigma^2/\nu)}\\
    &+2\,\exp\left(-\frac{N\alpha_N\varepsilon}{2\nu}\right) \mathbb{I}_{(\varepsilon >\sigma^2/\nu)}
\end{align*}
which is summable with respect to $N$ because $\frac{N\alpha_N}{\log(N)}\rightarrow \infty$.

To deal with the term $\left|\sum_{i=1}^N b^N_i m(\widehat{\mathbf{z}}_i)-m(\mathbf{z}^*)\right|$, we first note that
$$\left|\sum_{i=1}^N b^N_i m(\widehat{\mathbf{z}}_i)-m(\mathbf{z}^*)\right| \leqslant  \sum_{i=1}^N b^N_i \left| m(\widehat{\mathbf{z}}_i)-m(\mathbf{z}^*)\right|$$
where
\begin{align*}
\left| m(\widehat{\mathbf{z}}_i)-m(\mathbf{z}^*)\right| &= \left|\int \left\|\mathbf{y}\right\|^p \phi_{\mathbf{y}/\mathbf{z}=\mathbf{\widehat{z}}_i}(\mathbf{y}) d\mathbf{y}-\int \left\|\mathbf{y}\right\|^p \phi_{\mathbf{y}/\mathbf{z}=\mathbf{z}^*}(\mathbf{y}) d\mathbf{y}\right|\\
&\leqslant \int \left\|\mathbf{y}\right\|^p \left|(\phi_{\mathbf{y}/\mathbf{z}=\mathbf{\widehat{z}}_i}-\phi_{\mathbf{y}/\mathbf{z}=\mathbf{z}^*})(\mathbf{y})\right|d\mathbf{y}\\
&\leqslant \left\|\mathbf{\widehat{z}}_i-\mathbf{z}^*\right\| \int \left\|\mathbf{y}\right\|^p \ell(y) d\mathbf{y}\\
&\leqslant R\,\left\|\mathbf{\widehat{z}}_i-\mathbf{z}^*\right\|
\end{align*}
for any $\mathbf{\widehat{z}}_i$ such that $\left\|\mathbf{\widehat{z}}_i-\mathbf{z}^*\right\| \leqslant r_0$.

We finish the proof of case i) here, because, from this point on, the process is the same  as in Lemma~\ref{lemma_convergence_trimmed_distributions}, just replacing $L$ and $2\left\|h\right\|_{\infty}$ with $R$ and $\overline{m}$, respectively.

\paragraph{Case ii):} Based on the corollary to   \citet[Theorem 25.12]{Billingsley2}, it suffices to show that
$$\sup_N \int_{\mathbb{R}^{d_{\mathbf{y}}}} \|\mathbf{y}\|^a dQ^N_{\widetilde{\Xi}} < \infty$$

We first compute the integral for a fixed $N$.

\begin{align*}
\int_{\mathbb{R}^{d_{\mathbf{y}}}} \|\mathbf{y}\|^a dQ^N_{\widetilde{\Xi}} = \sum_{i=1}^{N}b_i^N \|\widehat{\mathbf{y}}_i\|^a &= \sum_{i=1}^{N}b_i^N \left( \|\widehat{\mathbf{y}}_i\|^a-m(\widehat{\mathbf{z}}_i)\right) + \sum_{i=1}^{N}b_i^N  m(\widehat{\mathbf{z}}_i)\\
&\leq \left|\sum_{i=1}^{N}b_i^N \left( \|\widehat{\mathbf{y}}_i\|^a-m(\widehat{\mathbf{z}}_i)\right)\right| + \overline{m}
\end{align*}

By Theorem~\ref{Bercu}, we have, for any arbitrary $\epsilon > 0$,
\begin{align*}
 \mathbb{P}\left( \left| \sum_{i=1}^N b^N_i \|\widehat{\mathbf{y}}_i\|^a
    -\sum_{i=1}^N b^N_i m(\widehat{\mathbf{z}}_i)
    \right| \geqslant \varepsilon\ \right)& \leqslant 2\,\exp\left(-\frac{N\alpha_N\varepsilon^2}{2\sigma^2}\right) \mathbb{I}_{(\varepsilon \leqslant \sigma^2/\nu)}\\&+2\,\exp\left(-\frac{N\alpha_N\varepsilon}{2\nu}\right) \mathbb{I}_{(\varepsilon >\sigma^2/\nu)}
\end{align*}
which is summable with respect to $N$,  because $\frac{N\alpha_N}{\log(N)}\rightarrow \infty$. Take $\varepsilon:= \varepsilon_0 >0$, there must then exist a sufficiently large $N_0$ such that
$$\left|\sum_{i=1}^{N}b_i^N \left( \|\widehat{\mathbf{y}}_i\|^a-m(\widehat{\mathbf{z}}_i)\right)\right| < \varepsilon_0$$
for $N \geqslant N_0$ with probability one.

Therefore,
$$\int_{\mathbb{R}^{d_{\mathbf{y}}}} \|\mathbf{y}\|^a dQ^N_{\widetilde{\Xi}} \leqslant \varepsilon_0 + \overline{m}$$ for large enough $N \geqslant N_0$.

Thus,
$$\sup_N \int_{\mathbb{R}^{d_{\mathbf{y}}}} \|y\|^a dQ^N_{\widetilde{\Xi}} \leqslant \text{max}\left\{\sup_{N < N_0} \int_{\mathbb{R}^{d_{\mathbf{y}}}} \|y\|^a dQ^N_{\widetilde{\Xi}},\  \varepsilon_0 + \overline{m}\right\} <\infty \text{ a.s. } $$
\qed
\endproof

\begin{remark}
The proof of Corollary~\ref{Extension_unbounded_y} is considerably simplified if it holds $$\left| m(\mathbf{z})-m(\mathbf{z}^*)\right|  \leqslant R\,\left\|\mathbf{z}-\mathbf{z}^*\right\|$$ for almost all $\mathbf{z} \in \Xi_{\mathbf{z}} $ and some $R >0$. In this case, for instance, we do not need the almost-everywhere boundedness condition on random variable $m(\mathbf{z})$.

\end{remark}

\section{Additional numerical experiments, case $\alpha>0$. Portfolio optimization}\label{NE_alpha_pos}
In this section, we present and discuss some numerical results for the case $\mathbb{Q}(\widetilde{\Xi}) = \alpha>0$. For this purpose, we use the same portfolio allocation problem described in the main manuscript. To this end, we assume instead that the feature vector lives in an uncertainty set $\mathcal{Z}$ such that $\mathbb{Q}(\widetilde{\Xi})>0$. In particular,  we consider $\mathcal{Z}:=\{ \mathbf{z}\in \mathbb{R}^3 :\|\mathbf{\tilde{z}}\|_{\infty}\leqslant r\}$, with $\mathbf{\tilde{z}}$ being the standardized feature vector.
Thus, we have that $\widetilde{\Xi}$ is given by
 $$\widetilde{\Xi}:=\{ (\mathbf{z},\mathbf{y})\in \mathbb{R}^{3+6}  :  \|\mathbf{\tilde{z}}\|_{\infty}\leqslant r\}$$
We take $r = 0.6$ for the simulation experiments.

We draw 50\,000 samples from the true joint data-generating distribution through the explicit form of $\mathbf{y}/\mathbf{z}$ given in the main text. We then use the conditional empirical distribution made up of those samples falling within $\widetilde{\Xi}$, specifically, 7306 data points, as a proxy of the true conditional distribution $\mathbb{Q}_{\widetilde{\Xi}}$. Consequently, we have that $\mathbb{Q}(\widetilde{\Xi}) \approx 0.14612$. We wish to solve the following optimization problem
\begin{equation}\label{Portf-alpha-pos}
    \min_{(\mathbf{x},\beta') \in X }  \mathbb{E}\left[\beta'+\frac{1}{\delta}\left(-\langle \mathbf{x},\mathbf{y} \rangle -\beta' \right)^+-\lambda\langle \mathbf{x},\mathbf{y} \rangle
    \; \mid \; (\mathbf{z},\mathbf{y})\in \widetilde{\Xi} \right]
  \end{equation}
 with the rest of the parameters being equal to the values taken in the instance $\alpha = 0$.

 We also compare here four data-driven approaches to solve problem~\eqref{Portf-alpha-pos}, namely:
 \begin{itemize}
     \item Our two approaches, i.e., problem $\text{P}_{(\alpha,\widetilde{\rho}_N)}$ with $\alpha:=\mathbb{Q}(\widetilde{\Xi})$ (that is approximately equal to $0.14612$, as we have just mentioned), denoted as ``DROTRIMM1'' and problem $\text{P}_{(\alpha_N,\widetilde{\rho}_N)}$, where
$\alpha_N:= \widehat{\mathbb{Q}}_N(\widetilde{\Xi})$ is an estimate  of $\alpha$. We refer to this approach as ```DROTRIMM2.'' In principle, this would be the natural approach that a decision-maker with no knowledge of $\alpha$ would use.
\item A sample average approximation (SAA) method that works with the samples falling in $\widetilde{\Xi}$.
\item The aforementioned SAA method followed by a standard Wasserstein-metric-based DRO approach to robustify it,
which we call ``SAADRO''.
 \end{itemize}

As in the previous numerical experiments, we employ a similar bootstrapping procedure based on the available data sample to tune the robustness parameter that each method $j$, with $j \in$ $\{\textrm{DROMTRIMM1, DROTRIMM2, SAADRO}\}$, uses.  More specifically, for each
$j \in \{\textrm{DROMTRIMM1, DROTRIMM2, SAADRO}\}$ and
a given value of reliability $1-\beta \in (0,1)$ (in our numerical experiments, we set $\beta$ to 0.15),  we seek an estimator  $param^{\beta,j}_N$ that leads to the best out-of-sample performance, while guaranteeing the desired level of confidence $1-\beta$. For each sample of size $N$, we use the following algorithm
to derive $param^{\beta,j}_N$ and the corresponding portfolio solution:
  \begin{enumerate}
    \item  We construct $kboot$ resamples (with replacement) of size $N$, each playing the role of a different training dataset. In our experiments we use $kboot = 50$.
    Moreover, we build a validation dataset (per resample) from those data points from the original sample of size $N$ that fall in $\widetilde{\Xi}$, but which have not been involved in the resample. We only consider resamples from which we can build a validation set of at least one data point. Furthermore, unlike DROTRIMM1 and DROTRIMM2,  SAADRO can only be implemented  if we have at least one data point falling
    within $\widetilde{\Xi}$ in the training set (the same occurs with SAA). Thus,
    we implicitly assume that the source sample has no fewer than two data points in $\widetilde{\Xi}$.

    \item For each resample $k=1,\ldots,kboot$ and each candidate value for $param$ (taken from the discrete set $\{b\cdot  10^c\;:\;  b \in \{0,\ldots, 9 \},\; c \in \{-3,-2,-1,0\}\}$),  we compute a  solution by method $j$ with parameter  $param$ on the $k$-th resample.
    The resulting optimal decision  is denoted as $\widehat{x}^{j,k}_N(param)$ and its corresponding objective value as $\widehat{J}^{j,k}_N(param)$. Thereafter, we calculate the out-of-sample performance $J(\widehat{x}^{j,k}_N(param))$  of the data-driven solution $\widehat{x}^{j,k}_N(param)$ over the  validation dataset.
\item From among the candidate values for $param$ such that $\widehat{J}^{j,k}_N(param)$
exceeds the value $J(\widehat{x}^{j,k}_N(param))$
in at least
$(1-\beta)\times kboot$ different resamples, we take as $param^{\beta,j}_N$  the one yielding the best cost performance averaged over the \emph{kboot} resamples.
\item Finally, we compute the solution given by method $j$ with parameter
$param^{\beta,j}_N$, $\widehat{x}^j_N:=\widehat{x}^{j}_N(param^{\beta,j}_N)$
and the respective certificate
$\widehat{J}^j_N:=\widehat{J}^{j}_N(param^{\beta,j}_N)$.
\end{enumerate}

\begin{figure}[htp]

\centering

\subfloat{%
  \includegraphics[width=0.5\textwidth]{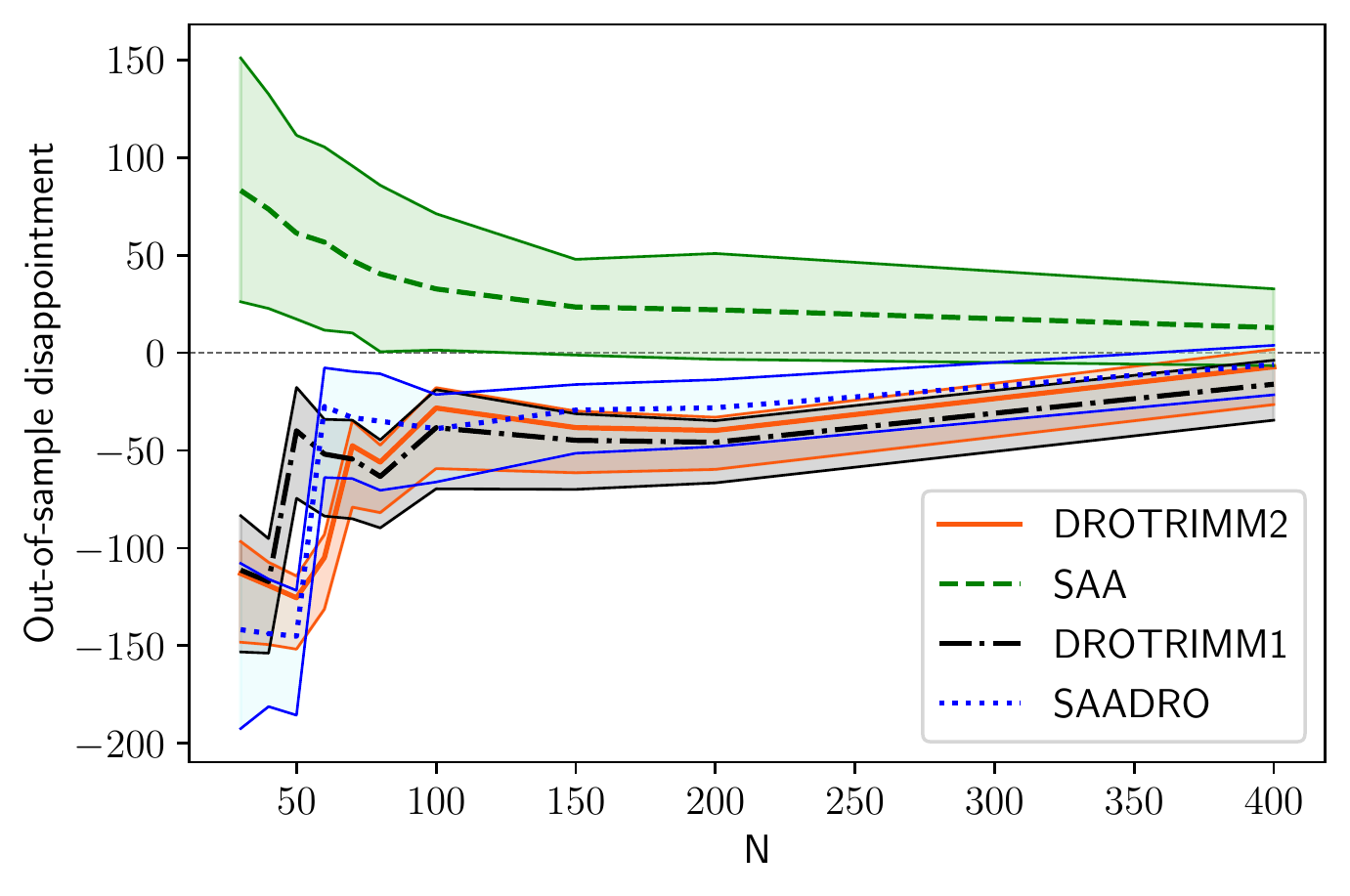}%
  \label{out-of-sample_portfolio_alpha_pos_radio_variante2.pdf}
}%
\subfloat{%
  \includegraphics[width=0.5\textwidth]{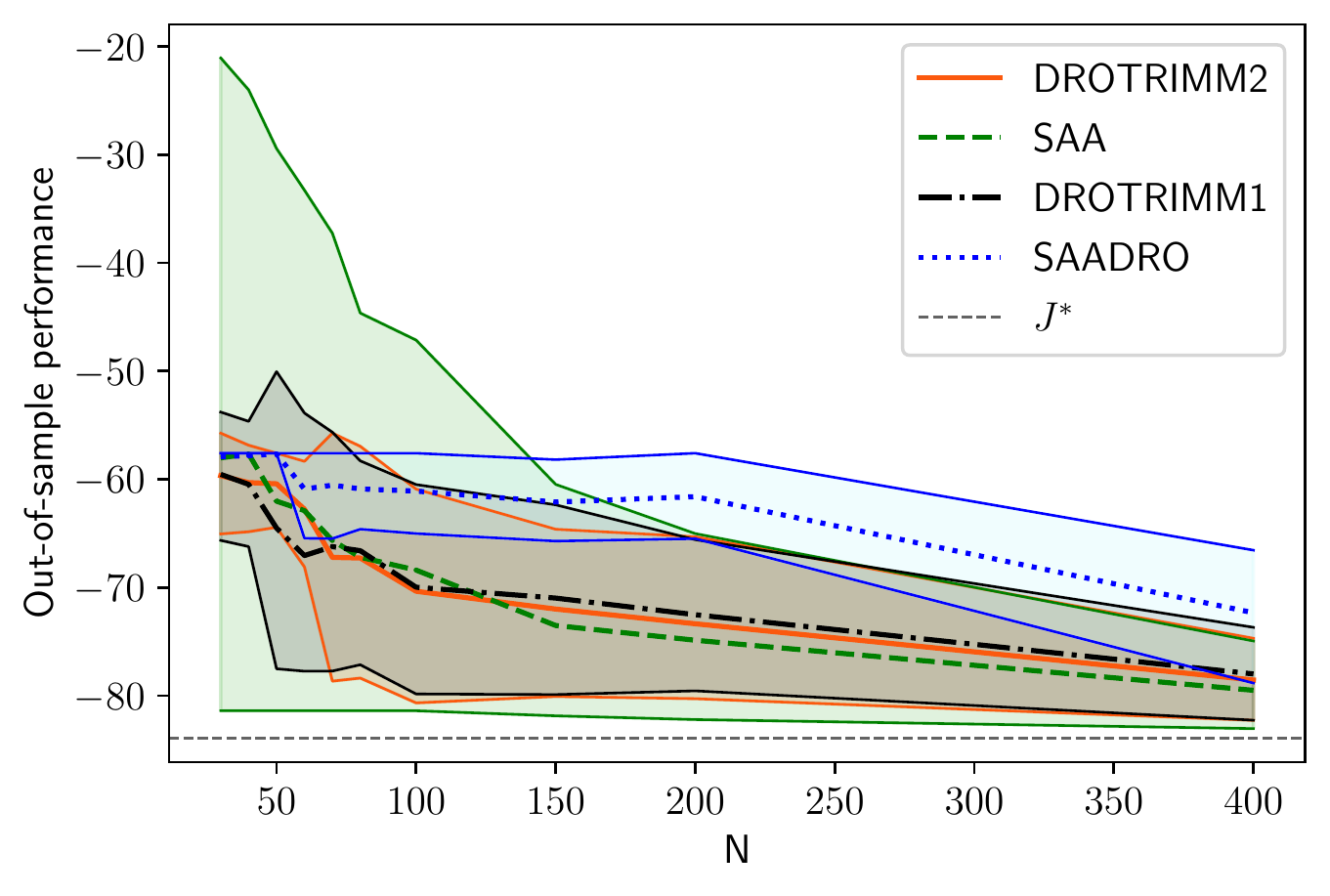}%
  \label{actual_expected_cost_portfolio_alphapos_radio_variante2.pdf}
}

\vspace{5mm}

\caption{ Portfolio problem  with features:  Performance metrics. Case $\alpha>0$ and $\delta=0.5,\ \lambda=0.1$}\label{Results_alpha_pos_performance}

\end{figure}

 Figure~\ref{Results_alpha_pos_performance} shows the box plots pertaining to the out-of-sample disappointment and performance associated with each of the considered data-driven approaches for various sample sizes. The box plots have been obtained from 200 independent runs per sample size $N$. The SAA method provides portfolios that, in expectation, perform reasonably well, especially when the sample size is large enough. However, SAA definitely fails to ensure the  desired level of reliability. As for the three approaches that incorporate robustness in the decision-making, DROTRIMM1 and DROTRIMM2 seem to systematically identify reliable portfolios with a better expected performance than those given by SAADRO.


\begin{figure}[htp]
\centering
\subfloat{%
  \includegraphics[width=0.5\textwidth]{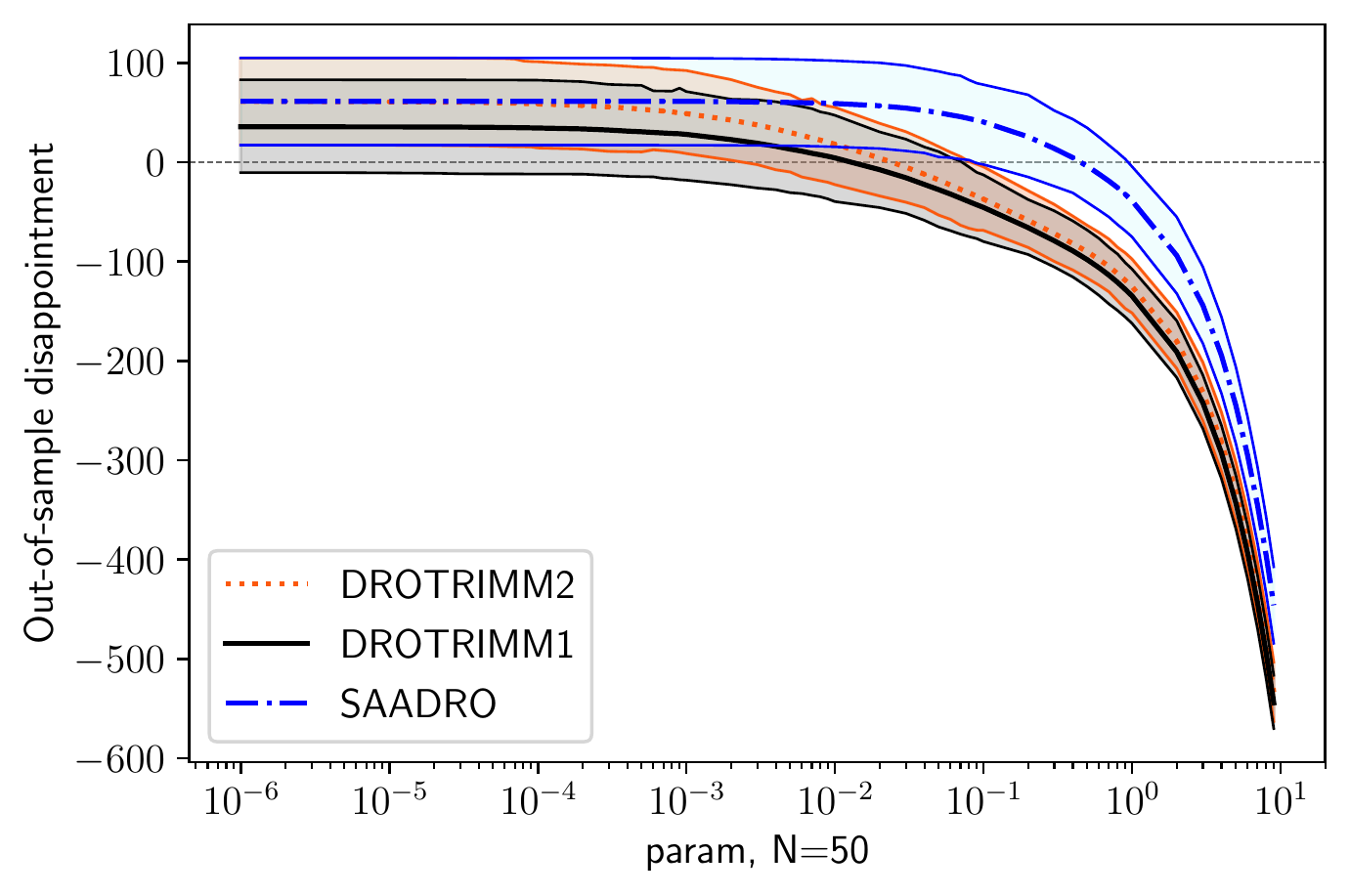}%
  \label{out-of-sample_portfolio_sens_50_DRO_variante1pos.pdf}
}%
\subfloat{%
  \includegraphics[width=0.5\textwidth]{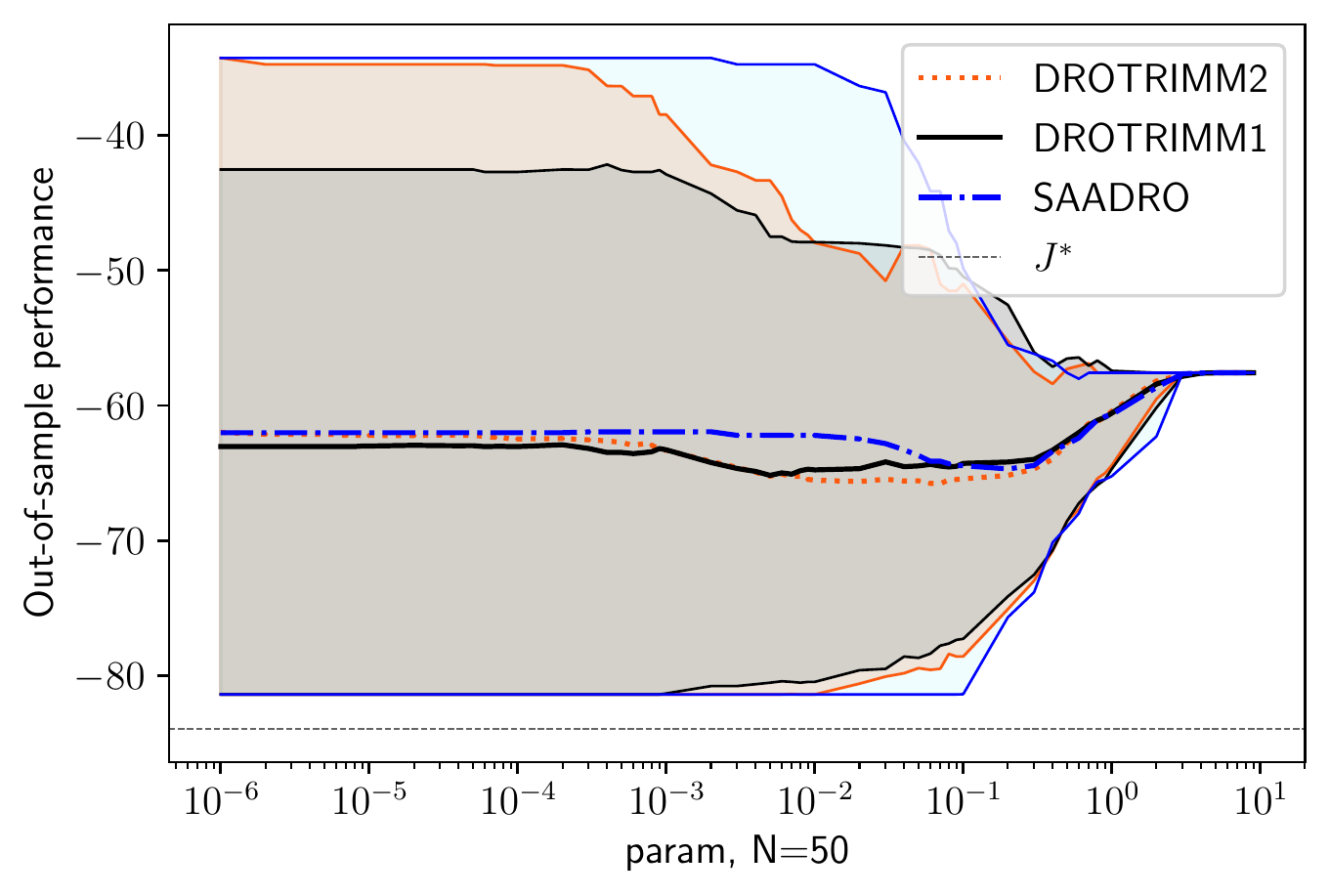}%
  \label{actualexpectedcost_portfolio_sens_50_DROvariante1pos.pdf}
}

\vspace{5mm}
\centering
\subfloat{%
  \includegraphics[width=0.5\textwidth]{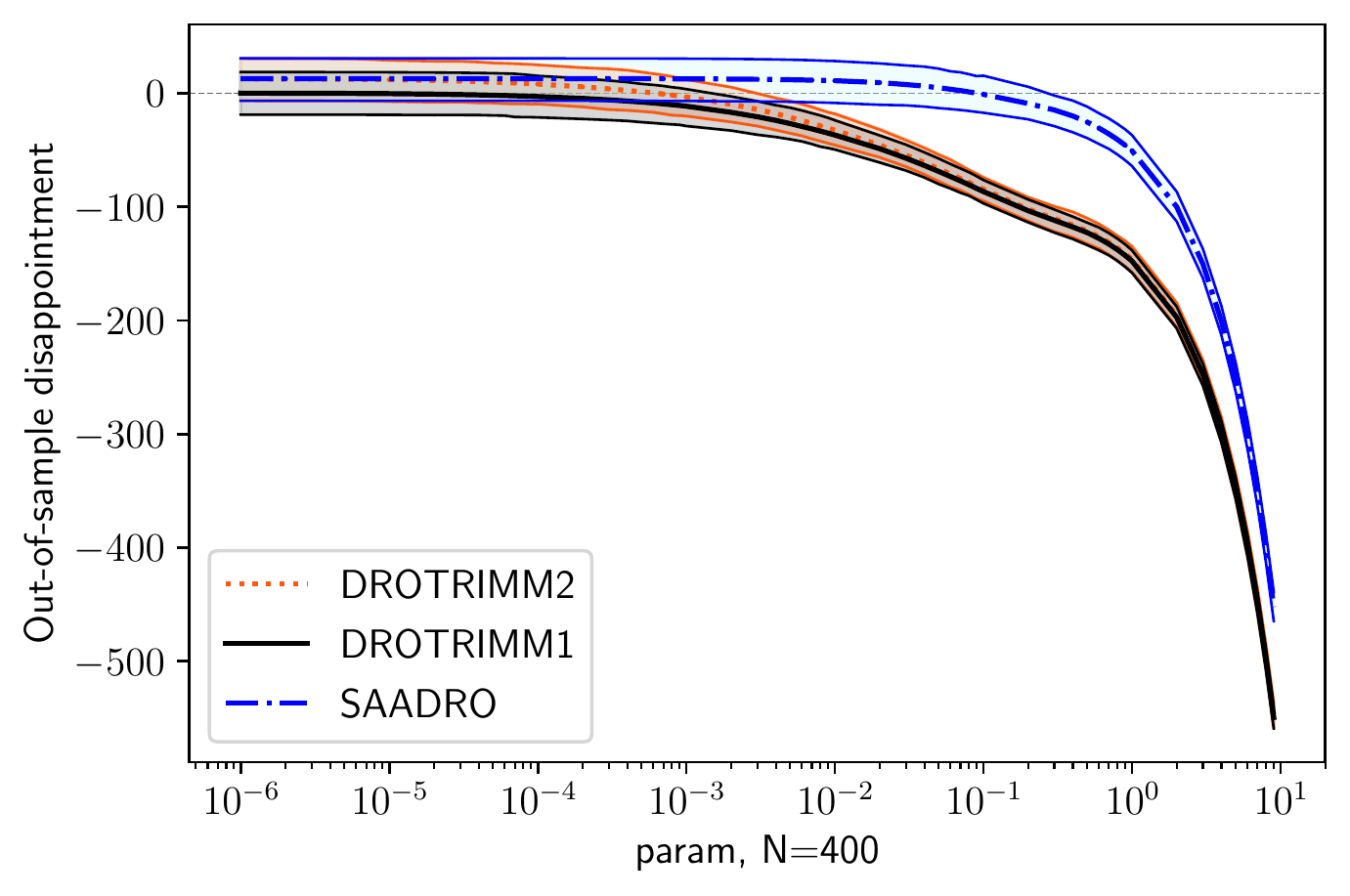}%
  \label{out-of-sample_portfolio_sens_400_DRO_variante1pos.pdf}
}%
\subfloat{%
  \includegraphics[width=0.5\textwidth]{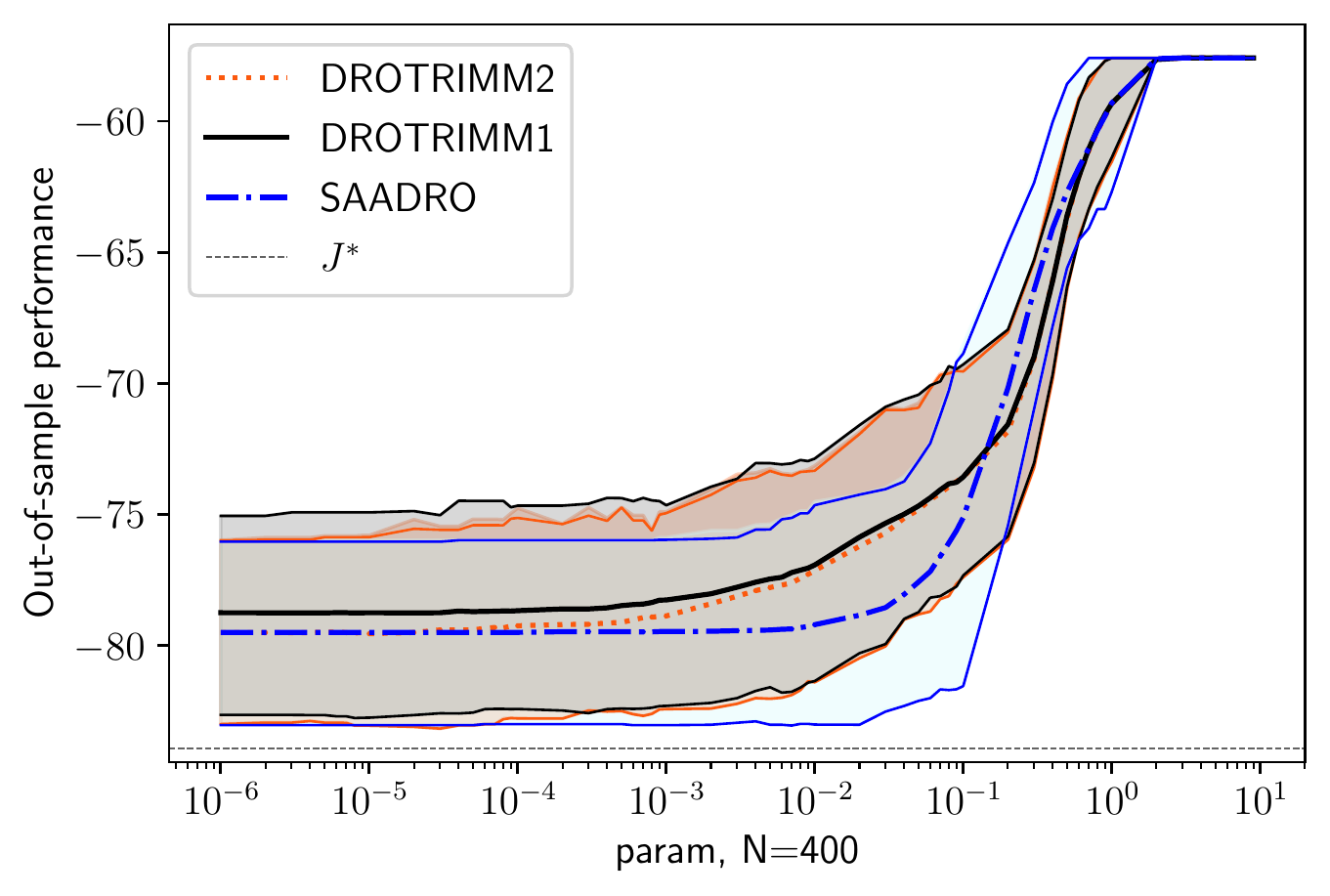}%
  \label{actualexpectedcost_portfolio_sens_400_DROvariante1pos.pdf}
}

\caption{Case $\alpha>0$, impact of the robustness parameter with 200 training samples  and $\delta=0.5,\lambda=0.1$}\label{alpha_pos_sens}
\end{figure}

 To investigate the ability of SAADRO, DROTRIMM1 and DROTRIMM2 to identify good portfolios, we provide Figure~\ref{alpha_pos_sens}, which is analogous to Figure~\ref{alpha_0_sens} in the case $\alpha = 0$. Observe that both DROTRIMM1 and DROTRIMM2 guarantee reliability for smaller values of their robustness parameter than SAADRO. This gives the former a competitive advantage over the latter, essentially because it appears that a better out-of-sample performance (in expectation) is, in general, aligned with a lower distributional robustness (this finding is consistent with the fact that the unreliable SAA solution performs fairly well in terms of the weighted mean-risk asset returns). To be more precise, taking a small sample size $N$ (say 50) and an equal value for each of their robustness parameters, DROTRIMM1 and DROTRIMM2 deliver portfolios with an actual expected cost (and variance) that is lower than or approximately equal to that of the portfolios provided by SAADRO. They do so for any value of their robustness parameter. Furthermore, when $N$ is increased, even though there exists a range of values of the robustness parameter for which SAADRO also identifies portfolios with a good performance out of sample, these are discarded by the method because they do not comply with the reliability specification. For instance, take $N = 400$. SAADRO needs a radius larger than $0.2$-$0.3$ to ensure reliability. However, for these values of the Wasserstein-ball radius, the portfolios given by SAADRO result in an actual expected cost above $-70$. On the other hand, DROTRIMM2 guarantees reliability with a value of its robustness parameter above $0.003$-$0.004$, for which, in addition, it provides solutions with an actual expected cost below $-77$.

To further support this finding, we conclude this section with Figure~\ref{Results_alpha_pos_performance_perfectValidation}, which is similar to Figure~\ref{Results_alpha_pos_performance}. However, Figure~\ref{Results_alpha_pos_performance_perfectValidation} has been obtained through a different experiment, in which the value of the robustness parameter that each method uses has been \emph{optimally} selected from the previously indicated discrete set. In other words, the results shown in that figure are those a decision-maker would obtain in the hypothetical case that  the true conditional distribution $\mathbb{Q}_{\widetilde{\Xi}}$ could be used to tune the robustness parameters of the DRO methods. Therefore, these results correspond to the best solutions that can be obtained from SAADRO, DROTRIMM1 and DROTRIMM2, and confirm that our approaches (especially, DROTRIMM2) can potentially identify portfolios that significantly outperform those delivered by SAADRO under the same reliability requirement.

\begin{figure}[htp]

\centering

\subfloat{%
  \includegraphics[width=0.5\textwidth]{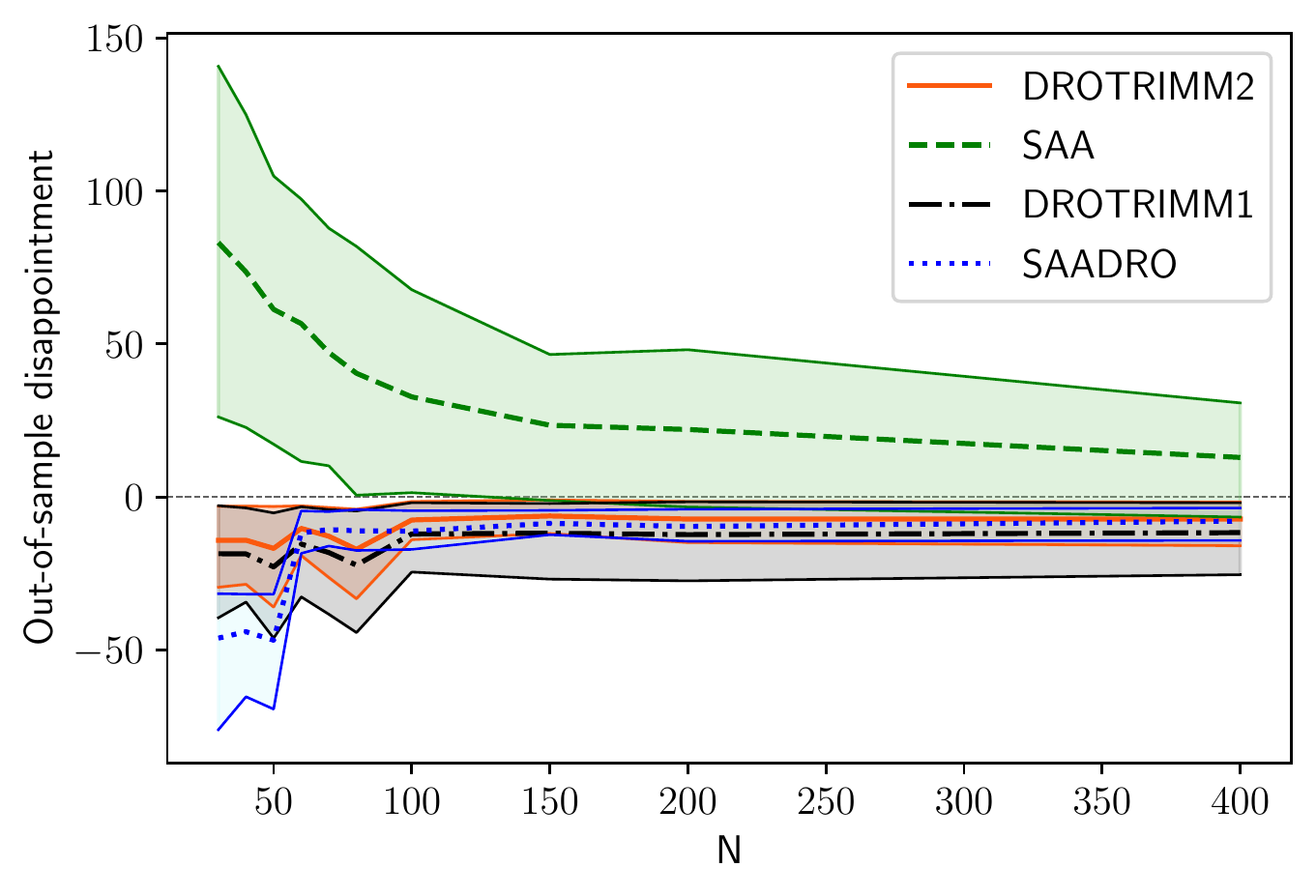}%
  \label{out_ofsample_portfolio_validperf_fijo.pdf}
}%
\subfloat{%
  \includegraphics[width=0.5\textwidth]{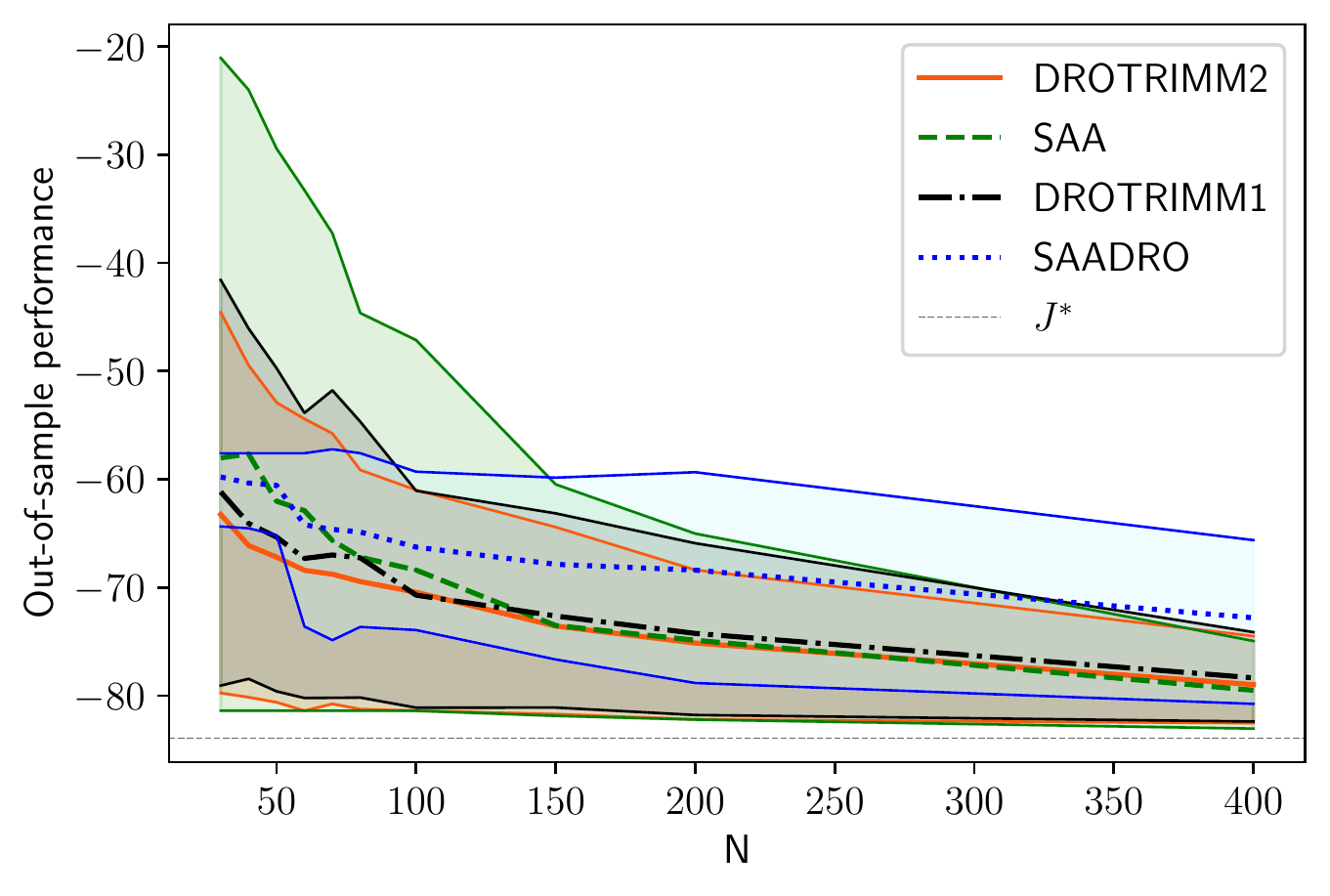}%
  \label{actualexpected_portfolio_validperf_fijo.pdf}
}

\vspace{5mm}

\caption{Portfolio problem  with features:  Performance metrics under an optimal selection of the robustness parameters. Case $\alpha>0$ and $\delta=0.5,\ \lambda=0.1$}\label{Results_alpha_pos_performance_perfectValidation}

\end{figure}

\end{document}